\documentclass[lettersize,journal]{IEEEtran}
\usepackage{amsmath,amsfonts}
\usepackage{algorithmic}
\usepackage{algorithm}
\usepackage{array}
\usepackage[caption=false,font=normalsize,labelfont=sf,textfont=sf]{subfig}
\usepackage{textcomp}
\usepackage{stfloats}
\usepackage{url}
\usepackage{verbatim}
\usepackage{graphicx}
\usepackage{ntheorem}
\theoremseparator{.}
\usepackage{cite}
\usepackage{multirow}
\usepackage{float}

\newtheorem{assumption}{Assumption}
\newtheorem{proposition}{Proposition}
\newtheorem{theorem}{Theorem}

\newtheorem{definition}{Definition}
\newtheorem{remark}{Remark}

\begin{document}

\title{A New Simple Stochastic Gradient Descent Type Algorithm With Lower Computational Complexity for Bilevel Optimization}

\author{Haimei Huo, Risheng Liu, \IEEEmembership{Member, IEEE}, Zhixun Su, \IEEEmembership{Member, IEEE}
\thanks{Haimei Huo is with the School of Mathematical Sciences, Dalian University of Technology, Dalian 116024, China (e-mail: ab1234@mail.dlut.edu.cn).} 
\thanks{Risheng Liu(Corresponding author) is with the DUT-RU International School of Information Science and Engineering, Dalian University of Technology, and also with the Key Laboratory for Ubiquitous Network and Service Software of Liaoning Province, Dalian 116024, China (e-mail: rsliu@dlut.edu.cn).}
\thanks{Zhixun Su(Corresponding author) is with the School of Mathematical Sciences, Dalian University of Technology, with the State Key Laboratory of Structural Analysis for Industrial Equipment, and also with the Key Laboratory for Computational Mathematics and Data Intelligence of Liaoning Province, Dalian 116024, China (e-mail: zxsu@dlut.edu.cn).}
}

\markboth{Journal of \LaTeX\ Class Files,~Vol.~14, No.~8, August~2021}
{Huo \MakeLowercase{\textit{et al.}}: A New Simple SGD Type Algorithm With Lower Computational Complexity for Bilevel Optimization}


\maketitle

\begin{abstract}
Bilevel optimization has been widely used in many machine learning applications such as hyperparameter optimization and meta learning. Recently, many simple stochastic gradient descent(SGD) type algorithms(without using momentum and variance techniques) have been proposed to solve the bilevel optimization problems. However, all the existing simple SGD type algorithms estimate the hypergradient via stochastic estimation of Neumann series. In the paper, we propose to estimate the hypergradient via SGD-based Estimation(i.e., solving the linear system with SGD). By using warm start initialization strategy, a new simple SGD type algorithm SSGD based on SGD-based Estimation is proposed. We provide the convergence rate guarantee for SSGD and show that SSGD outperforms the best known computational complexity achieved by the existing simple SGD type algorithms. Our experiments validate our theoretical results and demonstrate the efficiency of our proposed algorithm SSGD in hyperparameter optimization applications.
\end{abstract}

\begin{IEEEkeywords}
Bilevel optimization, stochastic gradient descent(SGD)-based Estimation, simple SGD type algorithms, warm start.
\end{IEEEkeywords}

\section{Introduction}
\IEEEPARstart{B}{ilevel} optimization has become a powerful tool for many machine learning applications such as meta learning\cite{1,2, 3, 4}, hyperparameter optimization\cite{5}, and reinforcement learning\cite{6, 7}. It has a nested structure which contains two levels of optimization tasks, and the feasible region of the upper level(UL) function is restricted by the optimal solutions of the lower level(LL) problem. Generally, bilevel optimization can be expressed as follows:
\begin{align}   \label{eq1}
&\min\limits_{\boldsymbol{x}\in \mathbb{R}^p} \Phi(\boldsymbol{x}) := f(\boldsymbol{x}, \boldsymbol{y}^*(\boldsymbol{x}))\\
&s.t. ~\boldsymbol{y}^*(\boldsymbol{x}) \in \mathcal{S}(\boldsymbol{x}) := { \underset {\boldsymbol{y}\in \mathbb{R}^q} { \operatorname {arg\,min} } \, g(\boldsymbol{x}, \boldsymbol{y})} \nonumber
\end{align}
where $f(\boldsymbol{x}, \boldsymbol{y})$ is the UL function, $g(\boldsymbol{x}, \boldsymbol{y})$ is the LL function, and $\mathcal{S}(\boldsymbol{x})$ is the solution set of the LL problem.

Many deterministic algorithms have been proposed to solve the bilevel problem (\ref{eq1}). In early works, the bilevel problem is reformulated into single level problem using the first order optimality conditions of the LL problem\cite{8, 9, 10, 29}. However, such type of methods often involve a large number of constraints, and are not applicable in complex machine learning problems. To address the problem, gradient based algorithms have been revisited to solve problem (\ref{eq1}). Depending on the way to approximate the hypergradient $\nabla \Phi(\boldsymbol{x})$, these algorithms can be generally categorized into the iterative differentiation(ITD) based approach\cite{1, 5, 11, 12, 13} and the approximate implicit differentiation(AID) based approach\cite{11, 14, 15, 16}.

ITD based approach treats the iterative optimization algorithm for the LL problem as a dynamical system, and approximates the hypergradient by back-propagation(BP) along the dynamical system \cite{11, 17}, this kind of hypergradient approximation method termed as BP in the paper. Since computing the gradient along the entire dynamical system is computationally challenging for high-dimensional problems, \cite{5} proposed to truncate the path to compute the gradient. On the other hand, the convergence of ITD based approach for solving problem $(\ref{eq1})$ has also been studied. \cite{1} assumed that the LL function $g(\boldsymbol{x},\cdot)$ is strongly convex, and provided the asymptotic convergence analysis. To address the case that the LL function $g(\boldsymbol{x}, \cdot)$ is convex, \cite{18} proposed BDA and the asymptotic convergence analysis is provided. For the nonasymptotic convergence analysis of ITD based approach, see \cite{13}.

AID based approach applies implicit function theorem to the first order optimality conditions of the LL problem, and approximates the hypergradient by solving a linear system\cite{14, 15}. There are mainly two methods to solve the linear system: conjugate gradient descent(CG)\cite{11} and approximating the Hessian inverse with Neumann series\cite{16}, these two hypergradient approximation methods termed as CG and Neumann series(NS) in the paper. For the nonasymptotic convergence analysis of the AID based approach, see \cite{13, 16}.

In fact, in many machine learning applications, the bilevel optimization takes the following form:
\begin{align}  \label{eq2}
&\min\limits_{\boldsymbol{x}\in \mathbb{R}^p} \Phi(\boldsymbol{x}) := f\left( \boldsymbol{x}, \boldsymbol{y}^*(\boldsymbol{x})\right) := \mathbb{E}_{\boldsymbol{\xi}}\left[F(\boldsymbol{x}, \boldsymbol{y}^*(\boldsymbol{x}); \boldsymbol{\xi})\right] \\
&s.t. ~ \boldsymbol{y}^*(\boldsymbol{x}) :=  { \underset {\boldsymbol{y}\in \mathbb{R}^q} { \operatorname {arg\,min} } \, g(\boldsymbol{x}, \boldsymbol{y})} := \mathbb{E}_{\boldsymbol{\zeta}}[G(\boldsymbol{x}, \boldsymbol{y}; \boldsymbol{\zeta)}] \nonumber
\end{align}
where UL function $f(\boldsymbol{x}, \boldsymbol{y})$ and LL function $g(\boldsymbol{x}, \boldsymbol{y})$ take the expected values with respect to(w.r.t.) random variables $\boldsymbol{\xi}$ and $\boldsymbol{\zeta}$, respectively; $g(\boldsymbol{x}, \cdot)$ is strongly convex w.r.t. $\boldsymbol{y}$ and $\Phi(\boldsymbol{x})$ is nonconvex w.r.t. $\boldsymbol{x}$. Notice that if $\boldsymbol{\xi}$ and $\boldsymbol{\zeta}$ are discrete random variables, the values of $\boldsymbol{\xi}$ and $\boldsymbol{\zeta}$ are chosen from the given data sets $\mathcal{D}_n := \{\boldsymbol{\xi}_1, \ldots, \boldsymbol{\xi}_n\}$ and $\mathcal{D}_m := \{\boldsymbol{\zeta}_1, \ldots, \boldsymbol{\zeta}_m\}$, respectively, and the probability to sample data from $\mathcal{D}_n$ and $\mathcal{D}_m$ is equal, respectively, then the objective functions in problem (\ref{eq2}) can be written as $f(\boldsymbol{x}, \boldsymbol{y}):= \frac{1}{n}\sum_{i=1}^{n}F(\boldsymbol{x}, \boldsymbol{y}; \boldsymbol{\xi}_i)$ and $g(\boldsymbol{x}, \boldsymbol{y}):= \frac{1}{m}\sum_{i=1}^{m}G(\boldsymbol{x},\boldsymbol{y}; \boldsymbol{\zeta}_i)$.

In order to achieve better efficiency than deterministic algorithms for solving problem (\ref{eq2}), many stochastic algorithms have been proposed. Without considering the loop to estimate the Hessian inverse, these existing stochastic algorithms for problem (\ref{eq2}) can be categorized into single loop algorithms and double loop algorithms according to the number of loops.

For single loop algorithms, the upper variable $\boldsymbol{x}$ and the LL variable $\boldsymbol{y}$ are updated simultaneously. As far as we know, the first single loop algorithm proposed for solving problem (\ref{eq2}) is the two-timescale stochastic approximation(TTSA) algorithm\cite{7}, which is a simple SGD type algorithm. To ensure convergence, a larger step size for the LL problem and a smaller step size for the UL problem are required. Moreover, to reach an $\epsilon$-stationary point, the computational complexity is $\widetilde{\mathcal{O}}(\epsilon^{-5/2})$ in terms of the target accuracy $\epsilon$. Recently, based on momentum and variance techniques, some single loop algorithms with improved computational complexity have been proposed, see e.g., \cite{19, 20}.

For a double loop algorithm, it contains two loops in which the inner loop is used to update the LL variable $\boldsymbol{y}$ to obtain an approximation to the optimal solution $\boldsymbol{y}^*(\boldsymbol{x})$ of the LL problem, and the outer loop is used to update the UL variable $\boldsymbol{x}$. Following the direction, \cite{16} proposed a bilevel stochastic approximation(BSA) algorithm and the computational complexity is $\widetilde{\mathcal{O}}(\epsilon^{-3})$ in terms of the target accuracy $\epsilon$. By introducing warm start initializaiton strategy, stocBiO\cite{13} achieves the computational complexity of $\widetilde{\mathcal{O}}(\epsilon^{-2})$ in terms of the target accuracy $\epsilon$. In addition, in order to improve the computational complexity of the simple SGD type double loop algorithms, variance reduction is used for both UL variable $\boldsymbol{x}$ and LL variable $\boldsymbol{y}$ in \cite{21}, and the computational complexity is $\widetilde{\mathcal{O}}(\epsilon^{-1.5})$ in terms of the target accuracy $\epsilon$.

Although many stochastic algorithms have been proposed to solve problem (\ref{eq2}), little attention is paid to the the hypergradient estimation process except a few attempts recently. \cite{22} noticed that to estimate the hypergradient, the existing single loop algorithms and double loop algorithms require an additional loop to estimate Hessian inverse. By studying the well known hypergradient approximation methods such as BP\cite{11}, CG\cite{11}, and NS\cite{16}, they identified a general approximation formulation of hypergradient computation. Based on the formulation and the momentum based variance reduction technique used in \cite{23}, they proposed a fully single loop algorithm(FSLA) in which an additional loop to estimate the Hessian inverse is not required.

In the paper, our goal is to design a new simple stochastic gradient descent(SGD) type algorithm for problem (\ref{eq2}), which can achieve lower computational complexity compared with the existing simple SGD type algorithms(simple SGD type algorithms refer to stochastic algorithms that do not use momentum and variance techniques). We notice that the existing simple SGD type algorithms(i.e., BSA, TTSA, and stocBiO) for solving problem (\ref{eq2}) estimate the hypergradient via stochastic estimation of Neumann series, this kind of hypergradient estimation method termed as Stochastic NS in the paper. Inspired by the work in \cite{22}, we study some hypergradient estimation methods: Stochastic BP, Stochastic NS, and SGD-based Estimation, and find that the algorithm has the potential to converge faster by using SGD-based Estimation to estimate the hypergradient. Based on the finding, by using warm start initialization strategy, a new simple SGD type algorithm SSGD based on SGD-based Estimation is proposed. Through our theoretical analysis, SSGD can achieve lower complexity compared with the existing simple SGD type algorithms. The main contributions of this paper are:
\begin{enumerate}
\item We propose to use SGD-based Estimation to estimate the hypergradient.

\item We propose a new simple SGD type algorithm SSGD.

\item We prove that SSGD can achieve lower computational complexity compared with the existing simple SGD type algorithms. As shown in Table \ref{table1}, the computational complexity of SSGD to reach an $\epsilon$-stationary point outperforms that of BSA and stocBiO by an order of $\mathcal{O}(\epsilon^{-1}\log({1}/{\epsilon}))$ and $\mathcal{O}(\log({1}/{\epsilon}))$, respectively. In addition, in terms of the target accuracy $\epsilon$, the computational complexity of SSGD outperforms that of TTSA by an order of $\mathcal{O}(\epsilon^{-0.5}\log({1}/{\epsilon}))$.

\item We perform experiments to validate our theoretical results and demonstrate the application of our algorithm in hyperparameter optimization.
\end{enumerate}

The rest of this paper is organized as follows. In Section \ref{section:2}, we first make some settings for problem (\ref{eq2}), and specify the problem to study in this paper. Then, for problem (\ref{eq2}), we introduce some hypergradient estimation methods and analyze them. In Section \ref{section:3}, we propose a new simple SGD type algorithm for problem (\ref{eq2}). The convergence and complexity results of our proposed algorithm are provided in Section \ref{section:4}. Experimental results are provided in Section \ref{section:5}. In Section \ref{section:6}, we conclude and summarize the paper.

\textbf{Notations.} We use $\|\cdot \|$ to denote the $l_2$ norm for vectors and spectral norm for matrices. For a twice differentiable function $f(\boldsymbol{x}, \boldsymbol{y}):\mathbb{R}^p \times \mathbb{R}^q\rightarrow \mathbb{R}$, $\nabla f(\boldsymbol{x}, \boldsymbol{y})$ denotes the gradient of $f(\boldsymbol{x}, \boldsymbol{y})$ taken w.r.t. all the variables $(\boldsymbol{x}, \boldsymbol{y})$, $\nabla_{\boldsymbol{x}} f(\boldsymbol{x}, \boldsymbol{y})$(resp. $\nabla_{\boldsymbol{y}} f(\boldsymbol{x}, \boldsymbol{y})$) denotes its partial derivate taken w.r.t.~$\boldsymbol{x}$(resp. $\boldsymbol{y}$), $\nabla_{\boldsymbol{x}} \nabla_{\boldsymbol{y}} f(\boldsymbol{x}, \boldsymbol{y})$  denotes the Jacobian of $\nabla_{\boldsymbol{x}} f(\boldsymbol{x}, \boldsymbol{y})$ w.r.t. $\boldsymbol{y}$, $\nabla_{\boldsymbol{y}}^2 f(\boldsymbol{x}, \boldsymbol{y})$ denotes the Hessian matrix of $f(\boldsymbol{x}, \boldsymbol{y})$ w.r.t. $\boldsymbol{y}$. Furthermore, for a twice differentiable function $f(\boldsymbol{x}, \boldsymbol{y}) := \mathbb{E}_{\boldsymbol{\xi}}[F(\boldsymbol{x}, \boldsymbol{y};\boldsymbol{\xi})]: \mathbb{R}^p \times \mathbb{R}^q\rightarrow \mathbb{R}$, which takes expected value w.r.t. random variable $\boldsymbol{\xi}$, we define 
\begin{equation*}
\nabla_{\boldsymbol{x}} F(\boldsymbol{x}, \boldsymbol{y}; \mathcal{D}) = \frac{1}{|\mathcal{D}|}\sum\limits_{i=1}^{|\mathcal{D}|}\nabla_{\boldsymbol{x}}F(\boldsymbol{x}, \boldsymbol{y}; \boldsymbol{\xi}_i)
\end{equation*}
where $\mathcal{D}=\{\boldsymbol{\xi}_1, \ldots, \boldsymbol{\xi}_{|\mathcal{D}|}\}$ is a sample set with samples $\boldsymbol{\xi}_i$, $i=1, \ldots, |\mathcal{D}|$ sampled over the distribution of $\boldsymbol{\xi}$, and $|\mathcal{D}|$ is the sample size of $\mathcal{D}$.  Similarly, the definitions of $\nabla_{\boldsymbol{y}} F(\boldsymbol{x}, \boldsymbol{y}; \mathcal{D})$, $\nabla_{\boldsymbol{x}} \nabla_{\boldsymbol{y}} F(\boldsymbol{x}, \boldsymbol{y}; \mathcal{D})$, and $\nabla_{\boldsymbol{y}}^2 F(\boldsymbol{x}, \boldsymbol{y}; \mathcal{D})$ can be obtained.

\begin{table}[!t]   
\caption{Computational complexity of the simple SGD type algorithms to reach an $\epsilon$-stationary point of $\Phi(\boldsymbol{x})$ for problem (\ref{eq2})(the computational complexity of BSA, TTSA, and stocBiO can be obtained from Table 1 in \cite{24}). Note that $\kappa$ denotes the condition number; the notation $\kappa^p$ denotes a polynomial function of $\kappa$ since \cite{7} does not provide the dependence on the condition number $\kappa$}
\label{table1}
\centering
\begin{tabular}{cccc}
\hline  Algorithm & loop & batch size & computational complexity\\   
\hline \textbf{BSA}\cite{16} &double & $\mathcal{O}(1)$ & $\mathcal{O}(\kappa^9 \epsilon^{-3}\log(1/\epsilon))$ \\   
 \textbf{TTSA} \cite{7} & single & $\mathcal{O}(1)$ & $\mathcal{O}(\kappa^p\epsilon^{-{5}/{2}}\log(1/\epsilon))$\\
 \textbf{stocBiO} \cite{13}& double & $\mathcal{O}(\epsilon^{-1})$ & $\mathcal{O}(\kappa^9 \epsilon^{-2}\log(1/\epsilon))$\\
 \textbf{SSGD(Ours)} & double & $\mathcal{O}(\epsilon^{-1})$ & $\mathcal{O}(\kappa^9 \epsilon^{-2})$ \\
\hline
\end{tabular}
\end{table}

\section{Analysis} \label{section:2}
In this section, we first make some settings for problem (\ref{eq2}), and then introduce some hypergradient estimation methods for problem (\ref{eq2}). Finally, we analyze these hypergradient estimation methods.

\subsection{Preliminaries} 
In the following, we make some assumptions on the objective functions in problem (\ref{eq2}).
\begin{assumption}\label{assum:1}
The functions $g(\boldsymbol{x}, \boldsymbol{y})$ and $G(\boldsymbol{x}, \boldsymbol{y}; \boldsymbol{\zeta})$ are $\mu$-strongly convex w.r.t. $\boldsymbol{y}$ for any $\boldsymbol{x}$ and $\boldsymbol{\zeta}$, and $\Phi(\boldsymbol{x})$ is nonconvex w.r.t. $\boldsymbol{x}$.
\end{assumption}

\begin{assumption}\label{assum:2}
The functions $f(\boldsymbol{x}, \boldsymbol{y})$ and $g(\boldsymbol{x}, \boldsymbol{y})$ satisfy
\begin{itemize}
\item $f(\boldsymbol{x}, \boldsymbol{y})$ is $M$-Lipschitz, i.e., for $\boldsymbol{z}_1 := [\boldsymbol{x}_1; \boldsymbol{y}_1]$, $\boldsymbol{z}_2 := [\boldsymbol{x}_2; \boldsymbol{y}_2]$, we have
\begin{equation*}
\|f(\boldsymbol{x}_1, \boldsymbol{y}_1) - f(\boldsymbol{x}_2, \boldsymbol{y}_2) \| \le M \|\boldsymbol{z}_1 - \boldsymbol{z}_2 \|.
\end{equation*}
\item $\nabla f(\boldsymbol{x}, \boldsymbol{y})$ and $\nabla g(\boldsymbol{x}, \boldsymbol{y})$ are $L$-Lipschitz, i.e., for $\boldsymbol{z}_1 := [\boldsymbol{x}_1; \boldsymbol{y}_1]$, $\boldsymbol{z}_2 := [\boldsymbol{x}_2; \boldsymbol{y}_2]$, we have 
\begin{align*}
&\| \nabla f(\boldsymbol{x}_1, \boldsymbol{y}_1)- \nabla f(\boldsymbol{x}_2, \boldsymbol{y}_2)\| \le L \|\boldsymbol{z}_1 -\boldsymbol{z}_2 \|,\\
&\| \nabla g(\boldsymbol{x}_1, \boldsymbol{y}_1) - \nabla g(\boldsymbol{x}_2, \boldsymbol{y}_2) \| \le L \|\boldsymbol{z}_1 - \boldsymbol{z}_2\|.
\end{align*}
\item $\nabla_{\boldsymbol{x}} \nabla_{\boldsymbol{y}} g(\boldsymbol{x}, \boldsymbol{y})$ and $\nabla_{\boldsymbol{y}}^{2}g(\boldsymbol{x}, \boldsymbol{y})$ are $\tau$- and $\rho$-Lipschitz, i.e., for $\boldsymbol{z}_1 := [\boldsymbol{x}_1; \boldsymbol{y}_1]$, $\boldsymbol{z}_2 := [\boldsymbol{x}_2; \boldsymbol{y}_2]$, we have 
\begin{align*}
& \| \nabla_{\boldsymbol{x}} \nabla_{\boldsymbol{y}} g(\boldsymbol{x}_1, \boldsymbol{y}_1) - \nabla_{\boldsymbol{x}} \nabla_{\boldsymbol{y}} g(\boldsymbol{x}_2, \boldsymbol{y}_2) \| \le \tau \|\boldsymbol{z}_1 - \boldsymbol{z}_2 \|, \\
& \| \nabla_{\boldsymbol{y}}^2 g(\boldsymbol{x}_1, \boldsymbol{y}_1) - \nabla_{\boldsymbol{y}}^2 g(\boldsymbol{x}_2, \boldsymbol{y}_2) \| \le \rho \| \boldsymbol{z}_1 - \boldsymbol{z}_2 \|.
\end{align*}
\end{itemize}
Furthermore, the same assumption holds for $F(\boldsymbol{x}, \boldsymbol{y}; \boldsymbol{\xi})$ and $G(\boldsymbol{x}, \boldsymbol{y}; \boldsymbol{\zeta})$ for any given $\boldsymbol{\xi}$ and $\boldsymbol{\zeta}$.
\end{assumption}

In fact, for problem $(\ref{eq2})$, Assumptions \ref{assum:1} and \ref{assum:2} are widely adopted to establish the convergence of the stochastic algorithms such as BSA, TTSA, and stocBiO. Based on Assumptions \ref{assum:1} and \ref{assum:2}, the differentiability of $\Phi(\boldsymbol{x})$, and the Lipschitz continuity of $\boldsymbol{y}^*(\boldsymbol{x})$ and $\nabla \Phi(\boldsymbol{x})$ can be ensured, as shown below.

\begin{proposition} \label{pro:1}
For problem (\ref{eq2}), suppose Assumptions $\ref{assum:1}$ and $\ref{assum:2}$ hold. Then $\Phi(\boldsymbol{x})$ is differentiable, and the hypergradient $\nabla \Phi(\boldsymbol{x})$ takes the form of 
\begin{equation}  \label{eq3}
\nabla \Phi(\boldsymbol{x}) = \nabla_{\boldsymbol{x}} f(\boldsymbol{x}, \boldsymbol{y}^*(\boldsymbol{x})) - \nabla_{\boldsymbol{x}} \nabla_{\boldsymbol{y}} g(\boldsymbol{x}, \boldsymbol{y}^*(\boldsymbol{x}))\boldsymbol{v}(\boldsymbol{x}, \boldsymbol{y}^*(\boldsymbol{x}))
\end{equation}
where $\boldsymbol{v}(\boldsymbol{x}, \boldsymbol{y}^*(\boldsymbol{x})) = [\nabla_{\boldsymbol{y}}^2 g(\boldsymbol{x}, \boldsymbol{y}^*(\boldsymbol{x}))]^{-1}\nabla_{\boldsymbol{y}} f(\boldsymbol{x}, \boldsymbol{y}^*(\boldsymbol{x}))$ is the solution of the following linear system:
\begin{equation*}
\nabla_{\boldsymbol{y}}^2 g\left(\boldsymbol{x}, \boldsymbol{y}^*(\boldsymbol{x})\right)\boldsymbol{v} = \nabla_{\boldsymbol{y}} f\left(\boldsymbol{x}, \boldsymbol{y}^*(\boldsymbol{x})\right).
\end{equation*}
Furthermore, $\boldsymbol{y}^*(\boldsymbol{x})$ and $\nabla \Phi(\boldsymbol{x})$ are Lipschitz continuous, i.e., for any $\boldsymbol{x}_1$, $\boldsymbol{x}_2\in \mathbb{R}^p$, we have 
\begin{align*}
&\| \boldsymbol{y}^*(\boldsymbol{x}_1) - \boldsymbol{y}^*(\boldsymbol{x}_2) \| \le \frac{L}{\mu} \|\boldsymbol{x}_1 - \boldsymbol{x}_2\|,\\
&\| \nabla \Phi(\boldsymbol{x}_1) - \nabla \Phi(\boldsymbol{x}_2) \| \le L_{\Phi} \| \boldsymbol{x}_1 - \boldsymbol{x}_2 \|,
\end{align*}
where $L_{\Phi} = L + (2L^2 + \tau M^2)/{\mu} +  (\rho L M + L^3 + \tau M L)/{\mu^2} + (\rho L^2 M)/{\mu^3}$.
\end{proposition}

The proof of Proposition \ref{pro:1} is in the supplementary material. From Assumption \ref{assum:1} and Proposition \ref{pro:1}, we know that $\Phi(\boldsymbol{x})$ is nonconvex and differentiable. Therefore, we want the algorithms to find an $\epsilon$-stationary point defined as follows.

\begin{definition}\label{def:1}
For problem (\ref{eq2}), we call $\bar{\boldsymbol{x}}$ an $\epsilon$-stationary point of $\Phi(\boldsymbol{x})$ if $\mathbb{E}\| \nabla \Phi(\bar{\boldsymbol{x}}) \|^2 \le \epsilon$.
\end{definition}

To facilitate our analysis of the complexity of the algorithm when it reaches the $\epsilon$ stationary point, we borrow the metrics of complexity as defined in Definition 2 of \cite{13}.

\begin{definition} \label{def:2}
For a stochastic function $F(\boldsymbol{x}, \boldsymbol{y};\boldsymbol{\xi})$ with $\boldsymbol{x}\in \mathbb{R}^p$, $\boldsymbol{y}\in \mathbb{R}^q$, and $\boldsymbol{\xi}$ being a random variable, and a vector $\boldsymbol{v}\in \mathbb{R}^q$, let Gc$(F, \epsilon)$ be the number of partial derivatives $\nabla_{\boldsymbol{x}} F(\boldsymbol{x}, \boldsymbol{y}; \boldsymbol{\xi})$ or $\nabla_{\boldsymbol{y}} F(\boldsymbol{x}, \boldsymbol{y};\boldsymbol{\xi})$, and let Jv$(F, \epsilon)$ and Hv$(F, \epsilon)$ be the number of Jacobian vector products $\nabla_{\boldsymbol{x}} \nabla_{\boldsymbol{y}} F(\boldsymbol{x}, \boldsymbol{y}; \boldsymbol{\xi})\boldsymbol{v}$ and Hessian vector products $\nabla_{\boldsymbol{y}}^2 F(\boldsymbol{x}, \boldsymbol{y};\boldsymbol{\xi})\boldsymbol{v}$, respectively.
\end{definition}

To prove the convergence of the stochastic algorithms, as adopted in \cite{7, 13, 16}, we make the following assumptions on the stochastic derivatives.

\begin{assumption} \label{assum:3}
For problem (\ref{eq2}), the stochastic derivatives $\nabla F(\boldsymbol{x}, \boldsymbol{y}; \boldsymbol{\xi})$, $\nabla G(\boldsymbol{x}, \boldsymbol{y}; \boldsymbol{\zeta})$, $\nabla_{\boldsymbol{x}} \nabla_{\boldsymbol{y}} G(\boldsymbol{x}, \boldsymbol{y}; \boldsymbol{\zeta})$, and $\nabla_{\boldsymbol{y}}^2 G(\boldsymbol{x}, \boldsymbol{y}; \boldsymbol{\zeta})$ are unbiased estimates of $\nabla f(\boldsymbol{x}, \boldsymbol{y})$, $\nabla g(\boldsymbol{x}, \boldsymbol{y})$, $\nabla_{\boldsymbol{x}} \nabla_{\boldsymbol{y}} g(\boldsymbol{x}, \boldsymbol{y})$, and $\nabla_{\boldsymbol{y}}^2 g(\boldsymbol{x}, \boldsymbol{y})$, respectively. Furthermore, their variances satisfy
\begin{align*}
&\mathbb{E} \| \nabla F(\boldsymbol{x}, \boldsymbol{y}; \boldsymbol{\xi}) - \nabla f(\boldsymbol{x}, \boldsymbol{y}) \|^2 \le \sigma_f^2, \\
&\mathbb{E} \|\nabla G(\boldsymbol{x}, \boldsymbol{y}; \boldsymbol{\zeta}) - \nabla g(\boldsymbol{x}, \boldsymbol{y}) \|^2 \le \sigma_g^2,\\
&\mathbb{E}\| \nabla_{\boldsymbol{x}} \nabla_{\boldsymbol{y}} G(\boldsymbol{x}, \boldsymbol{y}; \boldsymbol{\zeta}) - \nabla_{\boldsymbol{x}} \nabla_{\boldsymbol{y}} g(\boldsymbol{x}, \boldsymbol{y})\|^2  \le \sigma_{g, 1}^2, \\
&\mathbb{E}\| \nabla_{\boldsymbol{y}}^2 G(\boldsymbol{x},\boldsymbol{y}; \boldsymbol{\zeta}) - \nabla_{\boldsymbol{y}}^2 g(\boldsymbol{x}, \boldsymbol{y}) \|^2 \le \sigma_{g, 2}^2,
\end{align*}
where $\sigma_f^2$, $\sigma_g^2$, $\sigma_{g, 1}^2$, and $\sigma_{g, 2}^2$ are positive constants.
\end{assumption}

\subsection{Introduction of Some Hypergradient Estimation Methods}
In the paper, we study the hypergradient estimation methods: Stochastic BP, Stochastic NS, and SGD-based Estimation, where Stochastic NS is adopted to estimate the hypergradient for all the existing simple SGD type algorithms for solving problem (\ref{eq2}), and Stochastic BP and SGD-based Estimation are developed in this paper for problem (\ref{eq2}). Notice that Stochastic NS is actually a stochastic version of NS. In the similar way, we obtain Stochastic BP by introducing stochastic estimation into BP, and based on AID, we obtain SGD-based Estimation by solving linear system with SGD.

In the following, we introduce how each of these hypergradient estimation methods estimates the hypergradient in a general simple SGD type algorithm shown in Algorithm \ref{alg:1}, with problem (\ref{eq2}) satisfying Assumptions 1, 2, and 3.

As shown in Algorithm \ref{alg:1}, given $\boldsymbol{x}^k$, first, an approximation $\boldsymbol{y}^{k, T}$ to the optimal solution $\boldsymbol{y}^*(\boldsymbol{x}^k)$ of the LL problem is obtained through the SGD updates in line 5 of Algorithm \ref{alg:1}. Based on the output $\boldsymbol{y}^{k, T}$, these hypergradient estimation methods estimate the hypergradient as follows.

\begin{algorithm}[H]
\caption{A General Simple SGD Type Algorithm for Problem (\ref{eq2})}  \label{alg:1}
\begin{algorithmic}[1]
\REQUIRE $K$, $T$, $J$, stepsizes $\alpha$ and $\beta$, $0<\eta<1/L$, initialization $\boldsymbol{x^0}\in \mathbb{R}^p$\\
\FOR{$k = 0, \ldots, K-1$}
\STATE Initialize $\boldsymbol{y}^{k, 0}\in \mathbb{R}^q$\\
\FOR {$t = 0, \ldots, T-1$}
\STATE Draw a sample batch $\mathcal{S}_t$ from the distribution of $\boldsymbol{\zeta}$
\STATE Update $\boldsymbol{y}^{k, t+1} = \boldsymbol{y}^{k, t} - \beta \nabla_{\boldsymbol{y}} G(\boldsymbol{x}^{k}, \boldsymbol{y}^{k, t}; \mathcal{S}_t)$
\ENDFOR
\STATE Hypergradient estimation via one of the following :
\STATE \hspace{0.5cm}\textbf{Stochastic BP}: automatic differentiation through the iterations in line 5 of Algorithm \ref{alg:1}, and compute $\boldsymbol{h}_f^k$ via (\ref{eq6})\\
\STATE \hspace{0.5cm}\textbf{Stochastic NS}: \\
\STATE \hspace{0.8cm} 1) solve (\ref{eq7}) by approximating $[\nabla_{\boldsymbol{y}}^2 g(\boldsymbol{x}^k,\boldsymbol{y}^{k,T})]^{-1}$ with the first $J$ terms of the Neumann series in (\ref{eq8}), and obtain an approximate solution $\boldsymbol{v}^{k, J}$ in (\ref{eq9})\\
\STATE \hspace{0.8cm} 2) compute $\boldsymbol{h}_f^k$ via (\ref{eq10})\\
\STATE \hspace{0.5cm}\textbf{SGD-based Estimation}:\\
\STATE \hspace{0.8cm} 1) initialize $\boldsymbol{v}^{k,0}\in \mathbb{R}^q$\\
\STATE \hspace{0.8cm} 2) solve (\ref{eq11}) with $J$ steps of SGD, and obtain an approximate solution $\boldsymbol{v}^{k,J}$, as shown in (\ref{eq12}) \\
\STATE \hspace{0.8cm} 3) compute $\boldsymbol{h}_f^k$ via (\ref{eq13})\\
\STATE Update $\boldsymbol{x}^{k+1} = \boldsymbol{x}^k - \alpha \boldsymbol{h}_f^k$
\ENDFOR
\end{algorithmic}
\end{algorithm}

\textbf{Stochastic BP.} Notice that $\Phi(\boldsymbol{x})=f(\boldsymbol{x}, \boldsymbol{y}^*(\boldsymbol{x}))$, and that
\begin{align} \label{eq4}
\nabla \Phi(\boldsymbol{x}^k)&= \nabla_{\boldsymbol{x}} f(\boldsymbol{x}^k, \boldsymbol{y}^*(\boldsymbol{x}^k)) \\
&\qquad + (\nabla \boldsymbol{y}^*(\boldsymbol{x}^k))^\top \nabla_{\boldsymbol{y}} f(\boldsymbol{x}^k, \boldsymbol{y}^*(\boldsymbol{x}^k)).\nonumber
\end{align}
For $\nabla \Phi(\boldsymbol{x}^k)$ in the form of (\ref{eq4}), Stochastic BP estimates the hypergradient by automatic differentiation. Specifically, since $\boldsymbol{y}^{k, T}$ has a dependence on $\boldsymbol{x}^k$ through the SGD updates in line 5 of Algorithm \ref{alg:1}, Stochastic BP computes $(\nabla \boldsymbol{y}^{k, T}(\boldsymbol{x}^k))^\top$ by automatic differentiation along these SGD updates, and $(\nabla \boldsymbol{y}^{k, T}(\boldsymbol{x}^k))^\top$ can be obtained as follows:
\begin{equation} \label{eq5}
(\nabla \boldsymbol{y}^{k, T}(\boldsymbol{x}^k))^\top = -  \beta \sum\limits_{t=0}^{T-1}\nabla_{\boldsymbol{x}} \nabla_{\boldsymbol{y}} G(\boldsymbol{x}^k, \boldsymbol{y}^{k, t}; \mathcal{S}_t)\mathbf{A}_{T-1-t}
\end{equation}
where $\mathbf{A}_{T-1-t} = \prod_{i=t+1}^{T-1}(\mathbf{I} - \beta \nabla_{\boldsymbol{y}}^2 G(\boldsymbol{x}^{k},\boldsymbol{y}^{k, i}; \mathcal{S}_i))$ with $\mathbf{I}$ being the identity matrix, and $\mathbf{A}_{0}=\mathbf{I}$(for the derivation, see Appendix D.2 in \cite{13}). Then, 
\begin{align}  \label{eq6}
\boldsymbol{h}_f^{k}& = \nabla_{\boldsymbol{x}} F(\boldsymbol{x}^{k}, \boldsymbol{y}^{k, T}; \mathcal{D}_F) \\
& \qquad + (\nabla \boldsymbol{y}^{k, T}(\boldsymbol{x}^k))^\top \nabla_{\boldsymbol{y}} F(\boldsymbol{x}^{k}, \boldsymbol{y}^{k, T}; \mathcal{D}_F) \nonumber
\end{align}
is used to estimate the hypergradient $\nabla \Phi(\boldsymbol{x}^k)$ in (\ref{eq4}), where $\mathcal{D}_F$ is a sample set sampled from the distribution of $\boldsymbol{\xi}$ and 
$(\nabla \boldsymbol{y}^{k,T}(\boldsymbol{x}^k))^\top$ is given in (\ref{eq5}).

In fact, Stochastic BP is a variant of BP, and $\boldsymbol{h}_f^k$ in (\ref{eq6}) can be directly obtained by introducing stochastic estimation into the hypergradient approximation formula of BP in Proposition 2 in \cite{13}.

\textbf{Stochastic NS.} From Proposition \ref{pro:1}, we have 
\begin{align}  \label{eq6-1}
\nabla \Phi(\boldsymbol{x}^k)& = \nabla_{\boldsymbol{x}} f(\boldsymbol{x}^k, \boldsymbol{y}^*(\boldsymbol{x}^k)) \\
& \qquad - \nabla_{\boldsymbol{x}} \nabla_{\boldsymbol{y}} g(\boldsymbol{x}^k, \boldsymbol{y}^*(\boldsymbol{x}^k))\boldsymbol{v}(\boldsymbol{x}^k, \boldsymbol{y}^*(\boldsymbol{x}^k))\nonumber
\end{align}
with $\boldsymbol{v}(\boldsymbol{x}^k, \boldsymbol{y}^*(\boldsymbol{x}^k)) = [\nabla_{\boldsymbol{y}}^2 g(\boldsymbol{x}^k, \boldsymbol{y}^*(\boldsymbol{x}^k))]^{-1}\nabla_{\boldsymbol{y}} f(\boldsymbol{x}^k, \boldsymbol{y}^*(\boldsymbol{x}^k))$. For $\nabla \Phi(\boldsymbol{x}^k)$ in the form of (\ref{eq6-1}), Stochastic NS estimates the hypergradient by solving the following linear system: 
\begin{equation}  \label{eq7}
\nabla_{\boldsymbol{y}}^2 g(\boldsymbol{x}^{k}, \boldsymbol{y}^{k, T}) \boldsymbol{v} = \nabla_{\boldsymbol{y}} f(\boldsymbol{x}^{k}, \boldsymbol{y}^{k, T}).
\end{equation}
By approximating Hessian inverse $[\nabla_{\boldsymbol{y}}^2 g(\boldsymbol{x}^k, \boldsymbol{y}^{k, T})]^{-1}$ with the first $J$ terms of the following Neumann series: 
\begin{equation} \label{eq8}
\sum_{j=0}^{\infty}\eta(\mathbf{I}-\eta \nabla_{\boldsymbol{y}}^2 g(\boldsymbol{x}^k, \boldsymbol{y}^{k,T}))^j
\end{equation}
and computing $\nabla_{\boldsymbol{y}} F(\boldsymbol{x}^{k}, \boldsymbol{y}^{k, T}; \mathcal{D}_F)$ over a sample set $\mathcal{D}_F$ sampled from the distribution of $\boldsymbol{\xi}$, an approximate solution $\boldsymbol{v}^{k, J}$ of (\ref{eq7}) is obtained, and takes the form of 
\begin{equation}  \label{eq9}
\boldsymbol{v}^{k, J} = \eta \sum\limits_{j=0}^{J-1}\prod_{i=J-j}^{J-1}(\mathbf{I} - \eta \nabla_{\boldsymbol{y}}^2 G(\boldsymbol{x}^{k}, \boldsymbol{y}^{k, T}; \mathcal{B}_i))\boldsymbol{v}_0
\end{equation}
where $\prod_{i=J}^{J-1}(\cdot) = \mathbf{I}$, $\boldsymbol{v}_0 = \nabla_{\boldsymbol{y}} F(\boldsymbol{x}^{k}, \boldsymbol{y}^{k, T}; \mathcal{D}_F)$, and sample sets $\mathcal{B}_i, i=1, \ldots, J-1$ are sampled from the distribution of $\boldsymbol{\zeta}$. Then, 
\begin{equation} \label{eq10}
\boldsymbol{h}_f^{k} = \nabla_{\boldsymbol{x}} F(\boldsymbol{x}^{k}, \boldsymbol{y}^{k, T}; \mathcal{D}_F) -  \nabla_{\boldsymbol{x}} \nabla_{\boldsymbol{y}} G(\boldsymbol{x}^k, \boldsymbol{y}^{k, T}; \mathcal{D}_G)\boldsymbol{v}^{k, J}
\end{equation}
is used to estimate the hypergradient $\nabla \Phi(\boldsymbol{x}^k)$ in (\ref{eq6-1}), where $\mathcal{D}_F$ is the aforementioned sample set, $\mathcal{D}_G$ is a sample set sampled from the distribution of $\boldsymbol{\zeta}$, and $\boldsymbol{v}^{k,J}$ is given in (\ref{eq9}). Note that sample sets $\{\mathcal{D}_F, \mathcal{D}_G, \mathcal{B}_i(i=1, \ldots, J-1)\}$ are mutually independent.

In fact, $\boldsymbol{h}_f^k$ in (\ref{eq10}) can be directly obtained by introducing stochastic estimation into the hypergradient approximation formula of NS in eq.(3) in \cite{21}, which is the reason we call this hypergradient estimation method Stochastic NS. Furthermore, for more details about Stochastic NS, see subsection 2.2 in \cite{13}.

\textbf{SGD-based Estimation.} SGD-based Estimation estimates the hypergradient in the way similar to Stochastic NS except that the approximate solution to (\ref{eq7}) is obtained by using SGD. Notice that the solution to (\ref{eq7}) is equivalent to the optimal solution of the following optimization problem:
\begin{equation} \label{eq11}
\min\limits_{\boldsymbol{v}} \frac{1}{2} \boldsymbol{v}^\top \nabla_{\boldsymbol{y}}^2 g(\boldsymbol{x}^{k}, \boldsymbol{y}^{k, T})\boldsymbol{v} - \boldsymbol{v}^\top \nabla_{\boldsymbol{y}} f(\boldsymbol{x}^{k}, \boldsymbol{y}^{k, T}).
\end{equation}
SGD-based Estimation obtains an approximate solution $\boldsymbol{v}^{k, J}$ to (\ref{eq11}) by solving (\ref{eq11}) with $J$ steps of SGD as follows:
\begin{align}   \label{eq12}
\boldsymbol{v}^{k, j+1} &= (\mathbf{I} - \eta \nabla_{\boldsymbol{y}}^2 G(\boldsymbol{x}^{k}, \boldsymbol{y}^{k, T}; \mathcal{B}_j))\boldsymbol{v}^{k, j} \\
&\qquad + \eta \nabla_{\boldsymbol{y}} F(\boldsymbol{x}^{k}, \boldsymbol{y}^{k, T}; \mathcal{D}_{F, j}),\quad j=0,\ldots, J-1   \nonumber
\end{align}
where $\eta$ is the step size, and for each $j \in \{0, \ldots, J-1\}$, $\mathcal{D}_{F, j}$ and $\mathcal{B}_j$ are the sample sets sampled from the distributions of $\boldsymbol{\xi}$ and $\boldsymbol{\zeta}$, respectively. Then, 
\begin{equation} \label{eq13}
\boldsymbol{h}_f^k = \nabla_{\boldsymbol{x}} F(\boldsymbol{x}^k, \boldsymbol{y}^{k, T}; \mathcal{D}_F) - \nabla_{\boldsymbol{x}} \nabla_{\boldsymbol{y}} G(\boldsymbol{x}^k, \boldsymbol{y}^{k, T}; \mathcal{D}_G)\boldsymbol{v}^{k, J}
\end{equation}
is used to estimate the hypergradient, where $\mathcal{D}_F$ and $\mathcal{D}_G$ are the sample sets sampled from the distributions of $\boldsymbol{\xi}$ and $\boldsymbol{\zeta}$, respectively.

\subsection{Analysis of Hypergradient Estimation Methods}

Define $\mathcal{F}_{k,1}^{'} = \sigma\{ \boldsymbol{x}^0, \ldots, \boldsymbol{y}^{0, T}, \ldots, \boldsymbol{x}^k, \ldots, \boldsymbol{y}^{k,T}\}$, $\mathcal{F}_{k,2}^{'} = \sigma\{ \boldsymbol{x}^0, \ldots, \boldsymbol{y}^{0, T}, \ldots, \boldsymbol{v}^{0, J}, \ldots, \boldsymbol{x}^k, \ldots, \boldsymbol{y}^{k,T}\}$, where $\sigma\{\cdot\}$ denotes the $\sigma$-algebra generated by random variables.

Then, similar to the discussion in the proof of Lemma 3 in the supplementary material, from the iteration in line 16 of Algorithm \ref{alg:1}, we can obtain formula (19) in the supplementary material in which $\bar{\boldsymbol{h}}_f^k = \mathbb{E}[\boldsymbol{h}_f^k| \mathcal{F}_{k,1}^{'}]$ for Stochastic BP and Stochastic NS, and $\bar{\boldsymbol{h}}_f^k = \mathbb{E}[\boldsymbol{h}_f^k| \mathcal{F}_{k,2}^{'}]$ for SGD-based Estimation. By shifting the terms, we have 
\begin{align*}
&\mathbb{E}\| \nabla \Phi(\boldsymbol{x}^k) \|^2 \\
&\le \frac{2}{\alpha}\big(\mathbb{E}[\Phi(\boldsymbol{x}^k)] - \mathbb{E}[\Phi(\boldsymbol{x}^{k+1})]\big) - (1 - \alpha{L_{\Phi}})\mathbb{E}\| \bar{\boldsymbol{h}}_f^k\|^2 \\
&\qquad  + \mathbb{E}\| \nabla \Phi(\boldsymbol{x}^k) - \bar{\boldsymbol{h}}_f^k \|^2  + \alpha{L_{\Phi}} \mathbb{E}\| \boldsymbol{h}_f^k - \bar{\boldsymbol{h}}_f^k \|^2 
\end{align*}
where $\alpha$ is in Algorithm \ref{alg:1}, $L_{\Phi}$ is in Proposition \ref{pro:1}.

The upper bound of $\mathbb{E}\| \nabla \Phi(\boldsymbol{x}^k) \|^2$ involves the term $\| \nabla \Phi(\boldsymbol{x}^k) - \bar{\boldsymbol{h}}_f^k \|$. In the following, the upper bound of  $\| \nabla \Phi(\boldsymbol{x}^k) - \bar{\boldsymbol{h}}_f^k \|$ is provided. For the proof, see Lemma 10 in the supplementary material.

\begin{proposition}  \label{pro:2}
Apply Algorithm \ref{alg:1} to solve problem (\ref{eq2}). Suppose Assumptions \ref{assum:1}, \ref{assum:2}, and \ref{assum:3} hold. let $\mathbf{H}_0 = \mathbf{0}$, $\boldsymbol{v}_0 = \boldsymbol{v}^{k, 0}$, and define $\mathbf{H}^*  = \nabla_{\boldsymbol{x}}\nabla_{\boldsymbol{y}} g(\boldsymbol{x}^k, \boldsymbol{y}^*(\boldsymbol{x}^k))(\nabla_{\boldsymbol{y}}^2 g(\boldsymbol{x}^k, \boldsymbol{y}^*(\boldsymbol{x}^k)))^{-1}$, $\boldsymbol{v}^* = \nabla_{\boldsymbol{y}}^2 g(\boldsymbol{x}^k, \boldsymbol{y}^*(\boldsymbol{x}^k))^{-1}\nabla_{\boldsymbol{y}} f(\boldsymbol{x}^k, \boldsymbol{y}^*(\boldsymbol{x}^k))$, where $\boldsymbol{v}^{k, 0}$ is the initial value to solve problem (\ref{eq11})(see line 13 in Algorithm \ref{alg:1}), and $\boldsymbol{y}^*(\boldsymbol{x}^k)$ is the optimal solution to the LL problem of problem (\ref{eq2}) with $\boldsymbol{x} = \boldsymbol{x}^k$. For $\mu$ , $L$ in Assumptions \ref{assum:1}, \ref{assum:2} and $T$, $J$, $\beta$, $\eta$ in Algorithm \ref{alg:1}, let $\mu <1$, $L<1$, $T$, $J$ be any positive integers, and $\beta = \eta$. Then, 
\begin{itemize}
\item For Stochastic BP, for $\bar{\boldsymbol{h}}_f^k := \mathbb{E}[\boldsymbol{h}_f^k| \mathcal{F}_{k,1}^{'}]$ with $\boldsymbol{h}_f^k$ in (\ref{eq6}), we have 
\begin{align*}
&\|\nabla \Phi(\boldsymbol{x}^k) - \bar{\boldsymbol{h}}_f^k\| \\
& \le L\bigg(1 +\frac{L}{\mu}\bigg) \|\boldsymbol{y}^{k, T} - \boldsymbol{y}^*(\boldsymbol{x}^k)\| + M(1 - \eta \mu)^T \|\mathbf{H}^*\| \\
& + M\eta\bigg(\frac{L}{\mu}\rho + \tau\bigg) \sum\limits_{t=0}^{T-1}(1 - \eta \mu)^t \|\boldsymbol{y}^{k, T-1-t} - \boldsymbol{y}^*(\boldsymbol{x}^k)\|;
\end{align*}
\item For Stochastic NS, for $\bar{\boldsymbol{h}}_f^k := \mathbb{E}[\boldsymbol{h}_f^k| \mathcal{F}_{k,1}^{'}]$ with $\boldsymbol{h}_f^k$ in (\ref{eq10}), we have 
\begin{align*}
&\|\nabla \Phi(\boldsymbol{x}^k) - \bar{\boldsymbol{h}}_f^k\| \\
& \le L\bigg(1 +\frac{L}{\mu}\bigg) \|\boldsymbol{y}^{k, T} - \boldsymbol{y}^*(\boldsymbol{x}^k)\| + M(1 - \eta \mu)^J \|\mathbf{H}^*\| \\
& + M\eta\bigg(\frac{L}{\mu}\rho + \tau\bigg) \sum\limits_{t=0}^{J-1}(1 - \eta \mu)^t \|\boldsymbol{y}^{k, T} - \boldsymbol{y}^*(\boldsymbol{x}^k)\|;
\end{align*}
\item For SGD-based Estimation, for $\bar{\boldsymbol{h}}_f^k := \mathbb{E}[\boldsymbol{h}_f^k| \mathcal{F}_{k,2}^{'}]$ with $\boldsymbol{h}_f^k$ in (\ref{eq13}), we have 
\begin{align*}
&\|\nabla \Phi(\boldsymbol{x}^k) - \bar{\boldsymbol{h}}_f^k\|  \\
&\le \bigg(L +\frac{M}{\mu}\tau \bigg) \|\boldsymbol{y}^{k, T} - \boldsymbol{y}^*(\boldsymbol{x}^k)\| \\
&+ L(1 - \eta \mu)^J \| \boldsymbol{v}^{k, 0} - \boldsymbol{v}^* \| \\
&+ L \eta\bigg(\frac{M}{\mu}\rho + L \bigg) \sum\limits_{t=0}^{J-1}(1 - \eta \mu)^t \|\boldsymbol{y}^{k, T} - \boldsymbol{y}^*(\boldsymbol{x}^k)\|;
\end{align*}
\end{itemize}
$\nabla \Phi(\boldsymbol{x}^k)$ is in (3) with $\boldsymbol{x} = \boldsymbol{x}^k$, $\mu$, $L$, $M$, $\tau$, $\rho$ are given in Assumptions 1, 2.
\end{proposition} 

From Proposition \ref{pro:2}, it is easy to find that there are some differences for Stochastic BP, Stochastic NS, and SGD-based Estimation regarding the upper bound of $\| \nabla \Phi(\boldsymbol{x}^k) - \bar{\boldsymbol{h}}_f^k \|$. For Stochastic BP and Stochastic NS, the upper bound for $\| \nabla \Phi(\boldsymbol{x}^k) - \bar{\boldsymbol{h}}_f^k \|$ involves $M(1 - \eta \mu)^T \|\mathbf{H}^*\|$ and $M(1 - \eta \mu)^J \|\mathbf{H}^*\|$, respectively. While for SGD-based Estimation, the term $L(1-\eta\mu)^J\|\boldsymbol{v}^{k, 0}-\boldsymbol{v}^*\|$ is involved in the upper bound of $\| \nabla \Phi(\boldsymbol{x}^k) - \bar{\boldsymbol{h}}_f^k \|$.

Thus, to ensure that $\| \nabla \Phi(\boldsymbol{x}^k) - \bar{\boldsymbol{h}}_f^k \|$ is a sufficiently small number, a sufficiently large $T$ for Stochastic BP and a sufficiently large $J$ for Stochastic NS may be required. In contrast, for SGD-based Estimation, we notice that $\boldsymbol{v}^{k, 0}$ in the term $L(1 - \eta \mu)^J\|\boldsymbol{v}^{k, 0} - \boldsymbol{v}^*\|$ is the initial value to solve problem (\ref{eq11}). Apply warm start initialization strategy to $\boldsymbol{v}^{k, 0}$(i.e., $\boldsymbol{v}^{k, 0} = \boldsymbol{v}^{k-1, J}$, where $\boldsymbol{v}^{k-1, J}$ is the output of the $J$-th iteration to solve problem (\ref{eq11}) given $\boldsymbol{x}^{k-1}$) can allow us to track the errors, such as $\|\boldsymbol{v}^{k-1, 0} - \boldsymbol{v}^*\|$, \ldots, $\|\boldsymbol{v}^{0, 0}- \boldsymbol{v}^*\|$, in the preceeding loops. Thus, a smaller number $J$ may be able to ensure $\| \nabla \Phi(\boldsymbol{x}^k) - \bar{\boldsymbol{h}}_f^k \|$ to be a sufficiently small number, and there is no restriction on $T$. Notice that warm start initialization strategy is used in stocBiO\cite{13} to initialize the initial values for solving the LL problem, and the improved computational complexity is obtained. 

The above discussion inspires us to think that SGD-based Estimation may be able to allow $\mathbb{E}\| \nabla \Phi(\boldsymbol{x}^k) \|^2$ to be sufficiently small under weaker restrictions on $T$ than Stochastic BP and $J$ than Stochastic NS. Thus, SGD-based Estimation may be able to make Algorithm \ref{alg:1} converge faster than Stochastic BP and Stochastic NS.

Next, we perform experiments on a synthetic bilevel optimization problem, a special case of problem (\ref{eq2}), to evaluate the performance of Algorithm \ref{alg:1} under different hypergradient estimation methods.

\textbf{Synthetic Bilevel Optimization Problem.} We consider the following synthetic bilevel optimization problem:
\begin{equation} \label{eq14}
\min\limits_{\boldsymbol{x}\in \mathbb{R}^p} \Phi(\boldsymbol{x}) := f(\boldsymbol{x}, \boldsymbol{y}^*(\boldsymbol{x})) \quad s.t. ~\boldsymbol{y}^*(\boldsymbol{x}) = { \underset {\boldsymbol{y}\in \mathbb{R}^p} { \operatorname {arg\,min} } \, g(\boldsymbol{x}, \boldsymbol{y})}
\end{equation}
with 
\begin{align*}
&f(\boldsymbol{x}, \boldsymbol{y}) = \frac{1}{|\mathcal{D}_{\text{val}}|}\sum\limits_{(\boldsymbol{u}_i, v_i)\in \mathcal{D}_{\text{val}}}\frac{1}{2}({\boldsymbol{y}}^\top \boldsymbol{u}_i - v_i)^2 + \|\boldsymbol{x}\|^3, \\
&g(\boldsymbol{x}, \boldsymbol{y}) = \frac{1}{|\mathcal{D}_{\text{tr}}|}\sum\limits_{(\boldsymbol{u}_i, v_i)\in \mathcal{D}_{\text{tr}}}\frac{1}{2}(\boldsymbol{y}^\top \boldsymbol{u}_i - v_i)^2 + \frac{r}{2}\| \boldsymbol{y} - \boldsymbol{x} \|^2
\end{align*}
where $r = 0.5$, and $\mathcal{D}_{\text{tr}}$ and $\mathcal{D}_{\text{val}}$ are the data sets constructed as follows: given $\boldsymbol{w}_0\in \mathbb{R}^{p}$, first randomly sample $20000$ data points $\boldsymbol{e}_i\in \mathbb{R}^{p-1}$, $i=1, \ldots, 20000$ from normal distribution with a mean of 0 and a variance of 0.01; then, for each $i\in \{1, \ldots, 20000\}$, we set $\boldsymbol{u}_i = (\boldsymbol{e}_i, 1)$, and construct $v_i \in \mathbb{R}$ by $v_i = \boldsymbol{w}_0^\top \boldsymbol{u}_i + \sigma_i$, where $\sigma_i$ is the Gaussian noise with mean 0 and variance 1; finally, dividing the dataset $\{(\boldsymbol{u}_i, v_i)\}_{i=1}^{20000}$ equally into two parts, we obtain training set $\mathcal{D}_{\text{tr}}$ and validation set $\mathcal{D}_{\text{val}}$, with the training set size $|\mathcal{D}_{\text{tr}}|$ and validation set size $|\mathcal{D}_{\text{val}}|$ being 10000, respectively. 

It is easy to verify that $g(\boldsymbol{x}, \boldsymbol{y})$ is strongly convex w.r.t. $\boldsymbol{y}$. Moreover, for each $\boldsymbol{x}$, the minimum solution of the LL problem is 
\begin{equation*}
\boldsymbol{y}^*(\boldsymbol{x})= (\mathbf{A}_{\text{tr}}+ r\mathbf{I})^{-1}(\boldsymbol{b_{\text{tr}}}+ r\boldsymbol{x})
\end{equation*}
and the hypergradient is 
\begin{equation*}
\nabla \Phi(\boldsymbol{x})= r(\mathbf{A}_{\text{tr}}+ r\mathbf{I})^{-1}(\mathbf{A}_{\text{val}}\boldsymbol{y}^*(\boldsymbol{x}) - \boldsymbol{b_{\text{val}}}) + 3\|\boldsymbol{x}\|  \boldsymbol{x}
\end{equation*}
where 
\begin{align*}
&\mathbf{A}_{\text{val}} = \frac{1}{|\mathcal{D}_{\text{val}}|}\sum\limits_{(\boldsymbol{u}_i, v_i)\in \mathcal{D}_{\text{val}}}\boldsymbol{u}_i(\boldsymbol{u}_i)^\top,\\
&\boldsymbol{b_{\text{val}}} = \frac{1}{|\mathcal{D}_{\text{val}}|}\sum\limits_{(\boldsymbol{u}_i, v_i)\in \mathcal{D}_{\text{val}}} v_i \boldsymbol{u}_i, 
\end{align*}
and $\mathbf{A}_{\text{tr}}$, $\boldsymbol{b_{\text{tr}}}$ can be obtained by replacing $\mathcal{D}_{\text{val}}$ and $|\mathcal{D}_{\text{val}}|$ in $\mathbf{A}_{\text{val}}$ and $\boldsymbol{b_{\text{val}}}$ with $\mathcal{D}_{\text{tr}}$ and $|\mathcal{D}_{\text{tr}}|$.

In the following, we use Algorithm \ref{alg:1} to solve problem (\ref{eq14}), and compare the performance of Algorithm \ref{alg:1} under different hypergradient estimation methods. To be specific, the hypergradient estimation methods that we consider are Stochastic BP, Stochastic NS, and SGD-based Estimation, where for SGD-based Estimation, we use two initialization strategies to initialize $\boldsymbol{v}^{k, 0}$ in line 13 of Algorithm \ref{alg:1}, i.e., non-warm start initialization strategy and warm start initialization strategy. For non-warm start initialization strategy, we set $\boldsymbol{v}^{k, 0}$, $k=0, \ldots, K-1$ to be zero vectors. For warm-start initialization strategy, we set $\boldsymbol{v}^{k,0}=\boldsymbol{v}^{k-1, J}$, $k = 1, \ldots, K-1$, where for each $k \in \{1, \ldots, K-1\}$, $\boldsymbol{v}^{k-1, J}$ is obtained by the iteration steps starting from $\boldsymbol{v}^{k-1, 0}$ in (\ref{eq12}). The experimental details are in the supplementary material.

To simplify the narrative, we refer to Algorithm \ref{alg:1} as Stochastic BP(resp. Stochastic NS) if Algorithm \ref{alg:1} estimates the hypergradient via Stochastic BP(resp. Stochastic NS). Similarly, if Algorithm \ref{alg:1} estimates the hypergradient via SGD-based Estimation and using the non-warm start initialization strategy(resp. warm start initialization strategy) to initialize $\boldsymbol{v}^{k, 0}$, we refer to Algorithm \ref{alg:1} as SGD\_W\_Start\_False(resp. SGD\_W\_Start\_True). Furthermore, we use P\_$T$\_$i$(P\_$J$\_$i$) to denote that the iteration step $T$($J$) is set to be $i$ for P.

\begin{figure*}[!t]
\centering
\subfloat[]{\includegraphics[scale = 0.2]{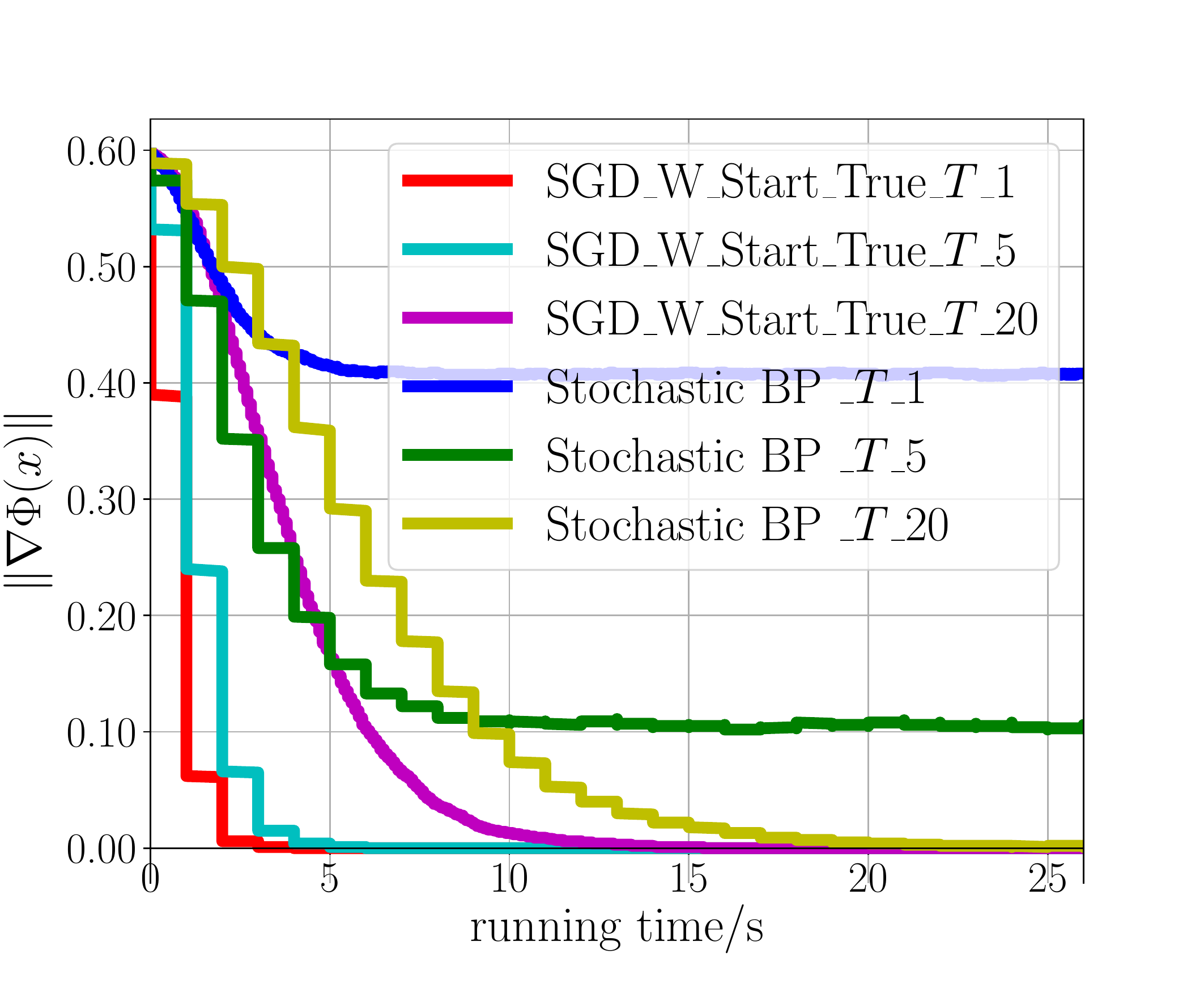}%
\label{fig_first_case}
}
\hfil
\subfloat[]{\includegraphics[scale = 0.2]{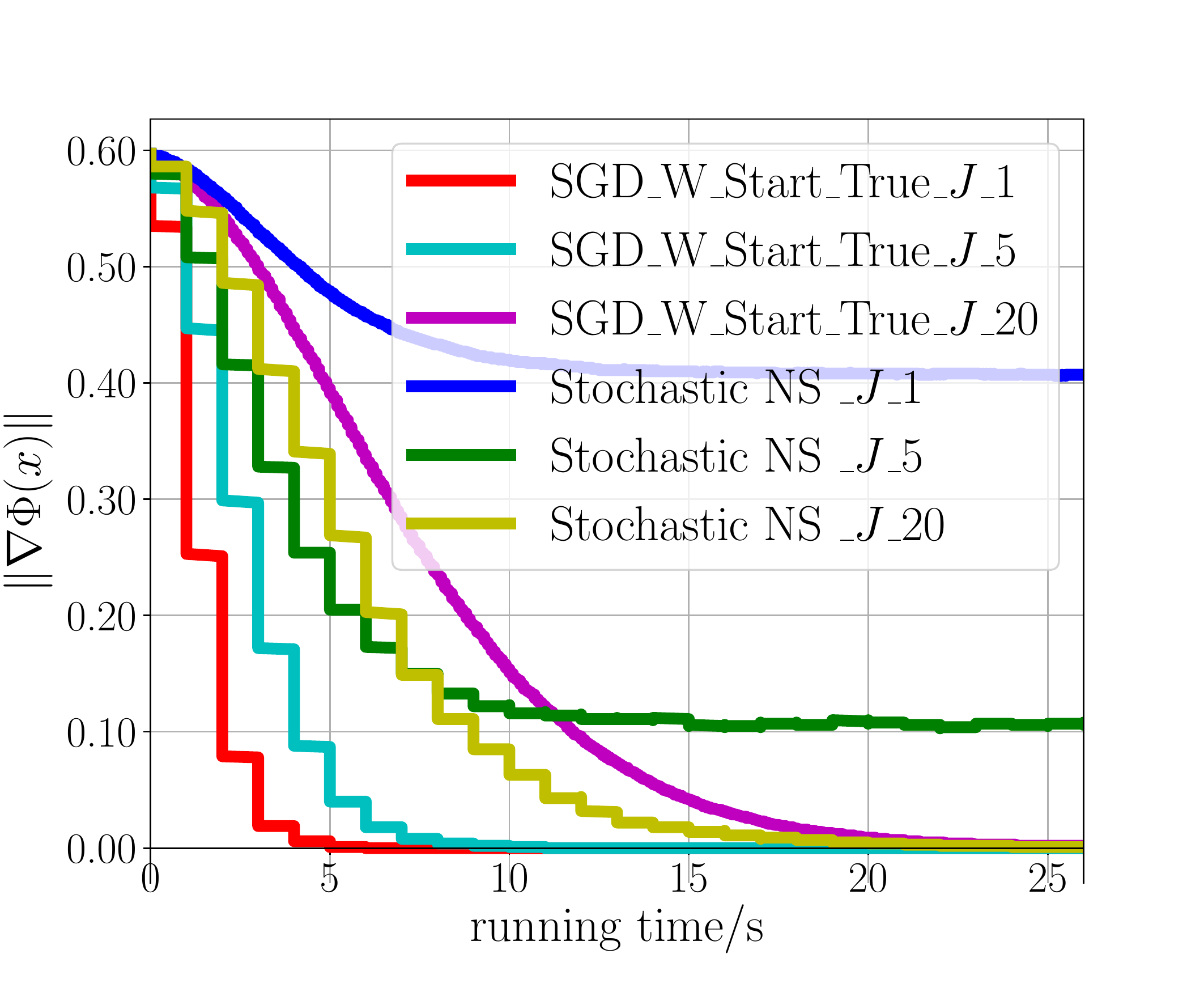}%
\label{fig_second_case}}
\hfil
\subfloat[]{\includegraphics[scale = 0.2]{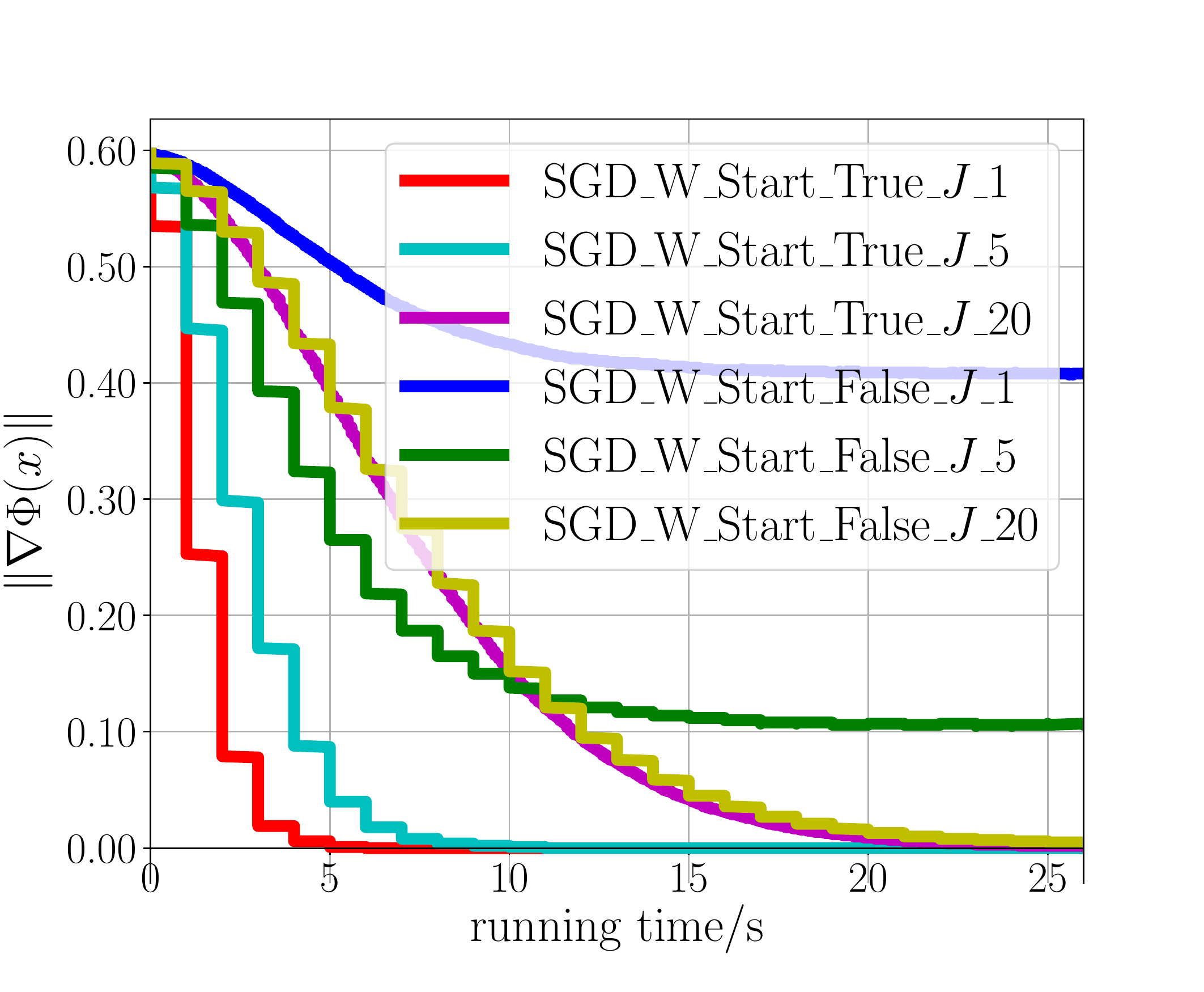}%
\label{fig_third_case}}
\caption{Comparing the performance of Algorithm \ref{alg:1} in solving problem (\ref{eq14}), with data sets constructed by setting $\boldsymbol{w}_0 =(2, 5, 7)$, under differnt hypergradient estimation methods.}
\label{fig:1}
\end{figure*}

In the experiments, data sets $\mathcal{D}_{\text{tr}}$ and $\mathcal{D}_{\text{val}}$ are constructed by setting $\boldsymbol{w}_0 = (2, 5, 7)\in \mathbb{R}^3$. From Fig. \ref{fig:1}\subref{fig_first_case}, we observe that to converge to a stationary point with $\|\nabla \Phi(\boldsymbol{x})\| = 0$, $T=1$ is sufficient for SGD\_W\_Start\_True, and a larger $T$($T=20$) is required for Stochastic BP. Furthermore, SGD\_W\_Start\_True always converges faster than Stochastic BP for $T\in \{1, 5, 20\}$. From Fig. \ref{fig:1}\subref{fig_second_case}, we observe that to converge to a stationary point with $\|\nabla \Phi(\boldsymbol{x})\| = 0$, $J=1$ is sufficient for SGD\_W\_Start\_True, and a larger $J$($J=20$) is required for Stochastic NS. Furthermore, SGD\_W\_Start\_True with $J=1$ has a faster convergence speed than Stochastic NS with $J=20$. From Fig. \ref{fig:1}\subref{fig_third_case}, we observe that to converge to a stationary point with $\|\nabla \Phi(\boldsymbol{x})\| = 0$, $J=1$ is sufficient for SGD\_W\_Start\_True, and a larger $J$($J=20$) is required for SGD\_W\_Start\_False. Furthermore, SGD\_W\_Start\_True with $J=1$ has a faster convergence speed than SGD\_W\_Start\_False with $J=20$.

In conclusion, the above experimental results show that compared with Stochastic BP and Stochastic NS, Algorithm \ref{alg:1} can obtain a faster convergence speed by using SGD-based Estimation to estimate the hypergradient and using warm start initialization strategy for SGD-based Estimation.

\section{Algorithm} \label{section:3}
In Section \ref{section:2}, we find that compared with Stochastic BP and Stochastic NS, SGD-based Estimation can make Algorithm \ref{alg:1} achieve a faster convergence speed. Notice that all the existing simple SGD type algorithms for solving problem (\ref{eq2}) use Stochastic NS to estimate the hypergradient. It is natural to think that using SGD-based Estimation to estimate the hypergradient may allow us to obtain an algorithm with lower computational complexity for problem (\ref{eq2}).

\begin{algorithm}[htb]
\caption{A New Simple SGD Type Algorithm SSGD for Problem (\ref{eq2})}  
\label{alg:2}
\begin{algorithmic}[1]
\REQUIRE  $K\ge 1$, $T\ge 1$, $J\ge 1$, stepsizes $\alpha$, $\beta$, and $\eta$, initializations $\boldsymbol{x}^0\in \mathbb{R}^p$, $\boldsymbol{y}^0\in \mathbb{R}^q$, $\boldsymbol{v}^0\in \mathbb{R}^q$\\
\FOR{$k = 0, \ldots, K-1$}
\STATE $\boldsymbol{x}^{k+1} = \boldsymbol{x}^k - \alpha \boldsymbol{h}_f^k$\\
\STATE Set $\boldsymbol{y}^{k+1, 0}=\boldsymbol{y}^k$\\
\FOR {$t = 0, \ldots, T-1$}
\STATE Draw a sample batch $\mathcal{S}_t$\\
\STATE $\boldsymbol{y}^{k+1, t+1} = \boldsymbol{y}^{k+1, t} - \beta \nabla_{\boldsymbol{y}} G(\boldsymbol{x}^{k+1}, \boldsymbol{y}^{k+1, t}; \mathcal{S}_t)$\\
\ENDFOR
\STATE Set $\boldsymbol{y}^{k+1}=\boldsymbol{y}^{k+1, T}$, $\boldsymbol{v}^{k+1, 0}=\boldsymbol{v}^k$\\
\FOR {$j = 0, \ldots, J-1$}
\STATE Draw sample batches $\mathcal{B}_j$, $\mathcal{D}_{F,j}$\\
\STATE $\boldsymbol{v}^{k+1, j+1} = (\mathbf{I} - \eta \nabla_{\boldsymbol{y}}^2 G(\boldsymbol{x}^{k+1}, \boldsymbol{y}^{k+1};\mathcal{B}_j))\boldsymbol{v}^{k+1, j}$ \\
 \hspace{1.3cm}$ \qquad + \eta \nabla_{\boldsymbol{y}} F(\boldsymbol{x}^{k+1}, \boldsymbol{y}^{k+1}; \mathcal{D}_{F,j})$\\
\ENDFOR 
\STATE Set $\boldsymbol{v}^{k+1}=\boldsymbol{v}^{k+1, J}$ \\
\STATE Draw sample batches $\mathcal{D}_F$, $\mathcal{D}_{G}$\\
\STATE Compute $\boldsymbol{h}_f^{k+1}$ via $(\ref{eq15})$\\
\ENDFOR
\end{algorithmic}
\end{algorithm}

As shown in Algorithm \ref{alg:2}, by using warm start initialization strategy, a new simple SGD type algorithm SSGD based on SGD-based Estimation is proposed, in which warm start is not only used for SGD-based Estimation, but also for solving LL problem.

In Algorithm \ref{alg:2}, given $\boldsymbol{x}^{k+1}$, SSGD runs $T$ steps of SGD to obtain an approximation $\boldsymbol{y}^{k+1, T}$ to the optimal solution $\boldsymbol{y}^*(\boldsymbol{x}^{k+1})$ of the LL problem(see line 6 in Algorithm \ref{alg:2}). Note that the initial value $\boldsymbol{y}^{k+1, 0}$ is initialized by using warm start and is set to be $\boldsymbol{y}^k$, where $\boldsymbol{y}^k$ is equal to $\boldsymbol{y}^{k, T}$, which is the output of the $T$-th iteration starting from $\boldsymbol{y}^{k, 0}$ in line 6 of Algorithm \ref{alg:2}.

After setting $\boldsymbol{y}^{k+1}=\boldsymbol{y}^{k+1, T}$, the hypergradient is estimated via SGD-based Estimation. Specifically, SSGD first obtains an approximation $\boldsymbol{v}^{k+1, J}$ to the optimal solution of the following optimization problem
\begin{equation*}
\min\limits_{\boldsymbol{v}} \frac{1}{2} \boldsymbol{v}^\top \nabla_{\boldsymbol{y}}^2 g(\boldsymbol{x}^{k+1}, \boldsymbol{y}^{k+1})\boldsymbol{v} - \boldsymbol{v}^\top \nabla_{\boldsymbol{y}} f(\boldsymbol{x}^{k+1}, \boldsymbol{y}^{k+1})
\end{equation*}
using $J$ steps of SGD starting from $\boldsymbol{v}^{k+1, 0}$(see line 11 in Algorithm \ref{alg:2}), where we also adopt a warm start with $\boldsymbol{v}^{k+1, 0}=\boldsymbol{v}^{k}$. Then, by setting $\boldsymbol{v}^{k+1}=\boldsymbol{v}^{k+1, J}$,
\begin{align} \label{eq15}
\boldsymbol{h}_f^{k+1} &= \nabla_{\boldsymbol{x}} F(\boldsymbol{x}^{k+1}, \boldsymbol{y}^{k+1};\mathcal{D}_F) \\
&  \qquad -\nabla_{\boldsymbol{x}} \nabla_{\boldsymbol{y}} G(\boldsymbol{x}^{k+1}, \boldsymbol{y}^{k+1}; \mathcal{D}_G)\boldsymbol{v}^{k+1}  \nonumber
\end{align}
is constructed as an estimate to $\nabla \Phi(\boldsymbol{x}^{k+1})$. Note that sample sets $\{\mathcal{D}_F, \mathcal{D}_G, \mathcal{B}_j(j=0, \ldots, J-1), \mathcal{D}_{F,j}(j=0, \ldots, J-1)\}$ are mutually independent in the hypergradient estimation process.

For the sample sets in Algorithm (\ref{alg:2}), we suppose the sample sets $\mathcal{S}_t$ for all $t=0, \ldots, T-1$, $\mathcal{B}_j$ for all $j=0, \ldots, J-1$, and $\mathcal{D}_G$ have the sizes of $S$, $D$, and $D_g$, respectively. Furthermore, we suppose sample set $\mathcal{D}_F$ and the sample sets $\mathcal{D}_{F, j}$ for all $j=0, \ldots, J-1$ have the same sample size $D_f$.

\begin{remark}
From Algorithm \ref{alg:2}, it is easy to observe that there are three loops. However, since in the classification of the existing stochastic algorithms for solving problem (\ref{eq2}), the loop to estimate the hypergradient, i.e., approximate the Hessian inverse, is typically not counted. To be consistent with the existing literature, when performing classification, the loop to solve linear system in lines 9-12 of Algorithm \ref{alg:2} is not counted. Then, it is obvious that SSGD can be either a single loop algorithm(i.e., $T=1$) or a double loop algorithm(i.e., $T>1$). Although in the analysis that follows, it can be concluded that the convergence of SSGD is guaranteed for any $T \ge 1$. To obtain a lower computational complexity than the existing simple SGD type algorithms, $T=\mathcal{O}(\kappa)$ is required for SSGD, where $\kappa$ is the condition number. Therefore, we call SSGD a double loop algorithm, as shown in Table \ref{table1}.
\end{remark}
\section{Theoretical Results} \label{section:4}
In this section, we provide the convergence and complexity results of SSGD. Detailed theoretical analysis can be found in the supplementary material.

We first show that the convergence of SSGD can be guaranteed for any $T\ge 1$ and $J\ge 1$. Furthermore, the complexity result is provided.

\begin{theorem} \label{tho:1}
Apply SSGD to solve problem (\ref{eq2}). Suppose Assumptions 1, 2, and 3 hold. Define $r_v = {2\mu L}/{(L^2 + \mu^2)} = \mathcal{O}(\kappa^{-1})$, $r_w = {\eta \mu}/{(7(2 - \eta \mu))}$, $\rho_2 = \mathcal{O}(\kappa^{-4})$, $\rho_1 = {(2\rho_2\left((1 + r_w)C_{23} + 2C_{24}\right))}/{(\rho_y (1 + r_v) - r_v)}$, and 
\begin{equation*}
L_{31} = \frac{L_{\Phi}}{2} + 2\left((\rho_1 + C_{24}\rho_2)\frac{2}{r_v} \frac{L^2}{\mu^2} + (1 + r_w)C_{22}\rho_2\right),
\end{equation*}
where $\kappa = {L}/{\mu}$ is the condition number, $\rho_y= 2\beta \mu L/(\mu + L)$, $L_{\Phi} = \mathcal{O}(\kappa^3)$ is defined in Proposition \ref{pro:1}, $C_{22}$
, $C_{23}$, and $C_{24}$ are defined in Lemma 7 in the supplementary material, and $\eta$, $\beta$ are the stepsizes. Let $\|\boldsymbol{v}^0\| \le M/\mu$, and choose stepsizes $\eta = {1}/{(2L)}$, $\beta = {3}/{(2(L + \mu))}=\mathcal{O}(1)$, and 
\begin{equation*}
\alpha = \min\left\{ \frac{1}{2L_{31}},~\frac{1}{4L^2}\rho_2 \eta \mu,~ \left((1 + r_w)C_{23} + 2C_{24}\right)\frac{\rho_2}{\bar{C}_1^2}\right\}, 
\end{equation*}
where $\bar{C}_1 = L + {M \tau}/{\mu} = \mathcal{O}({\kappa})$ is defined in Lemma 4 in the supplementary material. Then, for any $T\ge 1$ and $J\ge 1$, we have 
\begin{align*} 
&\frac{1}{K} \sum_{k=0}^{K-1} \mathbb{E}\| \nabla \Phi(\boldsymbol{x}^k) \|^2 \\
&\le \mathcal{O}\left(\frac{\kappa^5}{K} +\kappa^2 \frac{\sigma_f^2}{D_f} + \kappa^2 \frac{\sigma_{g, 1}^2}{D_g} + \kappa^8 \frac{1}{S}\sigma_g^2 + \kappa^4 \frac{\sigma_{g, 2}^2}{D}\right).
\end{align*}
Furthermore, to achieve an $\epsilon$-accurate stationary point, the computational complexity is $\mathcal{O}(\kappa^{13} \epsilon^{-2}\max\{T, J\})$.
\end{theorem}

Theorem \ref{tho:1} shows that SSGD converges sublinearly w.r.t. the number $K$ of iterations, the batch sizes $D_f$, $S$ for gradient estimation, $D_g$ for Jacobian estimation, and $D$ for Hessian matrix estimation. In addition, from Theorem \ref{tho:1}, we know that when $T=1$ and $J=1$, SSGD can achieve the lowest computational complexity $\mathcal{O}(\kappa^{13}\epsilon^{-2})$. Compared with the existing simple SGD type algorithms listed in Table \ref{table1}, this computational complexity of SSGD outperforms that of BSA, TTSA, and stocBiO by an order of $\mathcal{O}(\epsilon^{-1}\log({1}/{\epsilon}))$, $\mathcal{O}(\epsilon^{-0.5}\log({1}/{\epsilon}))$, and $\mathcal{O}(\log({1}/{\epsilon}))$ in terms of the target accuracy $\epsilon$, respectively. While, in terms of the condition number $\kappa$, it has a higher complexity compared with the algorithms listed in Table \ref{table1}.

Notice that the convergence result and the complexity result of Theorem \ref{tho:1} are established under the assumption that $T\ge 1$ and $J\ge 1$. In the following, we show that by imposing some restrictions on $T$ and $J$, SSGD can obtain a lower computational complexity.

\begin{figure*}[!t]
\centering
\subfloat[]{\includegraphics[scale = 0.2]{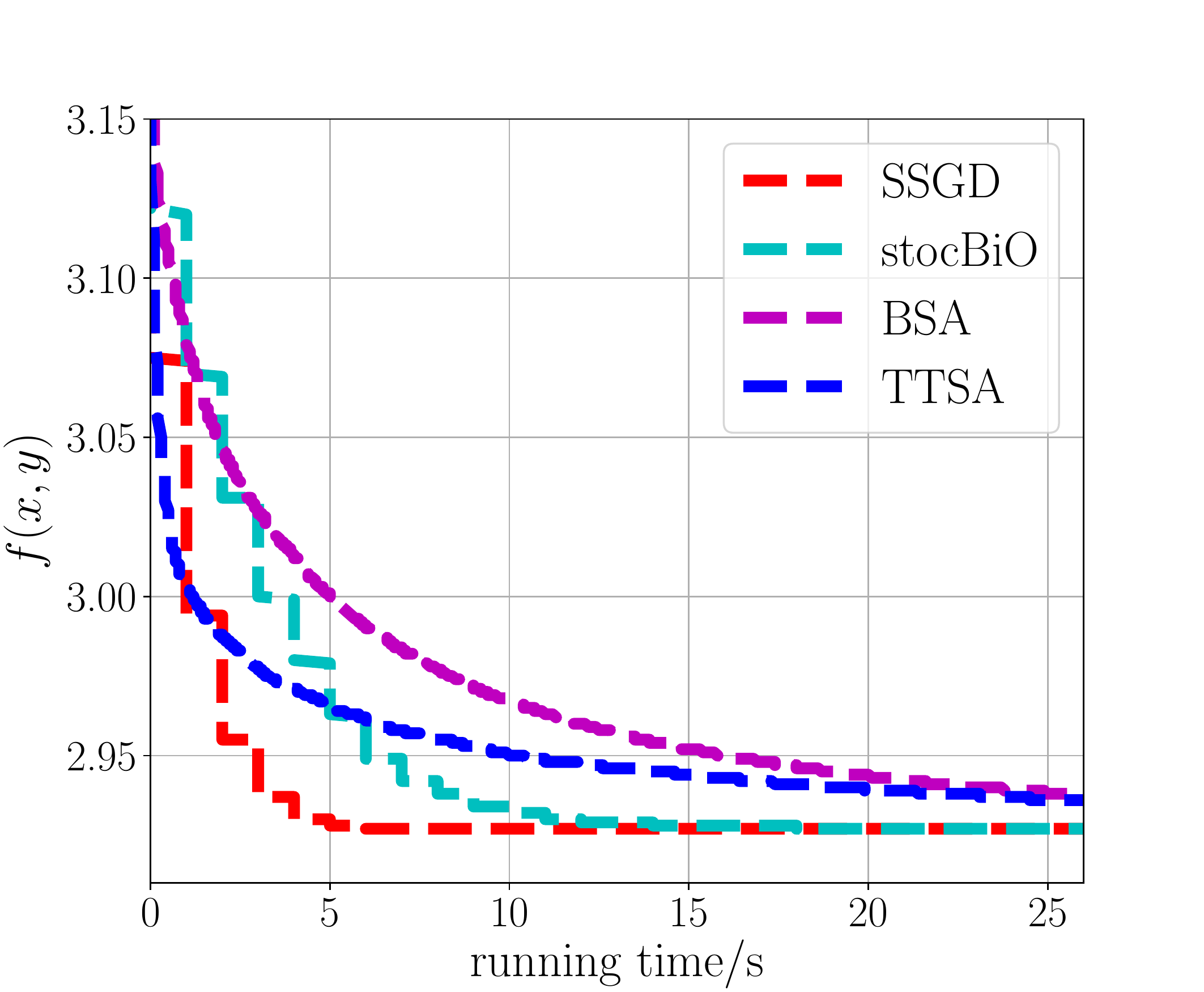}%
\label{fig_first_case}}
\hfil
\subfloat[]{\includegraphics[scale = 0.2]{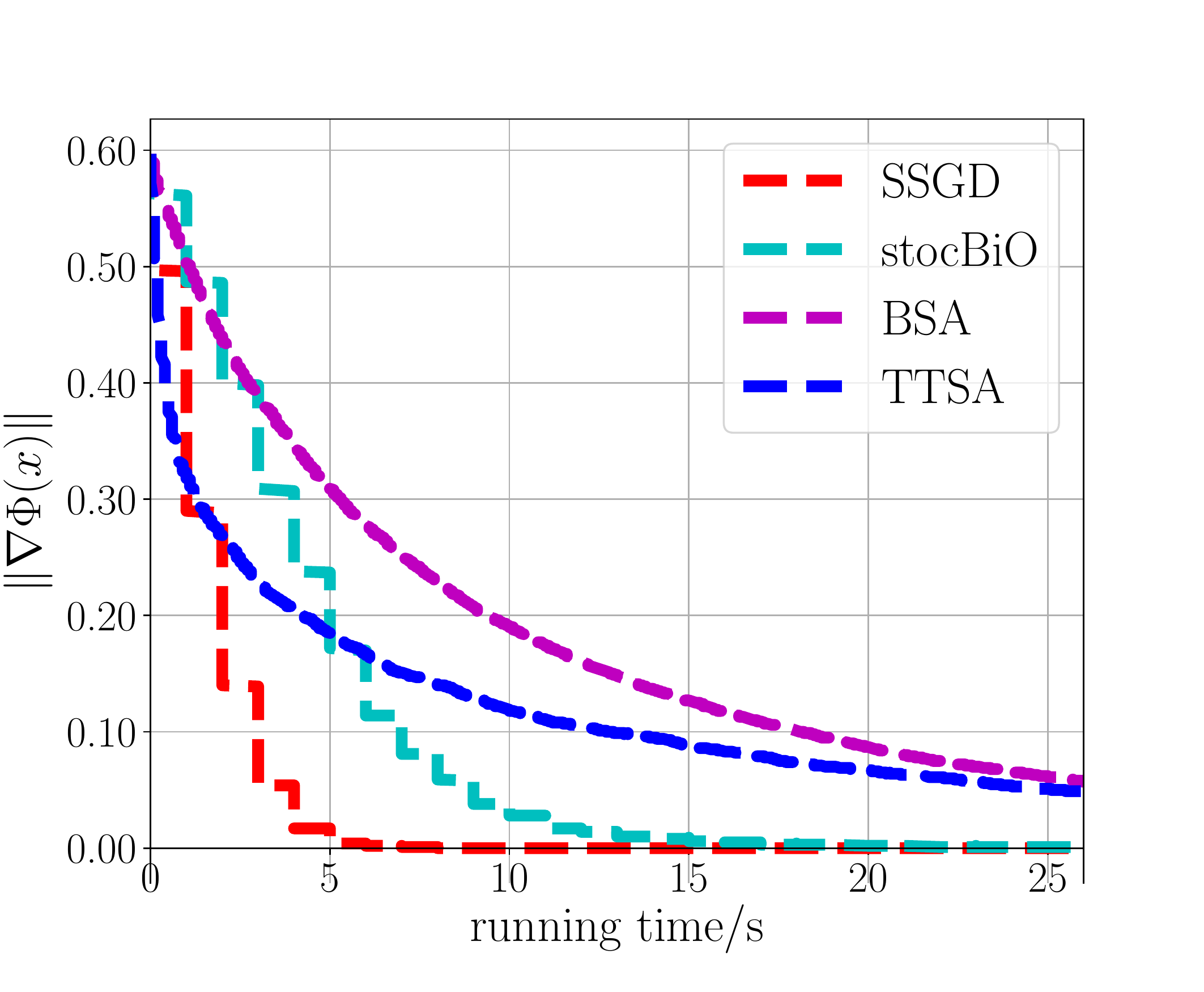}%
\label{fig_second_case}}
\hfil
\subfloat[]{\includegraphics[scale = 0.2]{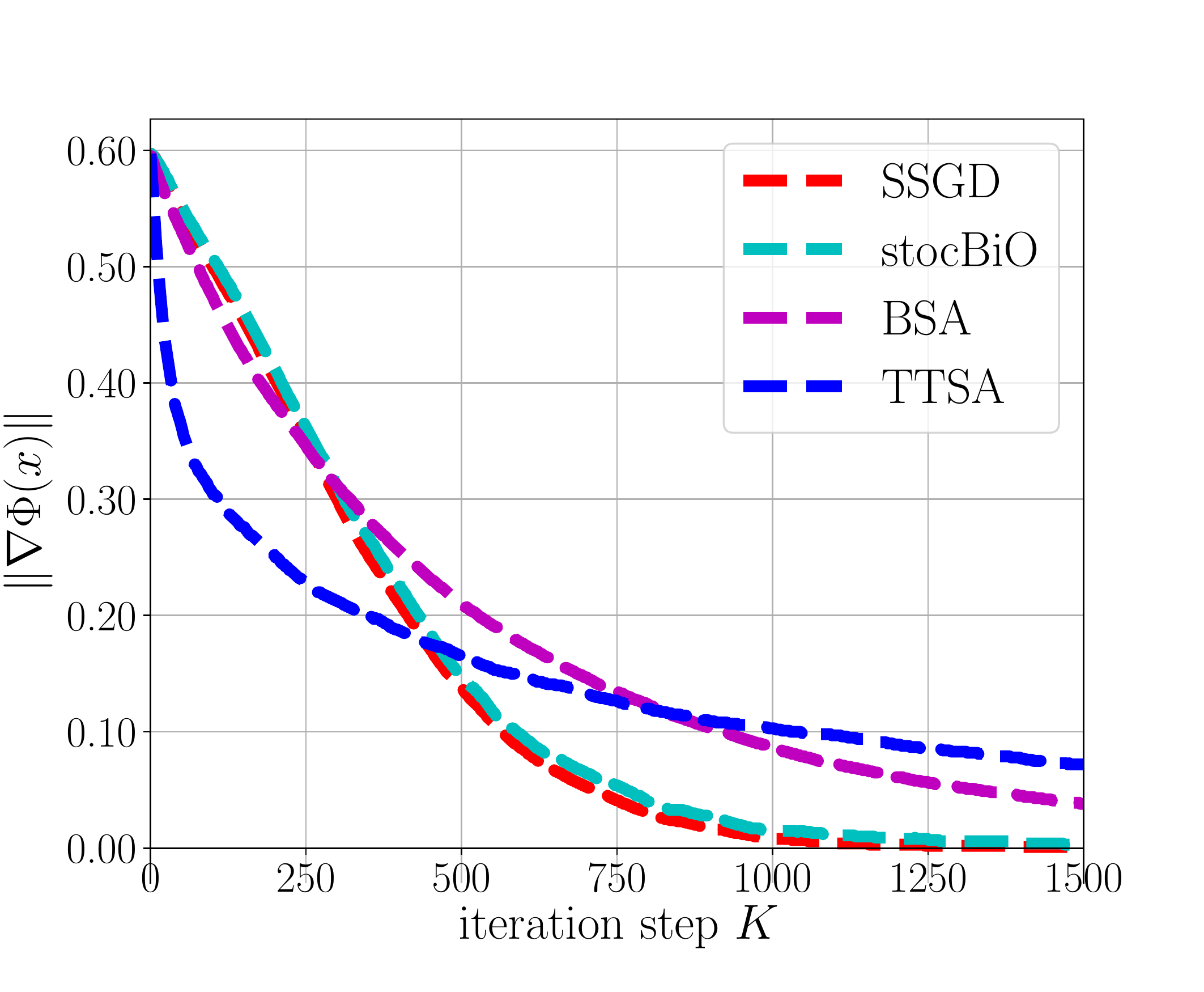}%
\label{fig_third_case}}
\caption{Comparison of the simple SGD type algorithms for solving problem (\ref{eq14}) whose data sets are constructed by setting $\boldsymbol{w}_0 =(2, 5, 7)$.}
\label{fig:2}
\end{figure*}

\begin{figure*}[!t]
\centering
\subfloat[]{\includegraphics[scale = 0.2]{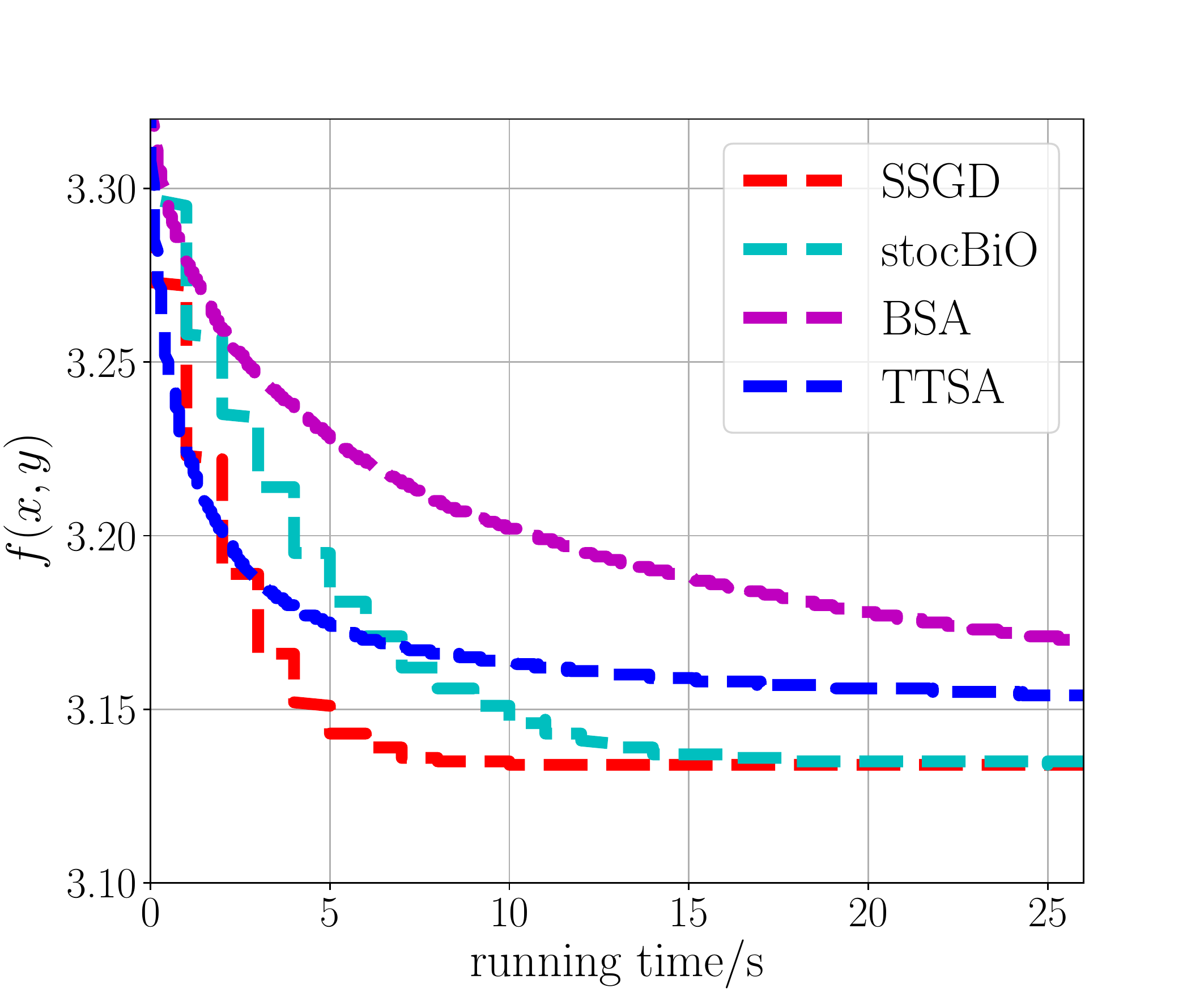}%
\label{fig_first_case}}
\hfil
\subfloat[]{\includegraphics[scale = 0.2]{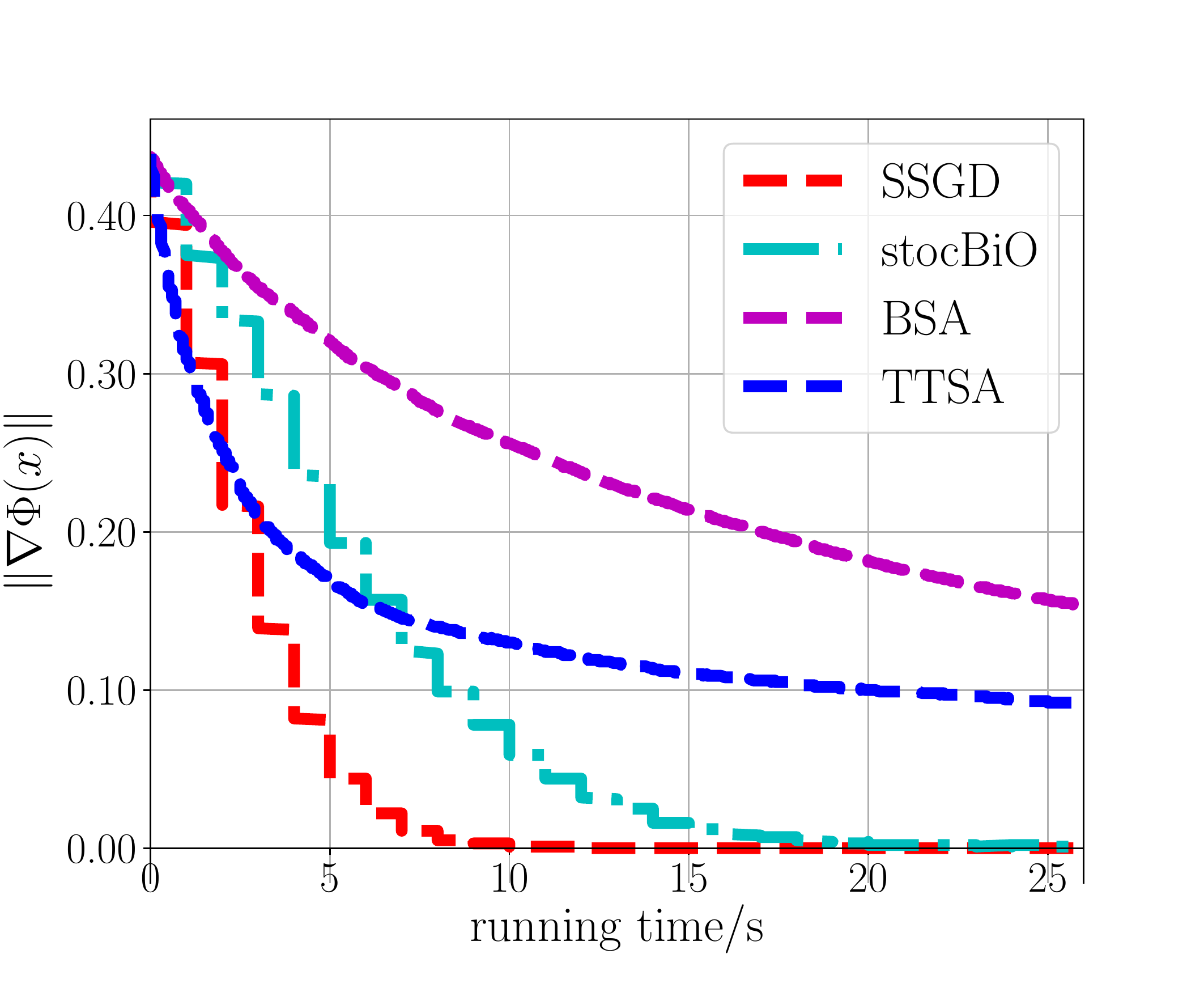}%
\label{fig_second_case}}
\hfil
\subfloat[]{\includegraphics[scale = 0.2]{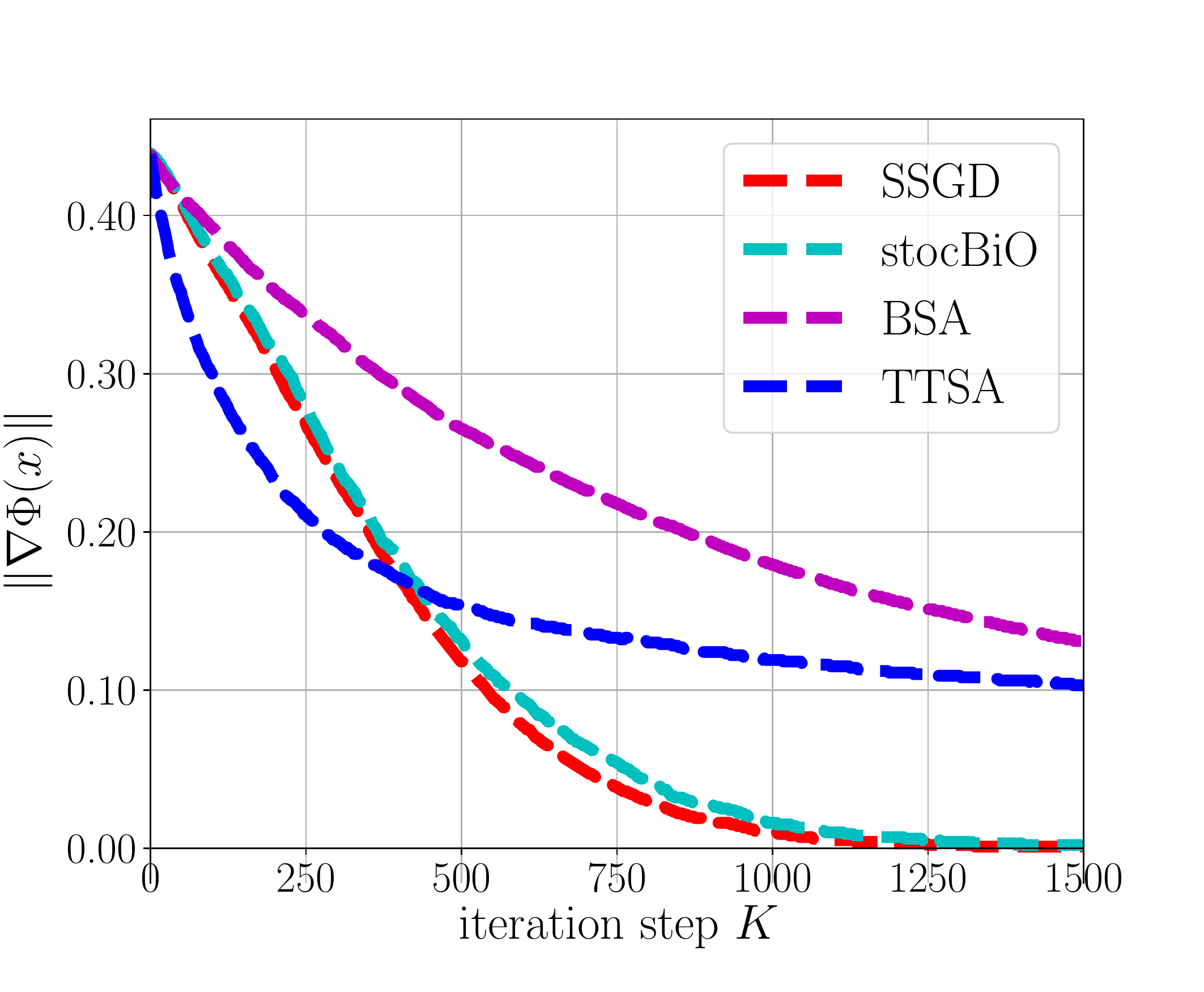}%
\label{fig_third_case}}
\caption{Comparison of the simple SGD type algorithms for solving problem (\ref{eq14}) whose data sets are constructed by setting $\boldsymbol{w}_0 =(1, 4, 6, \ldots, 6)\in \mathbb{R}^{500}$.}
\label{fig:3}
\end{figure*}

\begin{theorem} \label{tho:2}
Apply SSGD to solve problem (\ref{eq2}). Suppose Assumptions 1, 2, and 3 hold. Define $r_v = 1$, $r_w = 1$, $\rho_2 = \mathcal{O}(\kappa^{-3})$, $\rho_1 = 2\rho_2(1+r_w)C_{23}$, and 
\begin{equation*}
L_{31} = \frac{L_{\Phi}}{2} + 2\left(\frac{\rho_1}{4\mu^2}L^2 + \rho_2(1 + r_w)C_{22}\right),
\end{equation*}
where $\kappa = {L}/{\mu}$, $L_{\Phi} = \mathcal{O}(\kappa^3)$ is defined in Proposition \ref{pro:1}, $C_{22} = \left(1 + {1}/{r_w}\right)\left(2 + 8 {L^2}/{\mu^2}\right)C_1^2$, $C_{23} = \left(1 + {1}/{r_w}\right)(5 + r_w)C_1^2$, and $C_1 = {M}/{\mu^2}\rho + {L}/{\mu}$. Let $\|\boldsymbol{v}^0\| \le M/\mu$, and select $\eta = {1}/{(2L)}$, $\beta = {1}/{(2L)}$, 
\begin{equation*}
\alpha = \min\left\{ \frac{1}{2L_{31}},~\frac{1}{4L^2}\rho_2,~ \frac{1}{2}\rho_2(1 + r_w)\frac{C_{23}}{\bar{C}_1^2} \right\}, 
\end{equation*}
\begin{equation*}
T \ge \frac{\log(\rho_1/(8(\rho_1 + C_{24}\rho_2)))}{\log(L/(\mu + L))}=\mathcal{O}(\kappa),  
\end{equation*}
\begin{equation*}
J \ge \frac{\log(1/(4(1 + r_w)^2))}{\log(1-\eta \mu)} =\mathcal{O}(\kappa),
\end{equation*}
where $\bar{C}_1 = L + {M \tau}/{\mu} = \mathcal{O}({\kappa})$, and $C_{24} = \left(1 + {1}/{r_w}\right)(9 + 8r_w )C_1^2$. Then, we have 
\begin{align*} 
&\frac{1}{K} \sum_{k=0}^{K-1} \mathbb{E}\| \nabla \Phi(\boldsymbol{x}^k) \|^2  \\
&\le \mathcal{O}\left(\frac{\kappa^3}{K} +\kappa \frac{\sigma_f^2}{D_f} + \kappa^2 \frac{\sigma_{g, 1}^2}{D_g} + \kappa^5 \frac{1}{S}\sigma_g^2 + \kappa^3 \frac{\sigma_{g, 2}^2}{D}\right).
\end{align*}
Furthermore, a computational complexity of the order of $\mathcal{O}(\kappa^{9} \epsilon^{-2})$ is sufficient for SSGD to reach an $\epsilon$-accurate stationary point.
\end{theorem}

In Theorem \ref{tho:2}, the convergence of SSGD is established under the assumption that $T\ge \mathcal{O}(\kappa)$ and $J\ge \mathcal{O}(\kappa)$. Furthermore, a computational complexity by an order of $\mathcal{O}(\kappa^{9} \epsilon^{-2})$ is obtained in Theorem \ref{tho:2}, which is lower than the lowest computational complexity $\mathcal{O}(\kappa^{13} \epsilon^{-2})$ obtained in Theorem \ref{tho:1}. Compared with the algorithms listed in Table \ref{table1}, SSGD achieves the lowest computational complexity.

\section{Experimental Results}\label{section:5}

In this section, we first validate our theoretical results on synthetic bilevel optimization problems. Then, we apply our proposed algorithm on a hyperparameter optimization problem.

\subsection{Synthetic Bilevel Optimization Problems}
In the following, we perform experiments on the bilevel problem in (\ref{eq14}). The experimental details are in the supplementary material.

Firstly, we compare our proposed algorithm SSGD with the existing simple SGD type algorithms BSA, TTSA, and stocBiO. In the experiments, we consider two datasets, i.e., datasets constructed by setting $\boldsymbol{w}_0 = (2, 5, 7)$ and datasets constructed by setting $\boldsymbol{w}_0 =(1, 4, 6, \ldots, 6)\in \mathbb{R}^{500}$.

From Fig. \ref{fig:2}\subref{fig_third_case} and Fig. \ref{fig:3}\subref{fig_third_case}, we observe that SSGD and stocBiO run almost the same number of iteration steps $K$ to converge to the stationary point with $\|\nabla \Phi(\boldsymbol{x})\| = 0$, less than that of BSA and TTSA. However, compared with SSGD, a larger $J$ to estimate the hypergradient is required for stocBiO to converge to the stationary point with $\|\nabla \Phi(\boldsymbol{x})\| = 0$. As a result, Fig. \ref{fig:2}\subref{fig_second_case} and Fig. \ref{fig:3}\subref{fig_second_case} show that SSGD converges fastest among all the compared methods. Furthermore, from Fig. \ref{fig:2}\subref{fig_first_case} and Fig. \ref{fig:3}\subref{fig_first_case}, we observe that SSGD converges to the point with the smallest upper objective function value for both cases.

Then, we study the influence of iteration step $T$, iteration step $J$, batch size, and stepsize $\alpha$ on the convergence behavior of SSGD in Algorithm \ref{alg:2}. In the experiments, the datasets are constructed by setting $\boldsymbol{w}_0 = (2, 5, 7, \ldots, 7)\in \mathbb{R}^{1000}$.

\begin{figure}[htb]
\centering
\subfloat[]{\includegraphics[scale = 0.2]{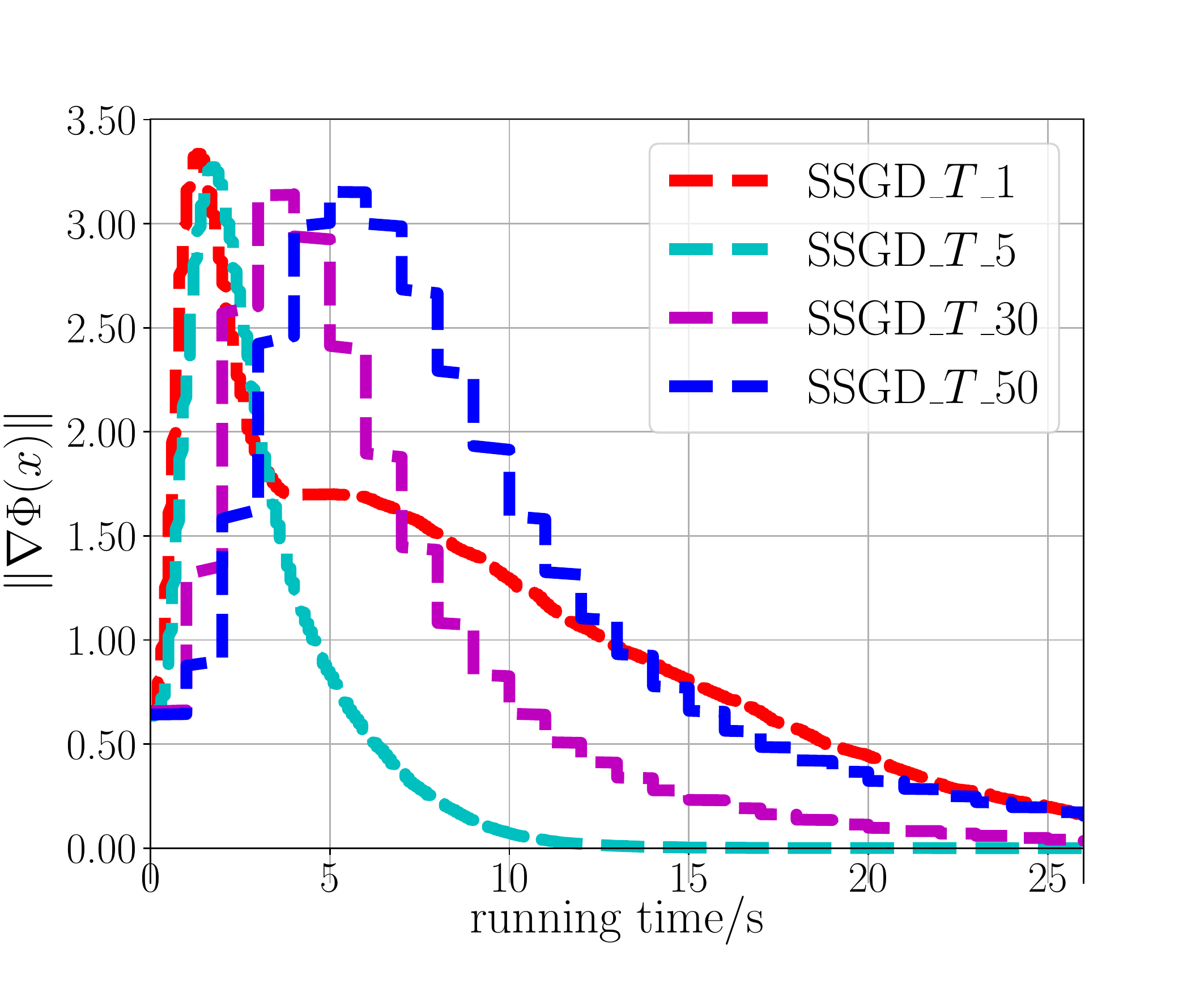}%
\label{fig_first_case}}
\subfloat[]{\includegraphics[scale = 0.2]{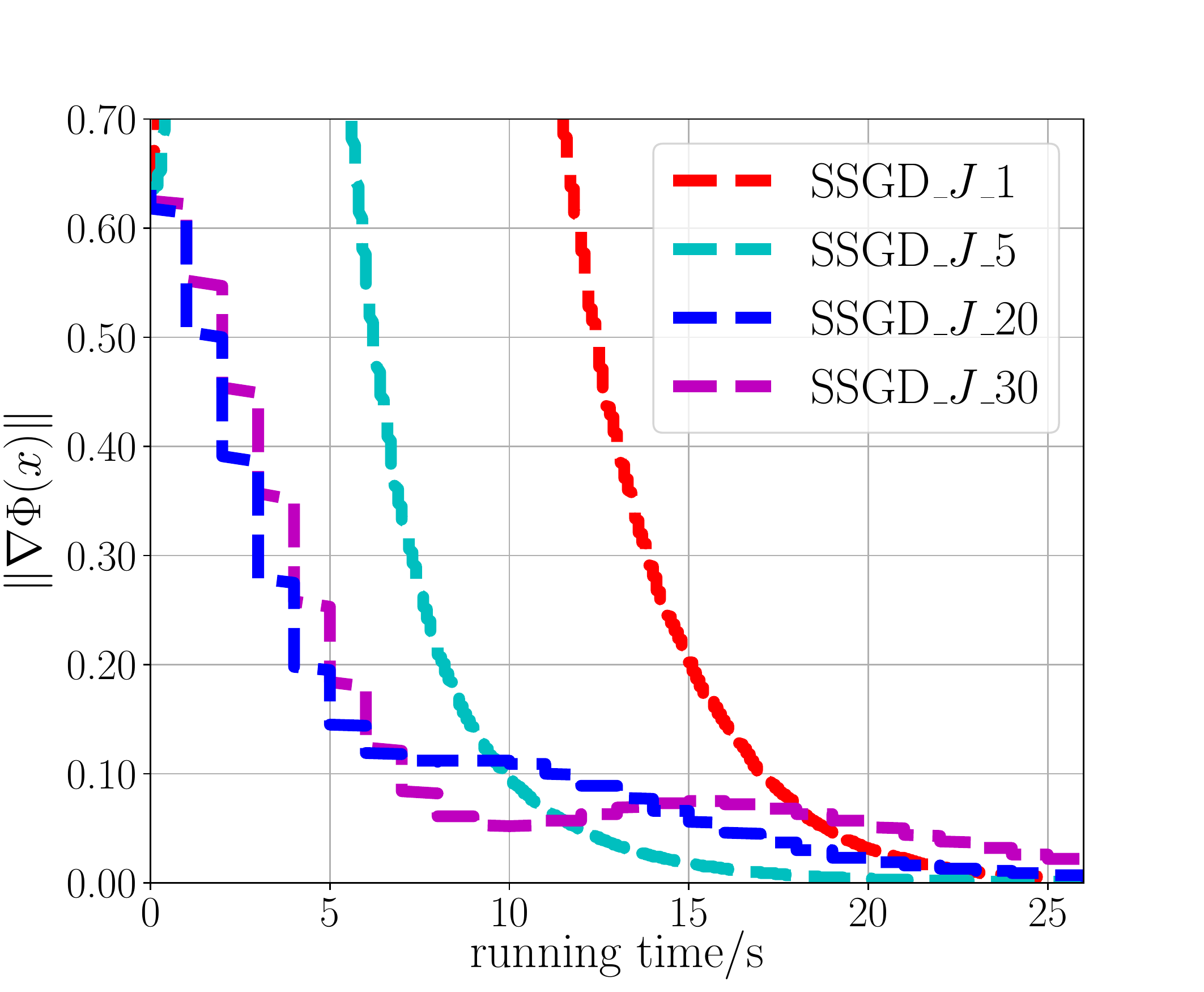}%
\label{fig_second_case}}
\caption{Illustrating the convergence behavior of SSGD under different iteration steps $T$ and iteration steps $J$. (a) shows the converence behavior under different $T$. (b) shows the convergence behavior under different $J$.}
\label{fig:4}
\end{figure}

Fig. \ref{fig:4}\subref{fig_first_case} shows the convergence behavior of SSGD under different choices of $T\in \{1, 5, 30, 50\}$, where SSGD\_$T$\_$i$ indicates that $T$ is set to be $i$. It shows that beginning from 1, increasing $T$ can speed up the convergence of SSGD, but when $T$ reaches $5$, increasing $T$ can slow down the convergence speed of SSGD. Fig. \ref{fig:4}\subref{fig_second_case} shows the convergence behavior of SSGD under different choices of $J\in \{1, 5, 20, 30\}$, where SSGD\_$J$\_$i$ indicates that $J$ is set to be $i$. We observe that SSGD converges fastest when $J$ is chosen to be 5.

\begin{figure}[htb]
\centering
\subfloat[]{\includegraphics[scale = 0.2]{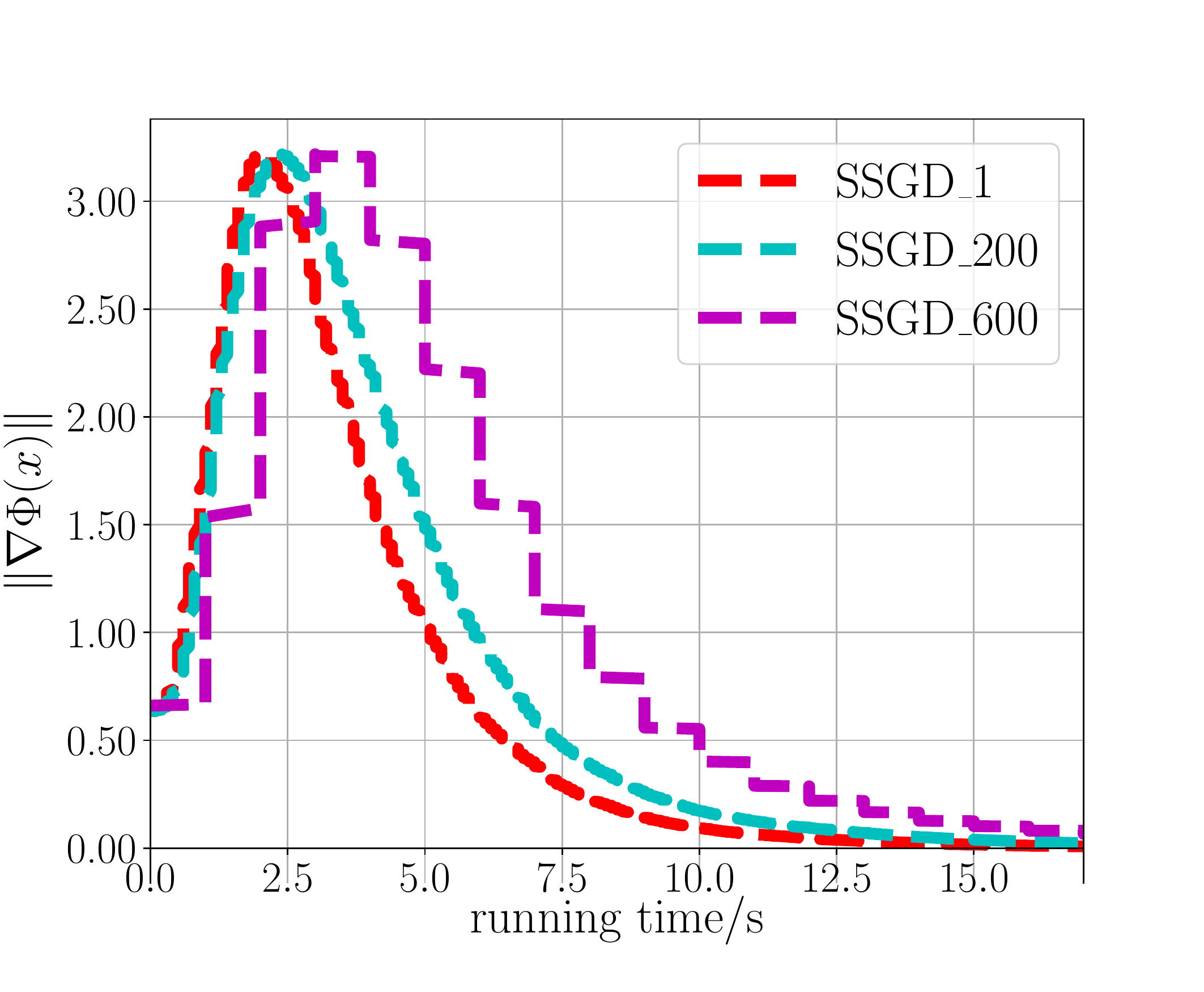}%
\label{fig_first_case}}
\hfil
\subfloat[]{\includegraphics[scale = 0.2]{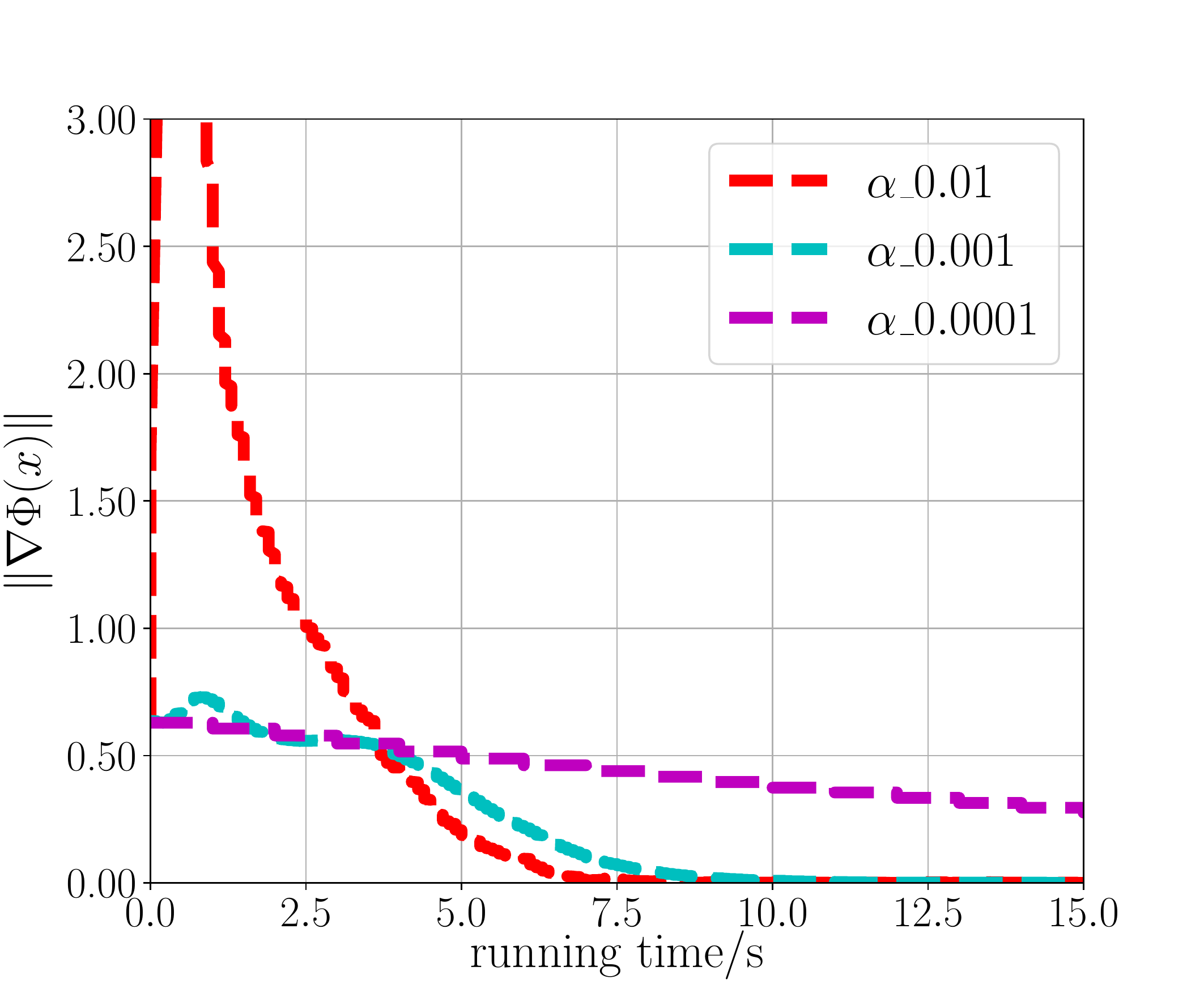}%
\label{fig_second_case}}
\caption{Illustrating the convergence behavior of SSGD under different batch sizes and stepsizes $\alpha$. (a) shows the convergence behavior under different batch sizes. (b) shows the convergence behavior under different stepsizes $\alpha$.}
\label{fig:5}
\end{figure}
\begin{figure}[htb]
\centering
\subfloat[]{\includegraphics[scale = 0.2]{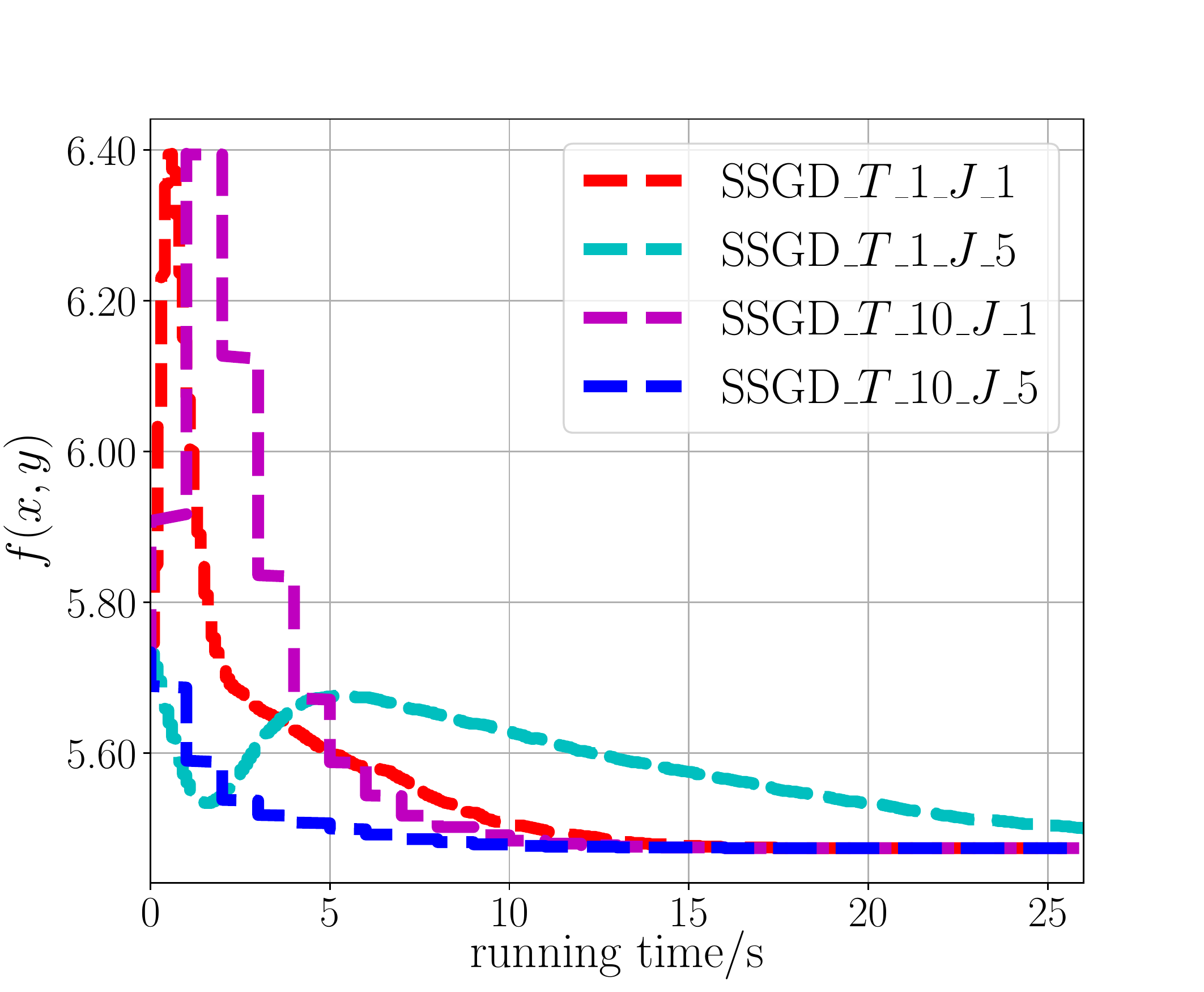}%
\label{fig_first_case}}
\hfil
\subfloat[]{\includegraphics[scale = 0.2]{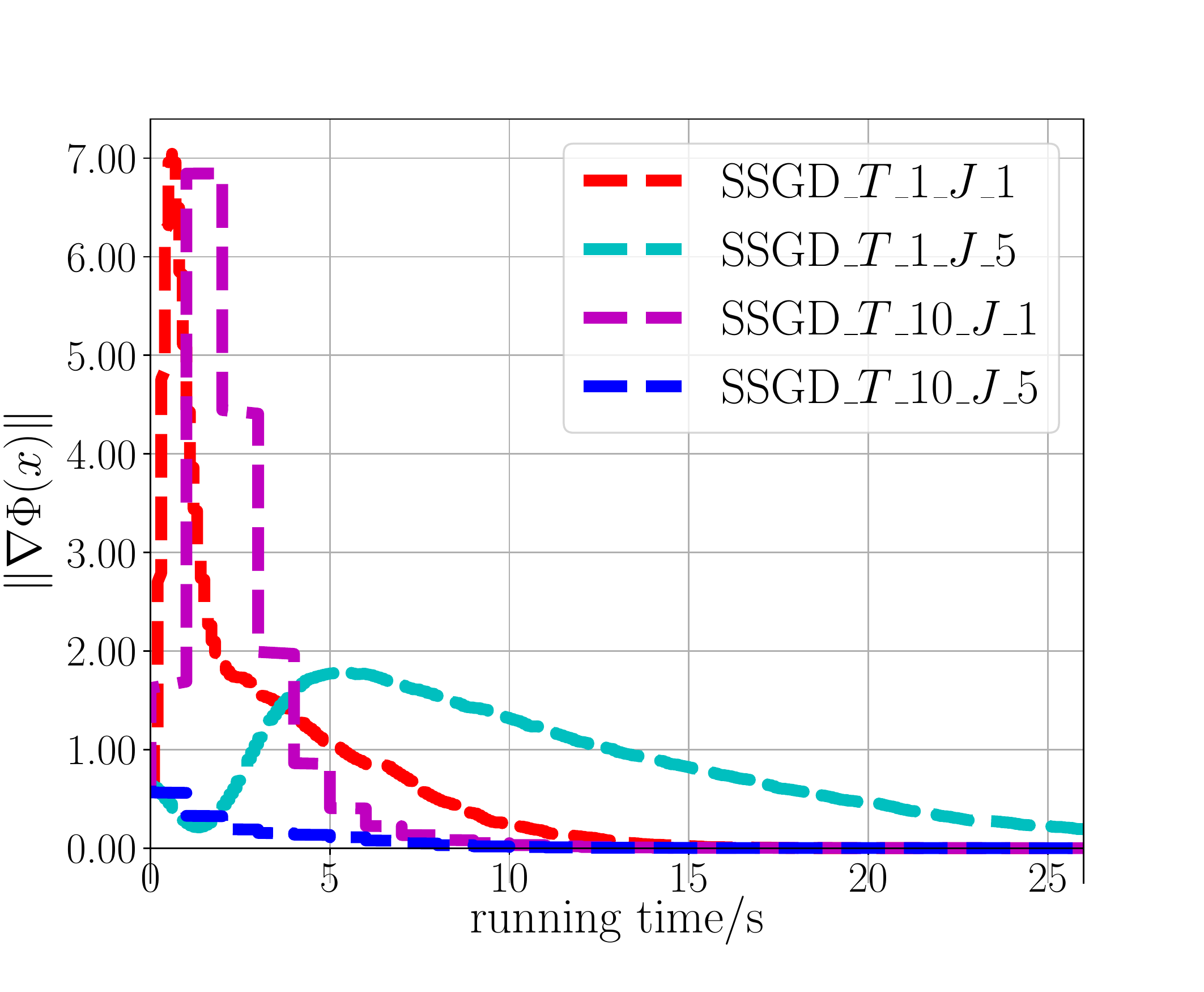}%
\label{fig_second_case}}
\caption{Illustrating the convergence behavior of SSGD under different choices of $T$ and $J$.}
\label{fig:6}
\end{figure}

Fig. \ref{fig:5}\subref{fig_first_case} shows the convergence behavior of SSGD under different batch sizes with SSGD\_$i$ indicating that the batch size is set to be $i$. It shows that SSGD converges fastest when batch size is set to be 1, and increasing batch size can slow down the convergence speed of SSGD. Fig. \ref{fig:5}\subref{fig_second_case} shows the convergence behavior of SSGD under different stepsizes $\alpha$. It shows that choosing a larger stepsize $\alpha$ allows SSGD to converge faster.

Fig. \ref{fig:6} shows the convergence behavior of SSGD under different choices of $T$ and $J$, where SSGD\_$T$\_$t$\_$J$\_$j$ indicates that $T$ is set to be $t$ and $J$ is set to be $j$. It shows that SSGD converges faster in the double loop case(i.e., $T>1$) than in the single loop case(i.e., $T=1$).

\subsection{Hyper-parameter Optimization}
In this section, we use a hyperparameter optimization problem, i.e., data hyper-cleaning\cite{5}, to evaluate the performance of our algorithm. Data hyper-cleaning is to train a classifier under the case that some data in the training set is corrupted. Following \cite{5}, when we train a linear classifier, data hyper-cleaning task can be expressed as a bilvel optimization problem whose UL function is given by 
\begin{equation*}
f(\boldsymbol{x}, \boldsymbol{y}) = \frac{1}{|\mathcal{D}_{\text{val}}|}\sum\limits_{(\boldsymbol{u}_i, v_i)\in \mathcal{D}_{\text{val}}}L({\boldsymbol{y}}^\top \boldsymbol{u}_i, v_i)
\end{equation*}
and whose LL function is given by 
\begin{equation*}
g(\boldsymbol{x}, \boldsymbol{y}) = \frac{1}{|\mathcal{D}_{\text{tr}}|}\sum\limits_{(\boldsymbol{u}_i, v_i)\in \mathcal{D}_{\text{tr}}}\sigma(x_i)L(\boldsymbol{y}^\top \boldsymbol{u}_i, v_i) + c\|\boldsymbol{y}\|^2
\end{equation*}
where $L$ denotes the cross-entropy loss, $\sigma(\cdot)$ is the sigmoid function, $c$ is the coefficient of the regular term and is set to be $0.001$, and $\mathcal{D}_{\text{val}}$ and $\mathcal{D}_{\text{tr}}$ denote the validation set and training set, respectively.

\begin{table}[htb]   
\caption{Comparison of the simple SGD type algorithms for data hyper-cleaning tasks on MNIST, FashionMNIST and CIFAR10 datasets. Acc.(\%) and F1 score(\%) denote the test accuracy and the harmonic mean of the precision and recall, respectively}
\label{table2}
\centering
\begin{tabular}{ccccccc}\hline
\multicolumn{1}{c}{\multirow{2}{*}{Method}} & \multicolumn{2}{c}{MNIST} & \multicolumn{2}{c}{FashionMNIST} & \multicolumn{2}{c}{CIFAR10} \\ 
 & Acc. & F1 score &Acc. & F1 score & Acc. & F1 score \\  \hline
stocBiO & 89.90 & 63.92 & 82.18 & 65.17 & \textbf{38.36} & 64.07 \\
BSA & 80.67 & 14.20 & 72.57 & 16.51& 22.72 &  21.09\\ 
TTSA & 88.88 & 33.94 & 80.62 & 32.54 & 35.69 & 40.03\\ 
SSGD & \textbf{90.32} & \textbf{83.81} & \textbf{82.32} & \textbf{80.45} & 38.35&\textbf{66.16} \\ \hline
\end{tabular}
\end{table}

In the following experiments, we compare our proposed algorithm SSGD with BSA, TTSA, stocBiO on three datasets: MNIST \cite{28}, FashionMNIST \cite{27} and CIFAR10. Each sample in the training set is corrupted with probability $0.3$. More results and details are in the supplementary material. 

Table \ref{table2} shows the experimental results on three datasets over 270s. We observe that SSGD achieves the highest test accuracy on MNIST and FashionMNIST datasets, and achieves the highest F1 score on all the datasets. It shows that SSGD can obtain competetive results on data hyper-cleaning problem.

\section{Conclusion} \label{section:6}
In the paper, we propose to use SGD-based Estimation to estimate the hypergradient for the nonconvex-strongly-convex bilevel optimization problems, and propose a novel simple SGD type algorithm. We provide the convergence guarantee for our algorithm, and prove that our proposed algorithm can achieve lower computational complexity compared with all the existing simple SGD type algorithms for bilevel optimization. Moreover, the experiments validate our theoretical results.

\section*{Acknowledgments}
The work is partially supported by National Key R\&D Program of China(2020YFB1313503 and 2018AAA0100300), National Natural Science Foundation of China(Nos. 61922019, 61733002, 61672125, and 61976041), and LiaoNing Revitalization Talents Program(XLYC1807088).

\bibliographystyle{IEEEtran}
\bibliography{IEEEabrv, mybible}

\begin{IEEEbiographynophoto}{Haimei Huo}
received the B.S. degree in information and computing science from Changchun University of Science and Technology, China, in 2017. She is currently pursuing the PhD degree in computational mathematics from Dalian University of Technology, Dalian, China. Her research interests include machine learning and optimization.
\end{IEEEbiographynophoto}
\begin{IEEEbiographynophoto}{Risheng Liu}
received his B.Sc. (2007) and Ph.D. (2012) from Dalian University of Technology, China. From 2010 to 2012, he was doing research as joint Ph.D. in robotics institute at Carnegie Mellon University. From 2016 to 2018, he was doing research as Hong Kong Scholar at the Hong Kong Polytechnic University. He is currently a full professor with the Digital Media Department at International School of Information Science \& Engineering, Dalian University of Technology. He was awarded the "Outstanding Youth Science Foundation" of the National Natural Science Foundation of China. His research interests include optimization, computer vision, and multimedia.
\end{IEEEbiographynophoto}

\begin{IEEEbiographynophoto}{Zhixun Su}
received the B.S. degree in mathematics from Jilin University, Changchun, China in 1987, the M.S. degree in computer science from Nankai University, Tianjin, China in 1990, and the Ph.D. degree in computational mathematics from Dalian University of Technology, Dalian, China in 1993. He has been a Professor in the School of Mathematical Sciences, Dalian University of Technology, since 1999. He is currently the director of the Key Laboratory for Computational Mathematics and Data Intelligence of Liaoning Province and the vise director of Liaoning Center for Applied Mathematics. His current research interests include computer graphics, computer vision, computational geometry, and machine learning.
\end{IEEEbiographynophoto}

\vfill

\end{document}


\title{Supplementary Material}
\maketitle

\appendices
\section{Proof of Theorem 1 and Theorem 2}\label{appendix:1}
Throughout the proof, we define 
\begin{align*}
&\mathcal{F}_{k, t} = \sigma\{\boldsymbol{x}^0, \boldsymbol{y}^0, \boldsymbol{v}^0, \ldots, \boldsymbol{x}^k, \boldsymbol{y}^{k,1}, \ldots, \boldsymbol{y}^{k,t}\}, \\
&\mathcal{F}_{k}^{'} = \sigma\{ \boldsymbol{x}^0, \boldsymbol{y}^0, \boldsymbol{v}^0, \ldots, \boldsymbol{x}^k, \boldsymbol{y}^{k,1}, \ldots, \boldsymbol{y}^{k, t}, \ldots, \boldsymbol{y}^k\}, \\
&\mathcal{F}_{k, j} = \sigma\{ \boldsymbol{x}^0, \boldsymbol{y}^0, \boldsymbol{v}^0, \ldots, \boldsymbol{x}^k, \boldsymbol{y}^{k,1}, \ldots, \boldsymbol{y}^k, \boldsymbol{v}^{k,1}, \ldots, \boldsymbol{v}^{k,j}\},
\end{align*}
where $\sigma\{\cdot\}$ denotes the $\sigma$-algebra generated by random variables.

\subsection{Proof of Proposition 1}
In Proposition 1, the differentiablity of $\Phi(\boldsymbol{x})$ and the form of $\nabla \Phi(\boldsymbol{x})$ can be obtained from Proposition 1 in \cite{13}. Furthermore, the Lipschitz continuity of $\boldsymbol{y}^*(\boldsymbol{x})$ and $\nabla \Phi(\boldsymbol{x})$ can be obtained from Lemma 2.2 in \cite{16} and Lemma 2 in the supplementary material of \cite{13}, respectively.

\subsection{Proof of Supporting Lemmas}

\begin{lemma}  \label{lemma:1}
For problem (2), we define $\boldsymbol{v}(\boldsymbol{x}, \boldsymbol{y}) = [\nabla_{\boldsymbol{y}}^2 g(\boldsymbol{x}, \boldsymbol{y})]^{-1} \nabla_{\boldsymbol{y}} f(\boldsymbol{x}, \boldsymbol{y})$. Suppose Assumptions 1 and 2 hold. Then, for $\boldsymbol{x}_1$, $\boldsymbol{x}_2 \in \mathbb{R}^p$, $\boldsymbol{y}_1$, $\boldsymbol{y}_2 \in \mathbb{R}^q$, we have 
\begin{equation*}
\| \boldsymbol{v}(\boldsymbol{x}_1, \boldsymbol{y}_1) - \boldsymbol{v}(\boldsymbol{x}_2, \boldsymbol{y}_2)\| \le C_1(\|\boldsymbol{x}_1 - \boldsymbol{x}_2 \| + \|\boldsymbol{y}_1 - \boldsymbol{y}_2 \|),
\end{equation*}
where $ C_1 = {M \rho}/{\mu^2} + {L}/{\mu}$.
\end{lemma}

\begin{IEEEproof}
Based on the definition of $\boldsymbol{v}(\boldsymbol{x}, \boldsymbol{y})$, we have 
\begin{align*}
&\| \boldsymbol{v}(\boldsymbol{x}_1, \boldsymbol{y}_1) - \boldsymbol{v}(\boldsymbol{x}_2, \boldsymbol{y}_2)\| \\     \nonumber
&= \left\| [\nabla_{\boldsymbol{y}}^2 g(\boldsymbol{x}_1, \boldsymbol{y}_1)]^{-1} \nabla_{\boldsymbol{y}} f(\boldsymbol{x}_1, \boldsymbol{y}_1) - [\nabla_{\boldsymbol{y}}^2 g(\boldsymbol{x}_2, \boldsymbol{y}_2)]^{-1} \nabla_{\boldsymbol{y}} f(\boldsymbol{x}_2, \boldsymbol{y}_2) \right\| \\   \nonumber
& \overset{(i)}{\le}\left\| \left\{[\nabla_{\boldsymbol{y}}^2 g(\boldsymbol{x}_1, \boldsymbol{y}_1)]^{-1} - [\nabla_{\boldsymbol{y}}^2 g(\boldsymbol{x}_2, \boldsymbol{y}_2)]^{-1}\right\} \nabla_{\boldsymbol{y}} f(\boldsymbol{x}_1, \boldsymbol{y}_1)\right\| \\   \nonumber
& \qquad + \left\| [\nabla_{\boldsymbol{y}}^2 g(\boldsymbol{x}_2, \boldsymbol{y}_2)]^{-1}\left(\nabla_{\boldsymbol{y}} f(\boldsymbol{x}_1, \boldsymbol{y}_1) - \nabla_{\boldsymbol{y}} f(\boldsymbol{x}_2, \boldsymbol{y}_2)\right)\right\|\\   \nonumber
& \overset{(ii)}{\le} M \left\| [\nabla_{\boldsymbol{y}}^2 g(\boldsymbol{x}_1, \boldsymbol{y}_1)]^{-1}\right\| \left\| [\nabla_{\boldsymbol{y}}^2 g(\boldsymbol{x}_2, \boldsymbol{y}_2)]^{-1}\right\| \left\| \nabla_{\boldsymbol{y}}^2 g(\boldsymbol{x}_1, \boldsymbol{y}_1) - \nabla_{\boldsymbol{y}}^2 g(\boldsymbol{x}_2, \boldsymbol{y}_2)\right\| \\   \nonumber
&\qquad + \frac{1}{\mu}\left\| \nabla_{\boldsymbol{y}} f(\boldsymbol{x}_1, \boldsymbol{y}_1) - \nabla_{\boldsymbol{y}} f(\boldsymbol{x}_2, \boldsymbol{y}_2)\right\|\\   \nonumber
&\overset{(iii)}{\le}\frac{M}{\mu^2}\rho(\left\| \boldsymbol{x}_1 - \boldsymbol{x}_2\right\| + \left\| \boldsymbol{y}_1 - \boldsymbol{y}_2\right\|) +\frac{L}{\mu}(\left\|\boldsymbol{x}_1 - \boldsymbol{x}_2\right\| + \left\|\boldsymbol{y}_1 -\boldsymbol{y}_2\right\|)\\   \nonumber
&= \left(\frac{M}{\mu^2}\rho + \frac{L}{\mu}\right) (\left\|\boldsymbol{x}_1 - \boldsymbol{x}_2\right\| + \left\|\boldsymbol{y}_1 - \boldsymbol{y}_2\right\|), \nonumber
\end{align*}
where $(i)$ uses the triangle inequality, $(ii)$ follows from Cauchy-Schwartz inequality, Assumptions 1 and 2, and the fact that $\|\mathbf{A}^{-1} - \mathbf{B}^{-1}\| = \|\mathbf{A}^{-1}(\mathbf{A}-\mathbf{B})\mathbf{B}^{-1}\|$ for any two matrices $\mathbf{A}$, $\mathbf{B}$, and $(iii)$ follows from Assumptions 1 and 2. Then, the proof is complete.
\end{IEEEproof}

In the following lemma, we prove the boundness of the sequence $\{\boldsymbol{v}^{k,j}\}$ generated by Algorithm 2.
\begin{lemma} \label{lemma:2}
Apply Algorithm 2 to solve problem (2). Suppose Assumptions 1 and 2 hold. Let $\|\boldsymbol{v}^0 \| \le {M}/{\mu}$, $\eta \le {1}/{L}$, then $\| \boldsymbol{v}^{k, j}\| \le {M}/{\mu}$ for any $k\in \{1, \ldots, K\}$ and $j\in \{ 0, \ldots, J\}$.
\end{lemma}

\begin{IEEEproof}
We first prove that, for any $k\in \{1, \ldots, K\}$, if $\|\boldsymbol{v}^{k-1}\| \le {M}/{\mu}$, then $\|\boldsymbol{v}^{k, j}\| \le {M}/{\mu}$ for $j=0, \ldots, J$. For any $k\in \{1, \ldots, K\}$, suppose $\|\boldsymbol{v}^{k-1}\| \le {M}/{\mu}$. Since $\boldsymbol{v}^{k, 0}=\boldsymbol{v}^{k-1}$(see line 8 in Algorithm 2), we have $\|\boldsymbol{v}^{k,0}\| \le {M}/{\mu}$. By using induction, suppose $\|\boldsymbol{v}^{k, j-1}\| \le {M}/{\mu}$ with $j\in \{1, \ldots, J\}$. Then, since $ \boldsymbol{v}^{k, j} = (\mathbf{I} - \eta \nabla_{\boldsymbol{y}}^2 G(\boldsymbol{x}^k, \boldsymbol{y}^k; \mathcal{B}_{j-1}))\boldsymbol{v}^{k, j-1} + \eta \nabla_{\boldsymbol{y}} F(\boldsymbol{x}^k, \boldsymbol{y}^k; \mathcal{D}_{F, j-1})$(see line 11 in Algorithm 2), we have 
\begin{align*}
\|\boldsymbol{v}^{k, j}\| &\overset{(i)}{\le} \|  (\mathbf{I} - \eta \nabla_{\boldsymbol{y}}^2 G(\boldsymbol{x}^k, \boldsymbol{y}^k; \mathcal{B}_{j-1}))\| \| \boldsymbol{v}^{k, j-1}\| \\   \nonumber
& \qquad + \eta \| \nabla_{\boldsymbol{y}} F(\boldsymbol{x}^k, \boldsymbol{y}^k; \mathcal{D}_{F, j-1})\|\\   \nonumber
&\overset{(ii)}{\le} (1 - \eta \mu)\| \boldsymbol{v}^{k, j-1}\| + \eta M \le \frac{M}{\mu},  \nonumber
\end{align*}
where $(i)$ uses the triangle inequality and Cauchy-Schwartz inequality, and $(ii)$ follows from Assumptions 1 and 2. It is proved.

Next, we prove $\| \boldsymbol{v}^k \| \le {M}/{\mu}$ for any $k\in \{0, \ldots, K-1\}$. From the conditions of Lemma $\ref{lemma:2}$, we have $\|\boldsymbol{v}^0 \| \le {M}/{\mu}$. By using induction, suppose $\|\boldsymbol{v}^{k-1}\| \le {M}/{\mu}$ with $k\in \{1, \ldots, K-1\}$. Then, based on the above discussion, we have $\|\boldsymbol{v}^{k, J}\| \le {M}/{\mu}$, and $\|\boldsymbol{v}^k\|=\|\boldsymbol{v}^{k, J} \| \le {M}/{\mu}$(see line 13 in Algorithm 2). It is proved.

Finally, based on the above two conclusions, the proof is complete.
\end{IEEEproof}

To go further, we first give some formulas. From the hypergradient formula in (3), the hypergradient at $\boldsymbol{x}^k$ can be expressed as follows:
\begin{equation}  \tag{17}\label{eq17}
\nabla \Phi(\boldsymbol{x}^k) = \nabla_{\boldsymbol{x}} f(\boldsymbol{x}^k, \boldsymbol{y}^*(\boldsymbol{x}^k)) - \nabla_{\boldsymbol{x}} \nabla_{\boldsymbol{y}} g(\boldsymbol{x}^k, \boldsymbol{y}^*(\boldsymbol{x}^k))\boldsymbol{v}(\boldsymbol{x}^k, \boldsymbol{y}^*(\boldsymbol{x}^k)),
\end{equation}
where $\boldsymbol{v}(\boldsymbol{x}^k, \boldsymbol{y}^*(\boldsymbol{x}^k)) = [\nabla_{\boldsymbol{y}}^2 g(\boldsymbol{x}^k, \boldsymbol{y}^*(\boldsymbol{x}^k))]^{-1}\nabla_{\boldsymbol{y}} f(\boldsymbol{x}^k, \boldsymbol{y}^*(\boldsymbol{x}^k))$. Furthermore, from the definition of $\boldsymbol{h}_f^{k+1}$ in (16), $\boldsymbol{h}_f^k$ in line 2 of Algorithm 2 has the following form:
\begin{equation} \tag{18}  \label{eq18}
\boldsymbol{h}_f^{k} = \nabla_{\boldsymbol{x}} F(\boldsymbol{x}^{k}, \boldsymbol{y}^{k};\mathcal{D}_F) 
-\nabla_{\boldsymbol{x}} \nabla_{\boldsymbol{y}} G(\boldsymbol{x}^{k}, \boldsymbol{y}^{k}; \mathcal{D}_G)\boldsymbol{v}^{k}.
\end{equation}

To prove the convergence of SSGD, we first analyze the variation of $\Phi(\boldsymbol{x})$ after the variable $\boldsymbol{x}$ is updated once.

\begin{lemma} \label{lemma:3}
Suppose Assumptions 1 and 2 hold. Let $\bar{\boldsymbol{h}}_f^k = \mathbb{E}[\boldsymbol{h}_f^k| \mathcal{F}_k^{'}]$, we have
\begin{align} \tag{19}  \label{eq19}
\mathbb{E}[\Phi(\boldsymbol{x}^{k+1})] &\le \mathbb{E}[\Phi(\boldsymbol{x}^k)] - \frac{\alpha}{2}\mathbb{E}\| \nabla \Phi(\boldsymbol{x}^k) \|^2 - \left(\frac{\alpha}{2} - \frac{L_{\Phi}}{2}\alpha^2\right)\mathbb{E}\| \bar{\boldsymbol{h}}_f^k\|^2 \\
&\qquad  + \frac{\alpha}{2}\mathbb{E}\| \nabla \Phi(\boldsymbol{x}^k) - \bar{\boldsymbol{h}}_f^k \|^2  + \frac{L_{\Phi}}{2}\alpha^2 \mathbb{E}\| \boldsymbol{h}_f^k - \bar{\boldsymbol{h}}_f^k \|^2 . \nonumber
\end{align}
\end{lemma}

\begin{IEEEproof}
 Using the Lipschitz property of $\nabla \Phi(\boldsymbol{x})$ in Proposition 1, and taking the conditional expectation on $\mathcal{F}_k^{'}$, we have
\begin{align*}
\mathbb{E}[\Phi(\boldsymbol{x}^{k+1})| \mathcal{F}_k^{'}] &\le \Phi(\boldsymbol{x}^k) + \mathbb{E}[ \langle \nabla \Phi(\boldsymbol{x}^k), \boldsymbol{x}^{k+1} - \boldsymbol{x}^k \rangle | \mathcal{F}_k^{'}] + \frac{L_{\Phi}}{2} \mathbb{E}[\|\boldsymbol{x}^{k+1} - \boldsymbol{x}^k\|^2 | \mathcal{F}_{k}^{'}]\\   \nonumber
&\overset{(i)}{=} \Phi(\boldsymbol{x}^k) - \alpha \mathbb{E}[\langle \nabla \Phi(\boldsymbol{x}^k), \boldsymbol{h}_f^k \rangle | \mathcal{F}_k^{'}] + \frac{L_{\Phi}}{2}\alpha^2 \mathbb{E}[\| \boldsymbol{h}_f^k \|^2 | \mathcal{F}_k^{'}]\\   \nonumber
&\overset{(ii)}{=} \Phi(\boldsymbol{x}^k) - \alpha \left \langle \nabla \Phi(\boldsymbol{x}^k), \bar{\boldsymbol{h}}_f^k \right \rangle + \frac{L_{\Phi}}{2} \alpha^2 \mathbb{E}[\| \boldsymbol{h}_f^k \|^2 | \mathcal{F}_k^{'}]\\   \nonumber
&\overset{(iii)}{=} \Phi(\boldsymbol{x}^k) - \frac{\alpha}{2} \left\| \nabla \Phi(\boldsymbol{x}^k)\right\|^2 - \left(\frac{\alpha}{2} - \frac{L_{\Phi}}{2}\alpha^2\right)\left\| \bar{\boldsymbol{h}}_f^k \right\|^2 + \frac{\alpha}{2}\| \nabla \Phi(\boldsymbol{x}^k) - \bar{\boldsymbol{h}}_f^k \|^2 \\   \nonumber
& \qquad +  \frac{L_{\Phi}}{2} \alpha^2 \mathbb{E}[\| \boldsymbol{h}_f^k - \bar{\boldsymbol{h}}_f^k \|^2 | \mathcal{F}_k^{'}], \nonumber
\end{align*}
where $(i)$ follows from the update step of $\boldsymbol{x}^k$ in line 2 of Algorithm 2, $(ii)$ holds since $\mathbb{E}[\langle \nabla \Phi(\boldsymbol{x}^k), \boldsymbol{h}_f^k \rangle | \mathcal{F}_k^{'}] = \langle  \nabla \Phi(\boldsymbol{x}^k), \mathbb{E}[\boldsymbol{h}_f^k | \mathcal{F}_k^{'}]\rangle$, and $(iii)$ holds since $\mathbb{E}[\langle \boldsymbol{h}_f^k - \bar{\boldsymbol{h}}_f^k, \bar{\boldsymbol{h}}_f^k\rangle | \mathcal{F}_k^{'}] = 0$.

Taking a full expectation, we have
\begin{align*}
\mathbb{E}[\Phi(\boldsymbol{x}^{k+1})] & \le\mathbb{E}[\Phi(\boldsymbol{x}^k)] - \frac{\alpha}{2}\mathbb{E}\| \nabla \Phi(\boldsymbol{x}^k) \|^2 - \left(\frac{\alpha}{2} - \frac{L_{\Phi}}{2}\alpha^2 \right)\mathbb{E}\|\bar{\boldsymbol{h}}_f^k\|^2 \\
&\qquad + \frac{\alpha}{2}\mathbb{E}\| \nabla \Phi(\boldsymbol{x}^k)-\bar{\boldsymbol{h}}_f^k\|^2 + \frac{L_{\Phi}}{2}\alpha^2 \mathbb{E}\|\boldsymbol{h}_f^k - \bar{\boldsymbol{h}}_f^k \|^2.
\end{align*}
Then, the proof is complete.

\end{IEEEproof}
Lemma \ref{lemma:3} shows that the variation of function value $\Phi(\boldsymbol{x})$ after one step of iteration \big(i.e., $\mathbb{E}[\Phi(\boldsymbol{x}^{k+1})] - \mathbb{E}[\Phi(\boldsymbol{x}^k)]$\big) depends on the errors $\mathbb{E}\| \nabla \Phi(\boldsymbol{x}^k)- \bar{\boldsymbol{h}}_f^k \|^2$ and $\mathbb{E}\|\boldsymbol{h}_f^k - \bar{\boldsymbol{h}}_f^k\|^2$. In the following, we analyze them.

\begin{lemma}\label{lemma:4}
Suppose Assumptions 1 and 2 hold. For $\boldsymbol{h}_f^k$ in (\ref{eq18}), $\bar{\boldsymbol{h}}_f^k = \mathbb{E}[\boldsymbol{h}_f^k | \mathcal{F}_k^{'}]$, and $\nabla \Phi(x^k)$ in (\ref{eq17}), we have 
\begin{equation}  \tag{20}  \label{eq20}
\mathbb{E}\| \nabla \Phi(\boldsymbol{x}^k) - \bar{\boldsymbol{h}}_f^k \|^2 \le 2 \bar{C}_1^2 \mathbb{E}\| \boldsymbol{y}^k -\boldsymbol{y}^*(\boldsymbol{x}^k) \|^2 +  2 L^2 \mathbb{E}\|\boldsymbol{v}^k - \boldsymbol{v}(\boldsymbol{x}^k, \boldsymbol{y}^*(\boldsymbol{x}^k))\|^2,
\end{equation}
where $\bar{C}_1=L + {M \tau}/{\mu}$.
\end{lemma}

\begin{IEEEproof}
Based on the definition of $\boldsymbol{h}_f^k$ in (\ref{eq18}), we have
\begin{equation} \tag{21} \label{eq21}
\bar{\boldsymbol{h}}_f^k = \nabla_{\boldsymbol{x}} f(\boldsymbol{x}^k, \boldsymbol{y}^k) - \nabla_{\boldsymbol{x}} \nabla_{\boldsymbol{y}} g(\boldsymbol{x}^k, \boldsymbol{y}^k) \mathbb{E}[\boldsymbol{v}^k| \mathcal{F}_k^{'}].
\end{equation}

By combining (\ref{eq21}) with $\nabla \Phi(\boldsymbol{x}^k)$ in (\ref{eq17}), we have

\begin{align*}
&\big\| \nabla \Phi(\boldsymbol{x}^k) - \bar{\boldsymbol{h}}_f^k\big\|  \\   \nonumber
&\overset{(i)}{\le} \big\| \nabla_{\boldsymbol{x}} f\left(\boldsymbol{x}^k, \boldsymbol{y}^*(\boldsymbol{x}^k)\right) - \nabla_{\boldsymbol{x}} f\left(\boldsymbol{x}^k, \boldsymbol{y}^k\right)\big\| \\
& \qquad +\big \| \nabla_{\boldsymbol{x}} \nabla_{\boldsymbol{y}} g\left(\boldsymbol{x}^k, \boldsymbol{y}^*(\boldsymbol{x}^k)\right) \boldsymbol{v}\left(\boldsymbol{x}^k, \boldsymbol{y}^*(\boldsymbol{x}^k)\right) - \nabla_{\boldsymbol{x}} \nabla_{\boldsymbol{y}} g(\boldsymbol{x}^k, \boldsymbol{y}^k) \mathbb{E}[\boldsymbol{v}^k| \mathcal{F}_k^{'}]\big\| \\   \nonumber
& \overset{(ii)}{\le}  L \| \boldsymbol{y}^k - \boldsymbol{y}^*(\boldsymbol{x}^k)\| + \big\| \boldsymbol{v}\left(\boldsymbol{x}^k, \boldsymbol{y}^*(\boldsymbol{x}^k)\right)\big\| \big\| \nabla_{\boldsymbol{x}} \nabla_{\boldsymbol{y}} g\left(\boldsymbol{x}^k, \boldsymbol{y}^*(\boldsymbol{x}^k)\right) - \nabla_{\boldsymbol{x}} \nabla_{\boldsymbol{y}} g(\boldsymbol{x}^k, \boldsymbol{y}^k)\big\|\\  \nonumber
&\qquad  + \| \nabla_{\boldsymbol{x}} \nabla_{\boldsymbol{y}} g(\boldsymbol{x}^k, \boldsymbol{y}^k)\| \big\| \boldsymbol{v}\left(\boldsymbol{x}^k,\boldsymbol{y}^*(\boldsymbol{x}^k)\right) - \mathbb{E}[\boldsymbol{v}^k| \mathcal{F}_k^{'}]\big\|\\  \nonumber
&\overset{(iii)}{\le} (L + {M \tau}/{\mu}) \| \boldsymbol{y}^k - \boldsymbol{y}^*(\boldsymbol{x}^k)\| +  L\big\| \mathbb{E}[\boldsymbol{v}^k - \boldsymbol{v}\left(\boldsymbol{x}^k, \boldsymbol{y}^*(\boldsymbol{x}^k)\right)| \mathcal{F}_k^{'}]\big\|,  \nonumber
\end{align*}
where $(i)$ follows from the triangle inequality, $(ii)$ follows from Assumption 2, triangle inequality and Cauchy-Schwartz inequality, and $(iii)$ follows from Assumption 2 and the fact that $\| \boldsymbol{v}(\boldsymbol{x}^k, \boldsymbol{y}^*(\boldsymbol{x}^k))\| \le {M}/{\mu}$. 

Then, notice that $(a+b)^2 \le 2a^2 + 2b^2$ for $a$, $b \in \mathbb{R}$, we have 

\begin{align*}
\| \nabla \Phi(\boldsymbol{x}^k) - \bar{\boldsymbol{h}}_f^k \|^2 &\le 2 \left(L + \frac{M}{\mu}\tau \right)^2 \|\boldsymbol{y}^k - \boldsymbol{y}^*(\boldsymbol{x}^k)\|^2 + 2L^2 \big\| \mathbb{E}[\boldsymbol{v}^k -\boldsymbol{v}(\boldsymbol{x}^k, \boldsymbol{y}^*(\boldsymbol{x}^k))| \mathcal{F}_k^{'}]\big\|^2 \\  \nonumber
&\overset{(i)}{\le} 2 \left(L +  \frac{M}{\mu}\tau \right)^2 \| \boldsymbol{y}^k - \boldsymbol{y}^*(\boldsymbol{x}^k)\|^2 + 2 L^2 \mathbb{E}[\|\boldsymbol{v}^k - \boldsymbol{v}(\boldsymbol{x}^k, \boldsymbol{y}^*(\boldsymbol{x}^k))\|^2 | \mathcal{F}_k^{'}],  \nonumber
\end{align*}
where $(i)$ holds since $\| \mathbb{E}[\boldsymbol{x}|\mathcal{F}_k^{'}]\|^2 \le \mathbb{E}[\|\boldsymbol{x}\|^2|\mathcal{F}_k^{'}]$ for any vector $\boldsymbol{x}\in \mathbb{R}^q$.

Taking a full expectation, we have 
\begin{equation*}
\mathbb{E}\| \nabla \Phi(\boldsymbol{x}^k) - \bar{\boldsymbol{h}}_f^k\|^2 \le 2 \left(L +  \frac{M}{\mu}\tau \right)^2 \mathbb{E}\|\boldsymbol{y}^k - \boldsymbol{y}^*(\boldsymbol{x}^k)\|^2 + 2 L^2 \mathbb{E}\|\boldsymbol{v}^k - \boldsymbol{v}(\boldsymbol{x}^k,\boldsymbol{y}^*(\boldsymbol{x}^k))\|^2.
\end{equation*}
Then, the proof is complete.
\end{IEEEproof}
\begin{lemma} \label{lemma:5}
Suppose Assumptions 1, 2, and 3 hold. Let $\|\boldsymbol{v}^0\|\le M/\mu$ and $\eta \le 1/L$. For $\boldsymbol{h}_f^k$ in (\ref{eq18}) and $\bar{\boldsymbol{h}}_f^k = \mathbb{E}[\boldsymbol{h}_f^k| \mathcal{F}_k^{'}]$, we have 
\begin{equation}  \label{eq22} \tag{22}
\mathbb{E}\|\boldsymbol{h}_f^k - \bar{\boldsymbol{h}}_f^k \|^2 \le \frac{\sigma_f^2}{D_f} + \frac{2M^2}{\mu^2}\frac{\sigma_{g, 1}^2}{D_g} + 2L^2 \mathbb{E}\|\boldsymbol{v}^k - \boldsymbol{v}(\boldsymbol{x}^k, \boldsymbol{y}^*(\boldsymbol{x}^k))\|^2,
\end{equation}
where $\boldsymbol{v}(\boldsymbol{x}^k,\boldsymbol{y}^*(\boldsymbol{x}^k)) = [\nabla_{\boldsymbol{y}}^2 g(\boldsymbol{x}^k, \boldsymbol{y}^*(\boldsymbol{x}^k))]^{-1} \nabla_{\boldsymbol{y}} f(\boldsymbol{x}^k, \boldsymbol{y}^*(\boldsymbol{x}^k))$.
\end{lemma}
\begin{IEEEproof}
Based on the definition of $\boldsymbol{h}_f^k$ in (\ref{eq18}), we have
\begin{equation*}
\bar{\boldsymbol{h}}_f^k = \nabla_{\boldsymbol{x}} f(\boldsymbol{x}^k, \boldsymbol{y}^k) - \nabla_{\boldsymbol{x}} \nabla_{\boldsymbol{y}} g(\boldsymbol{x}^k, \boldsymbol{y}^k) \mathbb{E}[\boldsymbol{v}^k| \mathcal{F}_k^{'}].
\end{equation*}
Then, taking the conditional expectation on $\mathcal{F}_{k}^{'}$, we have 
\begin{align}  \label{eq23}  \tag{23}
\mathbb{E}[\| \boldsymbol{h}_f^k - \bar{\boldsymbol{h}}_f^k \|^2 | \mathcal{F}_k^{'}] &\overset{(i)}{=} \mathbb{E}[\| \nabla_{\boldsymbol{x}} F(\boldsymbol{x}^k, \boldsymbol{y}^k; \mathcal{D}_F) - \nabla_{\boldsymbol{x}} f(\boldsymbol{x}^k, \boldsymbol{y}^k) \|^2|\mathcal{F}_k^{'}]\\
& \qquad +\mathbb{E}\big[\big\| \nabla_{\boldsymbol{x}} \nabla_{\boldsymbol{y}} G(\boldsymbol{x}^k, \boldsymbol{y}^k; \mathcal{D}_G)\boldsymbol{v}^k  - \nabla_{\boldsymbol{x}} \nabla_{\boldsymbol{y}} g(\boldsymbol{x}^k, \boldsymbol{y}^k) \mathbb{E}[\boldsymbol{v}^k| \mathcal{F}_k^{'}]\big\|^2 \big{|} \mathcal{F}_k^{'}\big],  \nonumber
\end{align}
where $(i)$ holds since $\{ \mathcal{D}_F, \mathcal{D}_G, \mathcal{B}_j(j=0, \ldots, J-1), \mathcal{D}_{F,j}(j=0, \ldots, J-1)\}$ are mutually independent sample sets.

As for the second term on the right of the equality in (\ref{eq23}), we have
\begin{align*}  \label{eq24}  
&\mathbb{E}\big[\big\| \nabla_{\boldsymbol{x}} \nabla_{\boldsymbol{y}} G(\boldsymbol{x}^k, \boldsymbol{y}^k; \mathcal{D}_G)\boldsymbol{v}^k - \nabla_{\boldsymbol{x}} \nabla_{\boldsymbol{y}} g(\boldsymbol{x}^k,\boldsymbol{y}^k) \mathbb{E}[\boldsymbol{v}^k| \mathcal{F}_k^{'}] \big\|^2 \big{|} \mathcal{F}_k^{'}\big]\\   \nonumber
&\overset{(i)}{\le} 2\mathbb{E}[\| \nabla_{\boldsymbol{x}} \nabla_{\boldsymbol{y}} G(\boldsymbol{x}^k, \boldsymbol{y}^k; \mathcal{D}_G)\boldsymbol{v}^k - \nabla_{\boldsymbol{x}} \nabla_{\boldsymbol{y}} g(\boldsymbol{x}^k, \boldsymbol{y}^k) \boldsymbol{v}^k \|^2 | \mathcal{F}_k^{'}]\\   \nonumber   
&\qquad  + 2 \mathbb{E}\big[\big\| \nabla_{\boldsymbol{x}} \nabla_{\boldsymbol{y}} g(\boldsymbol{x}^k, \boldsymbol{y}^k)\boldsymbol{v}^k - \nabla_{\boldsymbol{x}} \nabla_{\boldsymbol{y}} g(\boldsymbol{x}^k, \boldsymbol{y}^k)E[\boldsymbol{v}^k|\mathcal{F}_k^{'}]\big\|^2 \big{|} \mathcal{F}_k^{'}\big]\\   \nonumber
& \overset{(ii)}{\le} 2\mathbb{E}\big[\big\| \nabla_{\boldsymbol{x}} \nabla_{\boldsymbol{y}} G(\boldsymbol{x}^k, \boldsymbol{y}^k; \mathcal{D}_G) - \nabla_{\boldsymbol{x}} \nabla_{\boldsymbol{y}} g(\boldsymbol{x}^k, \boldsymbol{y}^k)\big\|^2 \big{|} \mathcal{F}_k^{'}] \mathbb{E}\big[\|\boldsymbol{v}^k\|^2 | \mathcal{F}_k^{'}\big]\\
& \qquad + 2 \| \nabla_{\boldsymbol{x}} \nabla_{\boldsymbol{y}} g(\boldsymbol{x}^k, \boldsymbol{y}^k)\|^2 \mathbb{E}[\| \boldsymbol{v}^k - \mathbb{E}[\boldsymbol{v}^k|\mathcal{F}_k^{'}] \|^2 | \mathcal{F}_k^{'}]\\  \tag{24}
&\overset{(iii)}{\le} \frac{2M^2}{\mu^2} \mathbb{E}[\| \nabla_{\boldsymbol{x}} \nabla_{\boldsymbol{y}} G(\boldsymbol{x}^k, \boldsymbol{y}^k; \mathcal{D}_G) - \nabla_{\boldsymbol{x}} \nabla_{\boldsymbol{y}} g(\boldsymbol{x}^k, \boldsymbol{y}^k)\|^2 | \mathcal{F}_k^{'}] + 2L^2 \mathbb{E}[\|\boldsymbol{v}^k - \mathbb{E}[\boldsymbol{v}^k| \mathcal{F}_k^{'}]\|^2 | \mathcal{F}_k^{'}],  \\   \nonumber
&\overset{(iv)}{\le} \frac{2M^2}{\mu^2} \mathbb{E}[\| \nabla_{\boldsymbol{x}} \nabla_{\boldsymbol{y}} G(\boldsymbol{x}^k,\boldsymbol{y}^k; \mathcal{D}_G) - \nabla_{\boldsymbol{x}} \nabla_{\boldsymbol{y}} g(\boldsymbol{x}^k, \boldsymbol{y}^k) \|^2 | \mathcal{F}_k^{'}] + 2L^2 \mathbb{E}[\|\boldsymbol{v}^k - \boldsymbol{v}(\boldsymbol{x}^k, \boldsymbol{y}^*(\boldsymbol{x}^k))\|^2 | \mathcal{F}_k^{'}],  \\    \nonumber
\end{align*}
where $(i)$ holds since $\| \boldsymbol{a}+\boldsymbol{b}\|^2 \le 2 \|\boldsymbol{a}\|^2 + 2\|\boldsymbol{b} \|^2$ for $\boldsymbol{a}$, $\boldsymbol{b}\in \mathbb{R}^p$, $(ii)$ follows from the Cauchy-Schwartz inequality and the fact that $\{ \mathcal{D}_G, \mathcal{B}_j(j=0, \ldots, J-1), \mathcal{D}_{F,j}(j=0, \ldots, J-1)\}$ are mutually independent sample sets, $(iii)$ follows from Assumption 2, Lemma \ref{lemma:2}, and the fact that $\boldsymbol{v}^k = \boldsymbol{v}^{k, J}$(see line 13 in Algorithm 2), and $(iv)$ holds since 
\begin{equation*}
\mathbb{E}[\|\boldsymbol{v}^k - \boldsymbol{v}(\boldsymbol{x}^k, \boldsymbol{y}^*(\boldsymbol{x}^k))\|^2 | \mathcal{F}_k^{'}] = \mathbb{E}[\|\boldsymbol{v}^k - \mathbb{E}[\boldsymbol{v}^k|\mathcal{F}_k^{'}]\|^2 |\mathcal{F}_k^{'}] + \mathbb{E}[\| \mathbb{E}[\boldsymbol{v}^k|\mathcal{F}_k^{'}] - \boldsymbol{v}(\boldsymbol{x}^k, \boldsymbol{y}^*(\boldsymbol{x}^k))\|^2 | \mathcal{F}_k^{'}].
\end{equation*}

By combining (\ref{eq23}) and (\ref{eq24}), we have
\begin{align*}  
\mathbb{E}[\| \boldsymbol{h}_f^k - \bar{\boldsymbol{h}}_f^k\|^2 |\mathcal{F}_k^{'}] &\le \mathbb{E}[\| \nabla_{\boldsymbol{x}} F(\boldsymbol{x}^k, \boldsymbol{y}^k; \mathcal{D}_F) - \nabla_{\boldsymbol{x}} f(\boldsymbol{x}^k, \boldsymbol{y}^k)\|^2 | \mathcal{F}_k^{'}] \\
&\qquad  + \frac{2M^2}{\mu^2}\mathbb{E}[\| \nabla_{\boldsymbol{x}} \nabla_{\boldsymbol{y}} G(\boldsymbol{x}^k, \boldsymbol{y}^k; \mathcal{D}_G) - \nabla_{\boldsymbol{x}} \nabla_{\boldsymbol{y}} g(\boldsymbol{x}^k, \boldsymbol{y}^k)\|^2 | \mathcal{F}_k^{'}] \\
&\qquad + 2L^2 \mathbb{E}[\|\boldsymbol{v}^k - \boldsymbol{v}(\boldsymbol{x}^k, \boldsymbol{y}^*(\boldsymbol{x}^k))\|^2 | \mathcal{F}_k^{'}].
\end{align*}

Taking a full expectation, we have 
\begin{align*}  
\mathbb{E}\| \boldsymbol{h}_f^k - \bar{\boldsymbol{h}}_f^k \|^2  &\le \mathbb{E}\| \nabla_{\boldsymbol{x}} F(\boldsymbol{x}^k,\boldsymbol{y}^k; \mathcal{D}_F) - \nabla_{\boldsymbol{x}} f(\boldsymbol{x}^k, \boldsymbol{y}^k)\|^2  \\
&\qquad  + \frac{2M^2}{\mu^2}\mathbb{E}\| \nabla_{\boldsymbol{x}} \nabla_{\boldsymbol{y}} G(\boldsymbol{x}^k,\boldsymbol{y}^k; \mathcal{D}_G) - \nabla_{\boldsymbol{x}} \nabla_{\boldsymbol{y}} g(\boldsymbol{x}^k, \boldsymbol{y}^k)\|^2  \\
& \qquad + 2L^2 \mathbb{E}\|\boldsymbol{v}^k - \boldsymbol{v}(\boldsymbol{x}^k, \boldsymbol{y}^*(\boldsymbol{x}^k))\|^2 \\
&\overset{(i)}{\le} \frac{\sigma_f^2}{D_f} + \frac{2M^2}{\mu^2}\frac{\sigma_{g, 1}^2}{D_g} + 2L^2 \mathbb{E}\|\boldsymbol{v}^k - \boldsymbol{v}(\boldsymbol{x}^k,\boldsymbol{y}^*(\boldsymbol{x}^k))\|^2,
\end{align*}
where $(i)$ follows from Assumption 3 and the fact that the sizes of $\mathcal{D}_F$ and $\mathcal{D}_G$ are $D_f$ and $D_g$, respectively. Then, the proof is complete.
\end{IEEEproof}

It is easy to observe that the upper bounds in Lemmas \ref{lemma:4} and \ref{lemma:5} depend on the errors of LL variable $\boldsymbol{y}^k$ and the approximate solution $\boldsymbol{v}^k$ to the linear system. In the following, we first analyze the error of the LL variable $\boldsymbol{y}^k$.

\begin{lemma} \label{lemma:6}
Suppose Assumptions 1, 2, and 3 hold. Let $\beta \le {2}/{(\mu + L)}$ and $r_v$ be any positive number, we have 
\begin{align}   \label{eq25}  \tag{25}
\mathbb{E}\| \boldsymbol{y}^{k+1} - \boldsymbol{y}^*(\boldsymbol{x}^{k+1})\|^2 & \le (1-\rho_y)^T(1+r_v)\mathbb{E}\| \boldsymbol{y}^k - \boldsymbol{y}^*(\boldsymbol{x}^k) \|^2 \\
&\qquad +(1 - \rho_y)^T\left(1 + \frac{1}{r_v}\right) \frac{L^2}{\mu^2}\alpha^2 \mathbb{E} \| \boldsymbol{h}_f^k \|^2  + \frac{\mu +L }{2\mu L}\beta \frac{1}{S}\sigma_g^2,  \nonumber
\end{align}
where $\rho_y = 2 \beta {\mu L}/{(\mu + L)}$.
\end{lemma}

\begin{IEEEproof}
For  any $0 \le t \le T-1$, taking the conditional expectation on $\mathcal{F}_{k+1, t}$, we have
\begin{align}  \label{eq26}     \nonumber
&\mathbb{E}[\| \boldsymbol{y}^{k+1, t+1} - \boldsymbol{y}^*(\boldsymbol{x}^{k+1})\|^2 | \mathcal{F}_{k+1, t}] \\    \nonumber
&\overset{(i)}{=} \mathbb{E}[\| \boldsymbol{y}^{k+1, t} - \beta \nabla_{\boldsymbol{y}} G(\boldsymbol{x}^{k+1}, \boldsymbol{y}^{k+1, t}; \mathcal{S}_t) - \boldsymbol{y}^*(\boldsymbol{x}^{k+1})\|^2 | \mathcal{F}_{k+1, t}]\\   \nonumber
&= \| \boldsymbol{y}^{k+1, t} - \boldsymbol{y}^*(\boldsymbol{x}^{k+1})\|^2 - 2 \beta \mathbb{E}\left[\left \langle  \nabla_{\boldsymbol{y}} G(\boldsymbol{x}^{k+1}, \boldsymbol{y}^{k+1, t}; \mathcal{S}_t), \boldsymbol{y}^{k+1, t} - \boldsymbol{y}^*(\boldsymbol{x}^{k+1}) \right\rangle \big| \mathcal{F}_{k+1, t}\right]\\  \tag{26}
& \qquad + \beta^2 \mathbb{E}[\| \nabla_{\boldsymbol{y}} G(\boldsymbol{x}^{k+1}, \boldsymbol{y}^{k+1, t}; \mathcal{S}_t)\|^2 | \mathcal{F}_{k+1, t}]\\   \nonumber
& = \| \boldsymbol{y}^{k+1, t} - \boldsymbol{y}^*(\boldsymbol{x}^{k+1})\|^2 - \underbrace{2 \beta \left\langle  \nabla_{\boldsymbol{y}} g(\boldsymbol{x}^{k+1}, \boldsymbol{y}^{k+1, t}), \boldsymbol{y}^{k+1, t} - \boldsymbol{y}^*(\boldsymbol{x}^{k+1}) \right\rangle}_{J_1}\\   \nonumber
&\qquad +\underbrace{\beta^2 \mathbb{E}[\| \nabla_{\boldsymbol{y}} G(\boldsymbol{x}^{k+1}, \boldsymbol{y}^{k+1, t}; \mathcal{S}_t)\|^2 | \mathcal{F}_{k+1, t}]}_{J_2},   \nonumber
\end{align}
where $(i)$ follows from the iteration step in line 6 of Algorithm 2.

For $J_1$, since $g$ is $L$-smooth $\mu$-strongly convex, from Theoem $2.1.12$ in \cite{25}, we have
\begin{equation} \label{eq27}  \tag{27}
J_1 \ge 2\beta \left( \frac{\mu L}{\mu + L}\| \boldsymbol{y}^{k+1, t} - \boldsymbol{y}^*(\boldsymbol{x}^{k+1})\|^2 + \frac{1}{\mu + L}\| \nabla_{\boldsymbol{y}} g(\boldsymbol{x}^{k+1}, \boldsymbol{y}^{k+1, t})\|^2 \right).
\end{equation}

For $J_2$, 
\begin{align}  \label{eq28} \tag{28}
& \beta^2 \mathbb{E}[ \| \nabla_{\boldsymbol{y}} G(\boldsymbol{x}^{k+1}, \boldsymbol{y}^{k+1, t}; \mathcal{S}_t)\|^2 | \mathcal{F}_{k+1, t} ]  \\   \nonumber
& \overset{(i)}{=} \beta^2 \mathbb{E}[ \| \nabla_{\boldsymbol{y}} G(\boldsymbol{x}^{k+1}, \boldsymbol{y}^{k+1, t}; \mathcal{S}_t) - \nabla_{\boldsymbol{y}} g(\boldsymbol{x}^{k+1}, \boldsymbol{y}^{k+1, t})\|^2 | \mathcal{F}_{k+1, t}] + \beta^2 \| \nabla_{\boldsymbol{y}} g(\boldsymbol{x}^{k+1}, \boldsymbol{y}^{k+1, t})\|^2, \nonumber
\end{align}
where $(i)$ holds since
\begin{equation*}
\mathbb{E}[\langle \nabla_{\boldsymbol{y}} G(\boldsymbol{x}^{k+1}, \boldsymbol{y}^{k+1, t}; \mathcal{S}_t) - \nabla_{\boldsymbol{y}} g(\boldsymbol{x}^{k+1}, \boldsymbol{y}^{k+1, t}), \nabla_{\boldsymbol{y}} g(\boldsymbol{x}^{k+1}, \boldsymbol{y}^{k+1, t})\rangle | \mathcal{F}_{k+1, t}] = 0. 
\end{equation*}

Plugging (\ref{eq27}) and (\ref{eq28}) into (\ref{eq26}), we have
\begin{align*}
& \mathbb{E}[ \| \boldsymbol{y}^{k+1, t+1} - \boldsymbol{y}^*(\boldsymbol{x}^{k+1})\|^2| \mathcal{F}_{k+1, t} ] \\
& \le \| \boldsymbol{y}^{k+1, t} - \boldsymbol{y}^*(\boldsymbol{x}^{k+1})\|^2 - 2 \beta \left( \frac{\mu L}{\mu + L}\| \boldsymbol{y}^{k+1, t} -\boldsymbol{y}^*(\boldsymbol{x}^{k+1})\|^2 + \frac{1}{\mu + L}\| \nabla_{\boldsymbol{y}} g(\boldsymbol{x}^{k+1}, \boldsymbol{y}^{k+1, t})\|^2 \right) \\
& \qquad + \beta^2 \mathbb{E}[\| \nabla_{\boldsymbol{y}} G(\boldsymbol{x}^{k+1}, \boldsymbol{y}^{k+1, t}; \mathcal{S}_t) - \nabla_{\boldsymbol{y}} g(\boldsymbol{x}^{k+1}, \boldsymbol{y}^{k+1, t})\|^2 | \mathcal{F}_{k+1, t}] + \beta^2 \| \nabla_{\boldsymbol{y}} g(\boldsymbol{x}^{k+1}, \boldsymbol{y}^{k+1, t})\|^2.
\end{align*}

By rearranging the terms, we have 
\begin{align*}
& \mathbb{E}[ \| \boldsymbol{y}^{k+1, t+1} - \boldsymbol{y}^*(\boldsymbol{x}^{k+1})\|^2 | \mathcal{F}_{k+1, t} ]\\
& \le \left(1 - 2 \beta \frac{\mu L}{\mu + L}\right)\| \boldsymbol{y}^{k+1, t} - \boldsymbol{y}^*(\boldsymbol{x}^{k+1})\|^2 + \beta\left(\beta - \frac{2}{\mu + L}\right)\| \nabla_{\boldsymbol{y}} g(\boldsymbol{x}^{k+1}, \boldsymbol{y}^{k+1, t})\|^2 \\
& \qquad +  \beta^2 \mathbb{E}[\| \nabla_{\boldsymbol{y}} G(\boldsymbol{x}^{k+1}, \boldsymbol{y}^{k+1, t}; \mathcal{S}_t) - \nabla_{\boldsymbol{y}} g(\boldsymbol{x}^{k+1}, \boldsymbol{y}^{k+1, t})\|^2 | \mathcal{F}_{k+1, t}] \\
&\overset{(i)}{\le} (1 - \rho_y) \| \boldsymbol{y}^{k+1, t} - \boldsymbol{y}^*(\boldsymbol{x}^{k+1})\|^2 + \beta^2 \mathbb{E}[\| \nabla_{\boldsymbol{y}} G(\boldsymbol{x}^{k+1}, \boldsymbol{y}^{k+1, t}; \mathcal{S}_t) - \nabla_{\boldsymbol{y}} g(\boldsymbol{x}^{k+1}, \boldsymbol{y}^{k+1, t})\|^2 | \mathcal{F}_{k+1, t}],
\end{align*}
where $(i)$ follows from the definition of $\rho_y= 2 \beta {\mu L}/{(\mu + L)}$ and the fact that $\beta \le 2/(\mu + L)$. 

Taking a full expectation, we have
\begin{equation*}
\mathbb{E}\| \boldsymbol{y}^{k+1, t+1}- \boldsymbol{y}^*(\boldsymbol{x}^{k+1})\|^2   \overset{(i)}{\le} (1 - \rho_y) \mathbb{E}\| \boldsymbol{y}^{k+1, t} - \boldsymbol{y}^*(\boldsymbol{x}^{k+1})\|^2 + \beta^2 \frac{1}{S}\sigma_g^2,
\end{equation*}
where $(i)$ follows from Assumption 3 and the fact that the size of $\mathcal{S}_t$ is $S$.

Since $\boldsymbol{y}^{k+1} = \boldsymbol{y}^{k+1, T}$(in line 8 of Algorithm 2) and $\boldsymbol{y}^k =\boldsymbol{y}^{k+1, 0}$(in line 3 of Algorithm 2), by using induction, we have 
\begin{align}   \label{eq29}  \tag{29}
\mathbb{E}\| \boldsymbol{y}^{k+1}- \boldsymbol{y}^*(x^{k+1})\|^2& \le (1 - \rho_y)^T \mathbb{E}\| \boldsymbol{y}^k - \boldsymbol{y}^*(\boldsymbol{x}^{k+1})\|^2+ \sum_{t=0}^{T-1}(1- \rho_y)^t \beta^2 \frac{1}{S}\sigma_g^2\\  \nonumber
& \overset{(i)}{\le}(1 - \rho_y)^T \mathbb{E}\| \boldsymbol{y}^k - \boldsymbol{y}^*(\boldsymbol{x}^{k+1})\|^2 + \frac{\mu + L}{2 \mu L}\beta \frac{1}{S}\sigma_g^2,   \nonumber
\end{align}
where $(i)$ follows from the fact that $\sum_{t=0}^{T-1}(1- \rho_y)^t \le \frac{\mu + L }{2\beta \mu L}$. 

Since
\begin{align}   \label{eq30} \tag{30}
&\mathbb{E}\|\boldsymbol{y}^k - \boldsymbol{y}^*(\boldsymbol{x}^{k+1})\|^2 \\   \nonumber
& \overset{(i)}{\le} (1 + r_v)\mathbb{E}\|\boldsymbol{y}^k - \boldsymbol{y}^*(\boldsymbol{x}^k)\|^2 + \left( 1 + \frac{1}{r_v}\right)\mathbb{E}\|\boldsymbol {y}^*(\boldsymbol{x}^k) - \boldsymbol{y}^*(\boldsymbol{x}^{k+1})\|^2 \\   \nonumber
& \overset{(ii)}{\le} (1 + r_v)\mathbb{E}\| \boldsymbol{y}^k - \boldsymbol{y}^*(\boldsymbol{x}^k)\|^2+ \left( 1 + \frac{1}{r_v}\right) \frac{L^2}{\mu^2}\mathbb{E}\| \boldsymbol{x}^k - \boldsymbol{x}^{k+1}\|^2  \\  \nonumber
&\overset{(iii)}{=} (1 + r_v)\mathbb{E}\| \boldsymbol{y}^k - \boldsymbol{y}^*(\boldsymbol{x}^k)\|^2 + \left( 1 + \frac{1}{r_v}\right) \frac{L^2}{\mu^2}\alpha^2 \mathbb{E}\|\boldsymbol{h}_f^k\|^2, \nonumber
\end{align}
where $(i)$ follows from the fact that $r_v \|\boldsymbol{a} \|^2 + {1}/{r_v} \|\boldsymbol{b} \|^2 \ge 2 \boldsymbol{a}^\top \boldsymbol{b}$ for $\boldsymbol{a}$, $\boldsymbol{b} \in \mathbb{R}^q$, $(ii)$ follows from Proposition 1, and $(iii)$ follows from the iteration step in line 2 of Algorithm 2.

By combining (\ref{eq29}) with (\ref{eq30}), we have 
\begin{align*}
\mathbb{E}\| \boldsymbol{y}^{k+1} - \boldsymbol{y}^*(\boldsymbol{x}^{k+1})\|^2  &\le (1 - \rho_y)^{T}(1 + r_v)\mathbb{E}\|\boldsymbol{y}^k - \boldsymbol{y}^*(\boldsymbol{x}^k)\|^2 \\
&\qquad + (1 - \rho_y)^T \left(1 + \frac{1}{r_v}\right)\frac{L^2}{\mu^2}\alpha^2 \mathbb{E}\| \boldsymbol{h}_f^k \|^2 
 + \frac{\mu + L}{2 \mu L}\beta \frac{1}{S}\sigma_g^2.
\end{align*}
Then, the proof is complete.
\end{IEEEproof}
Next, we analyze the error of $\boldsymbol{v}^k$.

\begin{lemma}  \label{lemma:7}
Suppose Assumptions 1, 2, 3 hold. Let $\eta < {1}/{L}$, $r_w$ be any positive number, $\|\boldsymbol{v}^0\|\le M/\mu$, and 
\begin{equation*}
\boldsymbol{v}(\boldsymbol{x}, \boldsymbol{y}) = [\nabla_{\boldsymbol{y}}^2 g(\boldsymbol{x}, \boldsymbol{y})]^{-1} \nabla_{\boldsymbol{y}} f(\boldsymbol{x}, \boldsymbol{y}).
\end{equation*}
Then, we have 
\begin{align}  \label{eq31}  \tag{31}
\mathbb{E}\| \boldsymbol{v}^{k+1} - \boldsymbol{v}(\boldsymbol{x}^{k+1}, \boldsymbol{y}^*(\boldsymbol{x}^{k+1}))\|^2 & \le (1+ r_w)C_{21} \mathbb{E}\|\boldsymbol{v}^k - \boldsymbol{v}(\boldsymbol{x}^k, \boldsymbol{y}^*(\boldsymbol{x}^k))\|^2  \\   \nonumber
&\qquad + (1 + r_w)C_{22}\alpha^2 \mathbb{E}\|\boldsymbol{h}_f^k \|^2  + (1 + r_w)C_{23}\mathbb{E}\|\boldsymbol{y}^k - \boldsymbol{y}^*(\boldsymbol{x}^k)\|^2 \\   \nonumber
&\qquad + C_{24}\mathbb{E}\left\|\boldsymbol{y}^{k+1} - \boldsymbol{y}^*(\boldsymbol{x}^{k+1})\right\|^2   
+ (1 + r_w)\eta \frac{M^2}{\mu^3}\frac{\sigma_{g, 2}^2}{D} + (1 + r_w)\frac{\eta}{\mu}\frac{\sigma_f^2}{D_f},  \nonumber
\end{align}
where 
\begin{align*}
C_{21} = (1 - \eta \mu)^J (1 + r_w)^2,  \qquad & C_{22} = \left(1 + \frac{1}{r_w}\right)\left(2 + 8 \frac{1}{\mu^2}L^2\right)C_1^2,\\
C_{23} = \left(1 + \frac{1}{r_w}\right)C_1^2(5 + r_w), \qquad & C_{24} = \left(1 + \frac{1}{r_w}\right)(9 + 8r_w )C_1^2,
\end{align*}
and $C_1 = \frac{M}{\mu^2}\rho + \frac{L}{\mu}$ is definied in Lemma \ref{lemma:1}.
\end{lemma}

\begin{IEEEproof}
Firstly, we have
\begin{align} \label{eq32}   \nonumber
& \mathbb{E}\| \boldsymbol{v}^{k+1} - \boldsymbol{v}(\boldsymbol{x}^{k+1}, \boldsymbol{y}^*(\boldsymbol{x}^{k+1}))\|^2 \\  \nonumber
& \overset{(i)}{\le} (1 + r_w)\mathbb{E}\| \boldsymbol{v}^{k+1} - \boldsymbol{v}(\boldsymbol{x}^{k+1}, \boldsymbol{y}^{k+1})\|^2  + \left(1 + \frac{1}{r_w}\right)\mathbb{E}\| \boldsymbol{v}(\boldsymbol{x}^{k+1}, \boldsymbol{y}^{k+1}) - \boldsymbol{v}(\boldsymbol{x}^{k+1}, \boldsymbol{y}^*(\boldsymbol{x}^{k+1}))\|^2    \\ \tag{32}
&\overset{(ii)}{\le} (1 + r_w)\mathbb{E}\| \boldsymbol{v}^{k+1} - \boldsymbol{v}(\boldsymbol{x}^{k+1}, \boldsymbol{y}^{k+1})\|^2 + \left(1 + \frac{1}{r_w}\right)C_1^2 \mathbb{E}\|\boldsymbol{y}^{k+1} - \boldsymbol{y}^*(\boldsymbol{x}^{k+1})\|^2, \nonumber
\end{align}
where $(i)$ follows from the fact that $r_w \|\boldsymbol{a}\|^2 + {1}/{r_w}\|\boldsymbol{b} \|^2 \ge 2 \boldsymbol{a}^\top \boldsymbol{b}$ for $\boldsymbol{a}$, $\boldsymbol{b}\in \mathbb{R}^q$, and $(ii)$ follows from Lemma $\ref{lemma:1}$. 

Then, we bound the first term on the right of the inequality in (\ref{eq32}). For any $0 \le j \le J-1$, taking the conditional expectation on $\mathcal{F}_{k+1, j}$, we have 
\begin{align*}
&\mathbb{E}[\|\boldsymbol{v}^{k+1, j+1} - \boldsymbol{v}(\boldsymbol{x}^{k+1}, \boldsymbol{y}^{k+1})\|^2 | \mathcal{F}_{k+1, j}]  \\ 
&\overset{(i)}{=}\mathbb{E}\big[\left\| (\mathbf{I} - \eta \nabla_{\boldsymbol{y}}^2 G(\boldsymbol{x}^{k+1}, \boldsymbol{y}^{k+1}; \mathcal{B}_j))\boldsymbol{v}^{k+1, j} 
 + \eta \nabla_{\boldsymbol{y}} F(\boldsymbol{x}^{k+1}, \boldsymbol{y}^{k+1};\mathcal{D}_{F, j}) \right. \\   
& \qquad \left. - (\mathbf{I} - \eta \nabla_{\boldsymbol{y}}^2 g(\boldsymbol{x}^{k +1}, \boldsymbol{y}^{k+1}))\boldsymbol{v}(\boldsymbol{x}^{k+1}, \boldsymbol{y}^{k+1}) - \eta \nabla_{\boldsymbol{y}} f(\boldsymbol{x}^{k+1}, \boldsymbol{y}^{k+1}) \right\|^2 \big| \mathcal{F}_{k+1, j}\big]  \\
&\overset{(ii)}{=} \left\| (\mathbf{I} - \eta \nabla_{\boldsymbol{y}}^2 g(\boldsymbol{x}^{k+1}, \boldsymbol{y}^{k+1}))(\boldsymbol{v}^{k+1, j} - \boldsymbol{v}(\boldsymbol{x}^{k+1}, \boldsymbol{y}^{k+1})) \right\|^2  \\
& \qquad + \mathbb{E}\big[\left\|\eta( \nabla_{\boldsymbol{y}}^2 g(\boldsymbol{x}^{k+1}, \boldsymbol{y}^{k+1}) - \nabla_{\boldsymbol{y}}^2 G(\boldsymbol{x}^{k+1}, \boldsymbol{y}^{k+1}; \mathcal{B}_j))\boldsymbol{v}^{k+1, j} + \eta(\nabla_{\boldsymbol{y}} F(\boldsymbol{x}^{k+1}, \boldsymbol{y}^{k+1}; \mathcal{D}_{F, j}) \right.   \\
&\qquad \left. - \nabla_{\boldsymbol{y}} f(\boldsymbol{x}^{k+1},\boldsymbol{y}^{k+1}))\right\|^2 \big| \mathcal{F}_{k+1, j}\big]   \\ 
&\overset{(iii)}{=} \big\| (\mathbf{I} - \eta \nabla_{\boldsymbol{y}}^2 g(\boldsymbol{x}^{k+1}, \boldsymbol{y}^{k+1}))(\boldsymbol{v}^{k+1, j} - \boldsymbol{v}(\boldsymbol{x}^{k+1}, \boldsymbol{y}^{k+1})) \big\|^2    \\
&\qquad  +\mathbb{E}[\| \eta (\nabla_{\boldsymbol{y}}^2 g(\boldsymbol{x}^{k+1}, \boldsymbol{y}^{k+1}) - \nabla_{\boldsymbol{y}}^2 G(\boldsymbol{x}^{k+1}, \boldsymbol{y}^{k+1}; \mathcal{B}_j))\boldsymbol{v}^{k+1, j}\|^2 | \mathcal{F}_{k+1, j}]   \\ 
&\qquad  +  \mathbb{E}[\| \eta (\nabla_{\boldsymbol{y}} F(\boldsymbol{x}^{k+1}, \boldsymbol{y}^{k+1}; \mathcal{D}_{F, j}) - \nabla_{\boldsymbol{y}} f(\boldsymbol{x}^{k+1}, \boldsymbol{y}^{k+1}))\|^2 | \mathcal{F}_{k+1, j}]    \\ 
&\overset{(iv)}{\le} \| (\mathbf{I} - \eta \nabla_{\boldsymbol{y}}^2 g(\boldsymbol{x}^{k+1}, \boldsymbol{y}^{k+1})) \|^2 \| \boldsymbol{v}^{k+1, j} - \boldsymbol{v}(\boldsymbol{x}^{k+1}, \boldsymbol{y}^{k+1}) \|^2  \\
&\qquad  + \eta^2 \| \boldsymbol{v}^{k+1, j}\|^2 \mathbb{E}[\| \nabla_{\boldsymbol{y}}^2 g(\boldsymbol{x}^{k+1}, \boldsymbol{y}^{k+1}) - \nabla_{\boldsymbol{y}}^2 G(\boldsymbol{x}^{k+1}, \boldsymbol{y}^{k+1}; \mathcal{B}_j)\|^2 | \mathcal{F}_{k+1, j}] \\ 
&\qquad  + \eta^2 \mathbb{E}[\| \nabla_{\boldsymbol{y}} F(\boldsymbol{x}^{k+1}, \boldsymbol{y}^{k+1}; \mathcal{D}_{F, j}) - \nabla_{\boldsymbol{y}} f(\boldsymbol{x}^{k+1}, \boldsymbol{y}^{k+1})\|^2 | \mathcal{F}_{k+1, j}],
\end{align*}
where $(i)$ uses the iteration step in line 11 of Algorithm 2 and the fact that 
\begin{equation*}
\boldsymbol{v}(\boldsymbol{x}^{k+1}, \boldsymbol{y}^{k+1}) = [\nabla_{\boldsymbol{y}}^2 g(\boldsymbol{x}^{k+1}, \boldsymbol{y}^{k+1})]^{-1} \nabla_{\boldsymbol{y}} f(\boldsymbol{x}^{k+1}, \boldsymbol{y}^{k+1}),
\end{equation*}
$(ii)$ holds since $\mathbb{E}[\| \boldsymbol{x} \|^2 | \mathcal{F}_{k+1, j}] = \| \mathbb{E}[\boldsymbol{x}| \mathcal{F}_{k+1, j}]\|^2 + \mathbb{E}[\| \boldsymbol{x} - \mathbb{E}[\boldsymbol{x}|\mathcal{F}_{k+1, j}]\|^2 | \mathcal{F}_{k+1, j}]$ for $\boldsymbol{x} \in \mathbb{R}^q$, $(iii)$ holds since $\mathcal{B}_j$ and $\mathcal{D}_{F,j}$ are independent, and $(iv)$ follows from the Cauchy-Schwartz inequality.

Taking a full expectation, we have
\begin{align*}
&\mathbb{E}\| \boldsymbol{v}^{k+1, j+1} - \boldsymbol{v}(\boldsymbol{x}^{k+1}, \boldsymbol{y}^{k+1})\|^2  \\
&\overset{(i)}{\le}  (1- \eta \mu) \mathbb{E}\| \boldsymbol{v}^{k+1, j} - \boldsymbol{v}(\boldsymbol{x}^{k+1}, \boldsymbol{y}^{k+1}) \|^2\\
& \qquad + \eta^2 \frac{M^2}{\mu^2} \mathbb{E}\| \nabla_{\boldsymbol{y}}^2 g(\boldsymbol{x}^{k+1}, \boldsymbol{y}^{k+1}) - \nabla_{\boldsymbol{y}}^2 G(\boldsymbol{x}^{k+1}, \boldsymbol{y}^{k+1}; \mathcal{B}_j)\|^2 \\ 
& \qquad + \eta^2 \mathbb{E}\| \nabla_{\boldsymbol{y}} F(\boldsymbol{x}^{k+1}, \boldsymbol{y}^{k+1}; \mathcal{D}_{F, j}) - \nabla_{\boldsymbol{y}} f(\boldsymbol{x}^{k+1}, \boldsymbol{y}^{k+1})\|^2 \\
& \overset{(ii)}{\le} (1 - \eta \mu)\mathbb{E}\| \boldsymbol{v}^{k+1, j} - \boldsymbol{v}(\boldsymbol{x}^{k+1}, \boldsymbol{y}^{k+1})\|^2 + \eta^2 \frac{M^2}{\mu^2} \frac{\sigma_{g, 2}^2}{D} + \eta^2 \frac{\sigma_f^2}{D_f},
\end{align*}
where $(i)$ follows from Lemma $\ref{lemma:2}$ and Assumptions 1, 2, and $(ii)$ follows from Assumption 3, and the fact that the sizes of $\mathcal{B}_j$ and $\mathcal{D}_{F, j}$ are $D$ and $D_f$, respectively.

Since $\boldsymbol{v}^{k+1} = \boldsymbol{v}^{k+1, J}$(see line 13 in Algorithm 2) and $\boldsymbol{v}^k = \boldsymbol{v}^{k+1, 0}$(see line 8 in Algorithm 2), by using induction, we have 

\begin{align} \label{eq33}    \nonumber
& \mathbb{E}\| \boldsymbol{v}^{k+1} - \boldsymbol{v}(\boldsymbol{x}^{k+1}, \boldsymbol{y}^{k+1})\|^2 \\   \nonumber
&\le (1 - \eta \mu)^J \mathbb{E}\| \boldsymbol{v}^k - \boldsymbol{v}(\boldsymbol{x}^{k+1}, \boldsymbol{y}^{k+1})\|^2 + \sum_{j=0}^{J-1}(1 - \eta \mu)^j \left(\eta^2 \frac{M^2}{\mu^2} \frac{\sigma_{g, 2}^2}{D} + \eta^2 \frac{\sigma_f^2}{Df}\right)\\   \nonumber
& \overset{(i)}{\le} (1 - \eta \mu)^J \mathbb{E}\| \boldsymbol{v}^k - \boldsymbol{v}(\boldsymbol{x}^{k+1}, \boldsymbol{y}^{k+1})\|^2 + \frac{1}{\eta \mu} \left(\eta^2 \frac{M^2}{\mu^2}\frac{\sigma_{g, 2}^2}{D} + \eta^2 \frac{\sigma_f^2}{D_f}\right)\\    \tag{33}
&\overset{(ii)}{\le} (1 - \eta \mu)^J (1 + r_w) \mathbb{E}\| \boldsymbol{v}^k - \boldsymbol{v}(\boldsymbol{x}^k, \boldsymbol{y}^k)\|^2\\   \nonumber
&\qquad + (1 - \eta \mu)^J \left(1 + \frac{1}{r_w}\right) \mathbb{E}\| \boldsymbol{v}(\boldsymbol{x}^k, \boldsymbol{y}^k) - \boldsymbol{v}(\boldsymbol{x}^{k+1}, \boldsymbol{y}^{k+1})\|^2    
 + \eta \frac{M^2}{\mu^3}\frac{\sigma_{g, 2}^2}{D} + \frac{\eta}{\mu}\frac{\sigma_f^2}{Df},   \nonumber
\end{align}
where $(i)$ holds since $\sum_{j=0}^{J-1}(1 - \eta \mu)^j \le {1}/{\eta \mu}$, and $(ii)$ holds since $r_w \|\boldsymbol{a}\|^2 + {1}/{r_w}\|\boldsymbol{b}\|^2 \ge 2 \boldsymbol{a}^\top \boldsymbol{b}$ for $\boldsymbol{a}$, $\boldsymbol{b}\in \mathbb{R}^q$.

For $\mathbb{E}\| \boldsymbol{v}(\boldsymbol{x}^k, \boldsymbol{y}^k) - \boldsymbol{v}(\boldsymbol{x}^{k+1}, \boldsymbol{y}^{k+1})\|^2$ in (\ref{eq33}), we have 
\begin{align}  \label{eq34}   \tag{34}
\mathbb{E}\| \boldsymbol{v}(\boldsymbol{x}^k, \boldsymbol{y}^k) - \boldsymbol{v}(\boldsymbol{x}^{k+1}, \boldsymbol{y}^{k+1})\|^2 &\overset{(i)}{\le} 2C_1^2 \mathbb{E}\| \boldsymbol{x}^k - \boldsymbol{x}^{k+1}\|^2 + 2 C_1^2 \mathbb{E}\| \boldsymbol{y}^k - \boldsymbol{y}^{k+1}\|^2\\   \nonumber
& \overset{(ii)}{=} 2 C_1^2 \alpha^2 \mathbb{E}\| \boldsymbol{h}_f^k \|^2 + 2 C_1^2 \mathbb{E}\| \boldsymbol{y}^k - \boldsymbol{y}^{k+1}\|^2, \nonumber
\end{align}  
where $(i)$ follows from Lemma \ref{lemma:1} and the fact that $\| \boldsymbol{a}+\boldsymbol{b}\|^2 \le 2 \|\boldsymbol{a}\|^2 + 2\| \boldsymbol{b}\|^2$ for $\boldsymbol{a}$, $\boldsymbol{b}\in \mathbb{R}^q$, and $(ii)$ follows from the iteration step in line 2 of Algorithm 2.

Since  
\begin{align}  \label{eq35}      \nonumber
&\mathbb{E}\| \boldsymbol{y}^k - \boldsymbol{y}^{k+1}\|^2 \\   \nonumber
& \overset{(i)}{\le} 2 \mathbb{E}\| \boldsymbol{y}^k - \boldsymbol{y}^*(\boldsymbol{x}^k)\|^2 + 4 \mathbb{E}\| \boldsymbol{y}^{k+1} - \boldsymbol{y}^*(\boldsymbol{x}^{k+1})\|^2 + 4 \mathbb{E}\|\boldsymbol{y}^*(\boldsymbol{x}^k) - \boldsymbol{y}^*(\boldsymbol{x}^{k+1})\|^2\\  \tag{35}
& \overset{(ii)}{\le} 2 \mathbb{E}\| \boldsymbol{y}^k - \boldsymbol{y}^*(\boldsymbol{x}^k)\|^2+ 4\mathbb{E}\| \boldsymbol{y}^{k+1} - \boldsymbol{y}^*(\boldsymbol{x}^{k+1})\|^2 + 4 \frac{L^2}{\mu^2} \mathbb{E}\|\boldsymbol{x}^k - \boldsymbol{x}^{k+1}\|^2\\  \nonumber
& \overset{(iii)}{\le} 2 \mathbb{E}\| \boldsymbol{y}^k - \boldsymbol{y}^*(\boldsymbol{x}^k)\|^2 + 4\mathbb{E}\| \boldsymbol{y}^{k+1} - \boldsymbol{y}^*(\boldsymbol{x}^{k+1})\|^2 + 4 \frac{L^2}{\mu^2}\alpha^2 \mathbb{E}\|\boldsymbol{h}_f^k\|^2,   \nonumber
\end{align}
where $(i)$ uses the inequality $\| \boldsymbol{a} + \boldsymbol{b}\|^2 \le 2 \|\boldsymbol{a} \|^2 + 2\|\boldsymbol{b}\|^2$ for $\boldsymbol{a}$, $\boldsymbol{b}\in \mathbb{R}^q$, $(ii)$ follows from Proposition 1, and $(iii)$ follows from the iteration step in line 2 of Algorithm 2.

Plugging (\ref{eq35}) into (\ref{eq34}), $\mathbb{E}\| \boldsymbol{v}(\boldsymbol{x}^k, \boldsymbol{y}^k) - \boldsymbol{v}(\boldsymbol{x}^{k+1}, \boldsymbol{y}^{k+1})\|^2$ is bounded as follows
\begin{align}  \label{eq36}  \tag{36}
& \mathbb{E}\| \boldsymbol{v}(\boldsymbol{x}^k, \boldsymbol{y}^k) - \boldsymbol{v}(\boldsymbol{x}^{k+1}, \boldsymbol{y}^{k+1})\|^2 \\   \nonumber
& \le \left(2+ \frac{8L^2}{\mu^2}\right)C_1^2\alpha^2 \mathbb{E}\|\boldsymbol{h}_f^k\|^2 + 4C_1^2 \mathbb{E}\| \boldsymbol{y}^k - \boldsymbol{y}^*(\boldsymbol{x}^k)\|^2+ 8 C_1^2\mathbb{E}\|\boldsymbol{y}^{k+1} - \boldsymbol{y}^*(\boldsymbol{x}^{k+1})\|^2.  \nonumber
\end{align}

Plugging (\ref{eq36}) into (\ref{eq33}), we have
\begin{align}  \label{eq37}   \nonumber
&\mathbb{E}\| \boldsymbol{v}^{k+1} - \boldsymbol{v}(\boldsymbol{x}^{k+1}, \boldsymbol{y}^{k+1})\|^2 \\   \nonumber
&\overset{(i)}{\le} (1 - \eta \mu)^J (1 + r_w)\mathbb{E}\| \boldsymbol{v}^k - \boldsymbol{v}(\boldsymbol{x}^k, \boldsymbol{y}^k)\|^2 + C_{22}\alpha^2 \mathbb{E}\|\boldsymbol{h}_f^k\|^2   \\   \nonumber
&\qquad + (1 - \eta \mu)^J \left(1 + \frac{1}{r_w}\right)4C_1^2 \mathbb{E}\|\boldsymbol{y}^k - \boldsymbol{y}^*(\boldsymbol{x}^k)\|^2  \\   \nonumber
 & \qquad + (1 - \eta \mu)^J \left(1 + \frac{1}{r_w}\right) 8 C_1^2 \mathbb{E}\| \boldsymbol{y}^{k+1} - \boldsymbol{y}^*(\boldsymbol{x}^{k+1})\|^2 + \eta \frac{M^2}{\mu^3}\frac{\sigma_{g, 2}^2}{D} + \frac{\eta}{\mu}\frac{\sigma_f^2}{D_f}\\   \nonumber
& \overset{(ii)}{\le} (1 - \eta \mu)^J (1 + r_w)^2 \mathbb{E}\| \boldsymbol{v}^k -  \boldsymbol{v}(\boldsymbol{x}^k, \boldsymbol{y}^*(\boldsymbol{x}^k))\|^2 \\    \tag{37}
&\qquad + (1 - \eta \mu)^J (1 + r_w)\left(1 + \frac{1}{r_w}\right)\mathbb{E}\| \boldsymbol{v}(\boldsymbol{x}^k, \boldsymbol{y}^*(\boldsymbol{x}^k)) - \boldsymbol{v}(\boldsymbol{x}^k, \boldsymbol{y}^k)\|^2  + C_{22}\alpha^2 \mathbb{E}\| \boldsymbol{h}_f^k\|^2\\   \nonumber
&\qquad + (1 - \eta \mu)^J \left(1 + \frac{1}{r_w}\right)4C_1^2 \mathbb{E}\| \boldsymbol{y}^k - \boldsymbol{y}^*(\boldsymbol{x}^k)\|^2 \\   \nonumber
&\qquad + (1 - \eta \mu)^J \left(1 + \frac{1}{r_w}\right)8C_1^2 \mathbb{E}\|\boldsymbol{y}^{k+1} - \boldsymbol{y}^*(\boldsymbol{x}^{k+1})\|^2 + \eta \frac{M^2}{\mu^3}\frac{\sigma_{g, 2}^2}{D} + \frac{\eta}{\mu}\frac{\sigma_f^2}{D_f}\\   \nonumber
&\overset{(iii)}{\le} C_{21} \mathbb{E}\| \boldsymbol{v}^k - \boldsymbol{v}(\boldsymbol{x}^k,\boldsymbol{y}^*(\boldsymbol{x}^k))\|^2 + C_{22}\alpha^2 \mathbb{E}\|\boldsymbol{h}_f^k\|^2 + C_{23} \mathbb{E}\|\boldsymbol{y}^k - \boldsymbol{y}^*(\boldsymbol{x}^k)\|^2 \\   \nonumber
&\qquad + \eta \frac{M^2}{\mu^3}\frac{\sigma_{g, 2}^2}{D} + \frac{\eta}{\mu}\frac{\sigma_f^2}{D_f} + (1 - \eta \mu)^J (1 + \frac{1}{r_w})8C_1^2 \mathbb{E}\|\boldsymbol{y}^{k+1} - \boldsymbol{y}^*(\boldsymbol{x}^{k+1})\|^2, \nonumber
\end{align}
where $(i)$ uses the definition of $C_{22}$, $(ii)$ holds by noticing that $r_w \|\boldsymbol{a}\|^2 + {1}/{r_w}\|\boldsymbol{b}\|^2 \ge 2 \boldsymbol{a}^\top \boldsymbol{b}$ for $\boldsymbol{a}$, $\boldsymbol{b}\in \mathbb{R}^q$, and $(iii)$ follows from Lemma \ref{lemma:1} and the definitions of $C_{21}$ and $C_{23}$. 

Plugging (\ref{eq37}) into (\ref{eq32}), we have 
\begin{align*} 
&\mathbb{E}\| \boldsymbol{v}^{k+1} - \boldsymbol{v}(\boldsymbol{x}^{k+1},\boldsymbol{y}^*(\boldsymbol{x}^{k+1}))\|^2 \\
& \overset{(i)}{\le} (1+ r_w)C_{21} \mathbb{E}\|\boldsymbol{v}^k - \boldsymbol{v}(\boldsymbol{x}^k, \boldsymbol{y}^*(\boldsymbol{x}^k))\|^2
\qquad + (1 + r_w)C_{22}\alpha^2 \mathbb{E}\|\boldsymbol{h}_f^k\|^2\\
&\qquad  + (1 + r_w)C_{23}\mathbb{E}\|\boldsymbol{y}^k - \boldsymbol{y}^*(\boldsymbol{x}^k)\|^2 
+ C_{24}\mathbb{E}\|\boldsymbol{y}^{k+1} - \boldsymbol{y}^*(\boldsymbol{x}^{k+1})\|^2\\
&\qquad + (1 + r_w)\eta \frac{M^2}{\mu^3}\frac{\sigma_{g, 2}^2}{D} + (1 + r_w)\frac{\eta}{\mu}\frac{\sigma_f^2}{D_f},
\end{align*}
where $(i)$ uses the definition of $C_{24}$. Then, the proof is complete.
\end{IEEEproof}

Define Lyapunov function
\begin{equation}  \label{eq38}  \tag{38}
H^{k} = \Phi(\boldsymbol{x}^{k}) + \rho_1 \| \boldsymbol{y}^{k} - \boldsymbol{y}^*(\boldsymbol{x}^{k})\|^2 + \rho_2 \| \boldsymbol{v}^{k} - \boldsymbol{v}(\boldsymbol{x}^{k}, \boldsymbol{y}^*(\boldsymbol{x}^{k}))\|^2,
\end{equation}
where $\rho_1$, $\rho_2$ are positive numbers. Then, based on the above lemmas, we analyze the difference between two Lyapunov functions $H^{k+1}$ and $H^k$ in the following lemma.

\begin{lemma}  \label{lemma:8}
Suppose Assumptions 1, 2, 3 hold. Define $\boldsymbol{v}(\boldsymbol{x}, \boldsymbol{y}) = [\nabla_{\boldsymbol{y}}^2 g(\boldsymbol{x}, \boldsymbol{y})]^{-1} \nabla_{\boldsymbol{y}} f(\boldsymbol{x}, \boldsymbol{y})$ and $\bar{\boldsymbol{h}}_f^k = \mathbb{E}[\boldsymbol{h}_f^k|\mathcal{F}_k^{'}]$. Let $\|\boldsymbol{v}^0\|\le M/\mu$, $\eta <{1}/{L}$, $\beta \le {2}/{(\mu + L)}$, and $\rho_1$, $\rho_2$ in (\ref{eq38}), $r_v$ in Lemma \ref{lemma:6}, and $r_w$ in Lemma \ref{lemma:7} be arbitrary positive numbers. Then for $H^k$ defined in (\ref{eq38}), we have 
\begin{align}  \label{eq39}  \nonumber
\mathbb{E}[H^{k+1}] - \mathbb{E}[H^k] &\le - \frac{\alpha}{2}\mathbb{E}\| \nabla \Phi(\boldsymbol{x}^k)\|^2 - L_{11}\mathbb{E}\| \bar{\boldsymbol{h}}_f^k \|^2  + L_{12}\mathbb{E}\| \boldsymbol{y}^k - \boldsymbol{y}^*(\boldsymbol{x}^k)\|^2  \\   \nonumber
&\qquad + L_{13}\mathbb{E}\| \boldsymbol{v}^k - \boldsymbol{v}(\boldsymbol{x}^k, \boldsymbol{y}^*(\boldsymbol{x}^k))\|^2 + \left(\frac{\alpha}{2} - L_{11} + \rho_2(1 + r_w)\frac{\eta}{\mu} \right)\frac{\sigma_f^2}{D_f} \\  \tag{39}
& \qquad + \left(\frac{\alpha}{2}-L_{11}\right) \frac{2M^2}{\mu^2}\frac{\sigma_{g, 1}^2}{D_g}  + (\rho_1 + C_{24}\rho_2)\frac{\mu + L}{2\mu L}\beta \frac{1}{S}\sigma_g^2 + \rho_2 (1 + r_w)\eta \frac{M^2}{\mu^3}\frac{\sigma_{g, 2}^2}{D} ,  \nonumber
\end{align}
where 
\begin{align*}
L_{11}& = \frac{\alpha}{2} - \frac{L_{\Phi}}{2}\alpha^2 - 2\left((\rho_1 + C_{24}\rho_2)(1 - \rho_y)^T\left( 1 + \frac{1}{r_v}\right) \frac{L^2}{\mu^2}\alpha^2 + \rho_2 (1 +r_w)C_{22}\alpha^2\right),\\
L_{12}& = \alpha \bar{C}_1^2 + (\rho_1 + C_{24}\rho_2)(1 - \rho_y)^T (1 + r_v) + \rho_2(1 + r_w)C_{23} - \rho_1,\\
L_{13}& = \alpha L^2 + 2L^2\left(\frac{\alpha}{2} - L_{11}\right) + \rho_2(1 + r_w)C_{21} - \rho_2,
\end{align*}
$\rho_y = 2 \beta {\mu L}/{(\mu + L)}$ is defined in Lemma \ref{lemma:6}, $\bar{C}_1 = L + {M \tau}/{\mu}$ is defined in Lemma \ref{lemma:4}, $L_{\Phi}$ is definied in Proposition 1, $C_{21}$, $C_{22}$, $C_{23}$, and $C_{24}$ are definied in Lemma \ref{lemma:7}.

\end{lemma}

\begin{IEEEproof}
Based on the definition of $H^k$ in (\ref{eq38}), taking a full expectation, we have 
\begin{align}  \label{eq40}  \nonumber
\mathbb{E}[H^{k+1}] - \mathbb{E}[H^k] &= \mathbb{E}[\Phi(\boldsymbol{x}^{k+1})] - \mathbb{E}[\Phi(\boldsymbol{x}^k)] + \rho_1 \mathbb{E}\| \boldsymbol{y}^{k+1} - \boldsymbol{y}^*(\boldsymbol{x}^{k+1})\|^2\\   \nonumber
& \qquad + \rho_2 \mathbb{E}\| \boldsymbol{v}^{k+1} - \boldsymbol{v}(\boldsymbol{x}^{k+1}, \boldsymbol{y}^*(\boldsymbol{x}^{k+1}))\|^2  \\   \tag{40}
& \qquad - \rho_1 \mathbb{E}\| \boldsymbol{y}^k - \boldsymbol{y}^*(\boldsymbol{x}^k)\|^2 - \rho_2 \mathbb{E}\| \boldsymbol{v}^k - \boldsymbol{v}(\boldsymbol{x}^k, \boldsymbol{y}^*(\boldsymbol{x}^k))\|^2.  \nonumber
\end{align}

By combining (\ref{eq40}) with the bound of $\mathbb{E}\| \boldsymbol{v}^{k+1} - \boldsymbol{v}(\boldsymbol{x}^{k+1}, \boldsymbol{y}^*(\boldsymbol{x}^{k+1}))\|^2$ in (\ref{eq31}), we have 
\begin{align}  \label{eq41}  \nonumber
\mathbb{E}[H^{k+1}] - \mathbb{E}[H^k] &\le \mathbb{E}[\Phi(\boldsymbol{x}^{k+1})] - \mathbb{E}[\Phi(\boldsymbol{x}^k)] + (\rho_1 + C_{24}\rho_2)\mathbb{E}\| \boldsymbol{y}^{k+1} - \boldsymbol{y}^*(\boldsymbol{x}^{k+1})\|^2 \\   \nonumber
&\qquad  + \left( \rho_2(1 + r_w)C_{21} - \rho_2 \right)\mathbb{E}\| \boldsymbol{v}^k - \boldsymbol{v}(\boldsymbol{x}^k, \boldsymbol{y}^*(\boldsymbol{x}^k))\|^2\\   \nonumber
&\qquad + \left( \rho_2 (1 + r_w)C_{23} - \rho_1 \right) \mathbb{E}\| \boldsymbol{y}^k - \boldsymbol{y}^*(\boldsymbol{x}^k)\|^2 \\      \tag{41}
&\qquad  + \rho_2 (1 + r_w)C_{22}\alpha^2 \mathbb{E}\| \boldsymbol{h}_f^k\|^2 + \rho_2 (1 + r_w)\eta \frac{M^2}{\mu^3}\frac{\sigma_{g, 2}^2}{D} + \rho_2(1 + r_w)\frac{\eta}{\mu}\frac{\sigma_f^2}{D_f}.\nonumber
\end{align}

By combining (\ref{eq41}) with the bound of $\mathbb{E}\| \boldsymbol{y}^{k+1} - \boldsymbol{y}^*(\boldsymbol{x}^{k+1})\|^2$ in (\ref{eq25}), we have
\begin{align}  \label{eq42}  \nonumber
\mathbb{E}[H^{k+1}] - \mathbb{E}[H^k] &\le \mathbb{E}[\Phi(\boldsymbol{x}^{k+1})] - \mathbb{E}[\Phi(\boldsymbol{x}^k)] \\   \nonumber
&\qquad  + \left( (\rho_1 + C_{24}\rho_2)(1 - \rho_y)^T (1 + r_v)+ \rho_2(1 + r_w)C_{23} - \rho_1 \right) \mathbb{E}\| \boldsymbol{y}^k - \boldsymbol{y}^*(\boldsymbol{x}^k)\|^2  \\   \nonumber
&\qquad  + \left( (\rho_1 + C_{24}\rho_2)(1 - \rho_y)^T \left( 1 + \frac{1}{r_v}\right) \frac{L^2}{\mu^2}\alpha^2 + \rho_2(1 + r_w)C_{22}\alpha^2 \right) \mathbb{E}\| \boldsymbol{h}_f^k \|^2 \\  \tag{42}
&\qquad  + \left( \rho_2(1 + r_w)C_{21} - \rho_2 \right) \mathbb{E}\| \boldsymbol{v}^k - \boldsymbol{v}(\boldsymbol{x}^k, \boldsymbol{y}^*(\boldsymbol{x}^k))\|^2 \\     \nonumber
&\qquad + (\rho_1 + C_{24}\rho_2)\frac{\mu + L}{2\mu L} \beta \frac{1}{S}\sigma_g^2 + \rho_2 (1 + r_w)\eta \frac{M^2}{\mu^3}\frac{\sigma_{g, 2}^2}{D} + \rho_2(1 + r_w)\frac{\eta}{\mu}\frac{\sigma_f^2}{D_f}.  \nonumber
\end{align}

By combining (\ref{eq42}) with the bound of $\mathbb{E}[\Phi(\boldsymbol{x}^{k+1})] - \mathbb{E}[\Phi(\boldsymbol{x}^k)]$ in (\ref{eq19}), we have
\begin{align} \label{eq43}  \nonumber
&\mathbb{E}[H^{k+1}] - \mathbb{E}[H^k] \\   \nonumber
& \le - \frac{\alpha}{2} \mathbb{E}\| \nabla \Phi(\boldsymbol{x}^k)\|^2- \left( \frac{\alpha}{2} - \frac{L_{\Phi}}{2}\alpha^2\right) \mathbb{E}\| \bar{\boldsymbol{h}}_f^k \|^2 + \frac{\alpha}{2}\mathbb{E}\| \nabla \Phi(\boldsymbol{x}^k) - \bar{\boldsymbol{h}}_f^k \|^2 + \frac{L_{\Phi}}{2}\alpha^2 \mathbb{E}\| \boldsymbol{h}_f^k - \bar{\boldsymbol{h}}_f^k\|^2 \\   \nonumber
& \qquad + \left( (\rho_1 + C_{24}\rho_2)(1 - \rho_y)^T (1 + r_v)+ \rho_2(1 + r_w)C_{23} - \rho_1 \right) \mathbb{E}\| \boldsymbol{y}^k - \boldsymbol{y}^*(\boldsymbol{x}^k)\|^2 \\   \nonumber
& \qquad + \left( (\rho_1 + C_{24}\rho_2)(1 - \rho_y)^T \left( 1 + \frac{1}{r_v}\right) \frac{L^2}{\mu^2}\alpha^2 + \rho_2(1 + r_w)C_{22}\alpha^2 \right) \mathbb{E}\| \boldsymbol{h}_f^k \|^2 \\  \nonumber
& \qquad + \left( \rho_2(1 + r_w)C_{21} - \rho_2 \right) \mathbb{E}\| \boldsymbol{v}^k - \boldsymbol{v}(\boldsymbol{x}^k, \boldsymbol{y}^*(\boldsymbol{x}^k))\|^2\\   \nonumber
& \qquad + (\rho_1 + C_{24}\rho_2)\frac{\mu + L}{2\mu L} \beta \frac{1}{S}\sigma_g^2 + \rho_2 (1 + r_w)\eta \frac{M^2}{\mu^3}\frac{\sigma_{g, 2}^2}{D} + \rho_2(1 + r_w)\frac{\eta}{\mu}\frac{\sigma_f^2}{D_f}\\  \tag{43}
& \overset{(i)}{\le} -\frac{\alpha}{2}\mathbb{E}\| \nabla \Phi(\boldsymbol{x}^k)\|^2- L_{11}\mathbb{E}\| \bar{\boldsymbol{h}}_f^k \|^2 + \frac{\alpha}{2}\mathbb{E}\| \nabla \Phi(\boldsymbol{x}^k)- \bar{\boldsymbol{h}}_f^k\|^2 + \left(\frac{\alpha}{2} - L_{11}\right) \mathbb{E}\|\boldsymbol{h}_f^k - \bar{\boldsymbol{h}}_f^k \|^2 \\  \nonumber
&\qquad  + \left( (\rho_1 + C_{24}\rho_2)(1 - \rho_y)^T (1 + r_v)+ \rho_2(1 + r_w)C_{23} - \rho_1 \right) \mathbb{E}\| \boldsymbol{y}^k - \boldsymbol{y}^*(\boldsymbol{x}^k)\|^2 \\   \nonumber
& \qquad + \left( \rho_2(1 + r_w)C_{21} - \rho_2 \right) \mathbb{E}\| \boldsymbol{v}^k - \boldsymbol{v}(\boldsymbol{x}^k, \boldsymbol{y}^*(\boldsymbol{x}^k))\|^2 \\   \nonumber
&\qquad + (\rho_1 + C_{24}\rho_2)\frac{\mu + L}{2\mu L} \beta \frac{1}{S}\sigma_g^2 + \rho_2 (1 + r_w)\eta \frac{M^2}{\mu^3}\frac{\sigma_{g, 2}^2}{D} + \rho_2(1 + r_w)\frac{\eta}{\mu}\frac{\sigma_f^2}{D_f},  \nonumber
\end{align}
where $(i)$ follows from the definition of $L_{11}$ and the fact that $\mathbb{E}\|\boldsymbol{h}_f^k\|^2 \le 2 \mathbb{E}\| \bar{\boldsymbol{h}}_f^k\|^2 + 2 \mathbb{E}\| \boldsymbol{h}_f^k - \bar{\boldsymbol{h}}_f^k\|^2$.

By combining (\ref{eq43}) with the bound of $\mathbb{E}\| \nabla \Phi(\boldsymbol{x}^k) - \bar{\boldsymbol{h}}_f^k\|^2$ in (\ref{eq20}), we have
\begin{align} \label{eq44}   \nonumber
&\mathbb{E}[H^{k+1}] - \mathbb{E}[H^k] \\   \nonumber
& \overset{(i)}{\le} -\frac{\alpha}{2}\mathbb{E}\| \nabla \Phi(\boldsymbol{x}^k)\|^2 - L_{11}\mathbb{E}\| \bar{\boldsymbol{h}}_f^k \|^2   \nonumber
 + L_{12}\mathbb{E}\| \boldsymbol{y}^k - \boldsymbol{y}^*(\boldsymbol{x}^k)\|^2  \\    \tag{44}
& \qquad + (\alpha L^2 + \rho_2(1 + r_w)C_{21}-\rho_2) \mathbb{E} \|\boldsymbol{v}^k - \boldsymbol{v}(\boldsymbol{x}^k, \boldsymbol{y}^*(\boldsymbol{x}^k))\|^2   \nonumber
+ \left(\frac{\alpha}{2} - L_{11}\right) \mathbb{E}\| \boldsymbol{h}_f^k - \bar{\boldsymbol{h}}_f^k \|^2\\   \nonumber
& \qquad + (\rho_1 + C_{24}\rho_2)\frac{\mu + L}{2\mu L} \beta \frac{1}{S}\sigma_g^2 + \rho_2 (1 + r_w)\eta \frac{M^2}{\mu^3}\frac{\sigma_{g, 2}^2}{D} + \rho_2(1 + r_w)\frac{\eta}{\mu}\frac{\sigma_f^2}{D_f},  \nonumber
\end{align}
where $(i)$ follows from the definiton of $L_{12}$.

By combining (\ref{eq44}) with the bound of $\mathbb{E}\|\boldsymbol{h}_f^k - \bar{\boldsymbol{h}}_f^k\|^2$ in (\ref{eq22}), and notice that ${\alpha}/{2}-L_{11}>0$, we have 

\begin{align*}  
\mathbb{E}[H^{k+1}] - \mathbb{E}[H^k] &\overset{(i)}{\le} - \frac{\alpha}{2}\mathbb{E}\| \nabla \Phi(\boldsymbol{x}^k)\|^2 - L_{11}\mathbb{E}\| \bar{\boldsymbol{h}}_f^k \|^2+ L_{12}\mathbb{E}\| \boldsymbol{y}^k - \boldsymbol{y}^*(\boldsymbol{x}^k)\|^2  \\
&\qquad + L_{13}\mathbb{E}\| \boldsymbol{v}^k - \boldsymbol{v}(\boldsymbol{x}^k, \boldsymbol{y}^*(\boldsymbol{x}^k))\|^2 + \left(\frac{\alpha}{2} - L_{11} + \rho_2(1 + r_w)\frac{\eta}{\mu}\right)\frac{\sigma_f^2}{D_f} \\
&\qquad + \left(\frac{\alpha}{2}-L_{11}\right) \frac{2M^2}{\mu^2}\frac{\sigma_{g, 1}^2}{D_g} 
+ (\rho_1 + C_{24}\rho_2)\frac{\mu + L}{2\mu L}\beta \frac{1}{S}\sigma_g^2 + \rho_2 (1 + r_w)\eta \frac{M^2}{\mu^3}\frac{\sigma_{g, 2}^2}{D},
\end{align*}
where $(i)$ follows from the definition of $L_{13}$.
\end{IEEEproof}
\subsection{Proof of Theorem 1}
For Theorem 1, since $\eta$, $\beta$, $\rho_1$, $\rho_2$, $r_v$, $r_w$, and $\boldsymbol{v}^0$ satisfy the conditions of Lemma \ref{lemma:8}, we have 
\begin{align}  \label{eq45}  \nonumber
\mathbb{E}[H^{k+1}] - \mathbb{E}[H^k] &\le - \frac{\alpha}{2}\mathbb{E}\| \nabla \Phi(\boldsymbol{x}^k)\|^2 - L_{11}\mathbb{E}\| \bar{\boldsymbol{h}}_f^k \|^2 + L_{12}\mathbb{E}\| \boldsymbol{y}^k - \boldsymbol{y}^*(\boldsymbol{x}^k)\|^2  \\ \nonumber
&\qquad + L_{13}\mathbb{E}\| \boldsymbol{v}^k - \boldsymbol{v}(\boldsymbol{x}^k, \boldsymbol{y}^*(\boldsymbol{x}^k))\|^2 + \left(\frac{\alpha}{2} - L_{11}\right)\frac{\sigma_f^2}{D_f} + \left(\frac{\alpha}{2}-L_{11}\right) \frac{2M^2}{\mu^2}\frac{\sigma_{g, 1}^2}{D_g} \\   \tag{45}
&\qquad  + (\rho_1 + C_{24}\rho_2)\frac{\mu + L}{2\mu L}\beta \frac{1}{S}\sigma_g^2 + \rho_2 (1 + r_w)\eta \frac{M^2}{\mu^3}\frac{\sigma_{g, 2}^2}{D} + \rho_2(1 + r_w)\frac{\eta}{\mu}\frac{\sigma_f^2}{D_f},  \nonumber
\end{align}
where $L_{11}$, $L_{12}$, and $L_{13}$ are defined in Lemma \ref{lemma:8}.

In the following, we will prove $L_{11}>0$, $L_{12}<0$, and $L_{13}<0$.

Firstly, we prove $L_{11}>0$. Based on the definition of $L_{11}$, we have 
\begin{align*}
\hspace{0.7cm}L_{11} &= \frac{\alpha}{2} - \frac{L_{\Phi}}{2}\alpha^2 - 2\left((\rho_1 + C_{24}\rho_2)(1 - \rho_y)^T\left(1 + \frac{1}{r_v}\right)\frac{L^2}{\mu^2}\alpha^2 + \rho_2(1 + r_w)C_{22}\alpha^2\right)\\
&\overset{(i)}{>} \frac{\alpha}{2} - L_{31}\alpha^2 \overset{(ii)}{\ge} 0,
\end{align*}
where $(i)$ holds since $r_v = {2\mu L}/{(\mu^2 + L^2)}\le 1$ and 
\begin{equation*}
\frac{L_{\Phi}}{2} + 2\left((\rho_1 + C_{24}\rho_2)(1 - \rho_y)^T\left(1 + \frac{1}{r_v}\right)\frac{L^2}{\mu^2} + \rho_2(1 + r_w)C_{22}\right) < L_{31}, 
\end{equation*}
and $(ii)$ uses the definition of $\alpha$.

Secondly, we prove $L_{12}<0$. Based on the definition of $L_{12}$, we have 
\begin{equation}  \label{eq46}  \tag{46}
L_{12} = \alpha\bar{C}_1^2 + \left(\rho_1 + C_{24}\rho_2\right)(1 - \rho_y)^T (1 + r_v) + \rho_2 (1 + r_w)C_{23} - \rho_1. 
\end{equation}
Notice that 
\begin{align*}
& \left(\rho_1 + C_{24}\rho_2\right)(1 - \rho_y)^T(1 + r_v) + \rho_2(1 + r_w)C_{23} - \rho_1 \\
& \overset{(i)}{<} \rho_1 r_v + C_{24}\rho_2(1 + r_v) + \rho_2 (1 + r_w)C_{23} - (\rho_1 + C_{24}\rho_2)\rho_y (1 + r_v) \\
& \overset{(ii)}{=} \rho_1 \left( r_v - \rho_y (1 + r_v)\right) + \rho_2 \left( C_{24}(1 + r_v) + (1 + r_w)C_{23} - C_{24}\rho_y (1 + r_v)\right)\\
& \overset{(iii)}{<}-\rho_2 \left((1 + r_w)C_{23} + 2C_{24}\right),
\end{align*}
where $(i)$ follows from $(1 - \rho_y)^T \le 1 - \rho_y$, $(ii)$ is obtained by rearranging the terms, and $(iii)$ follows from the definition of $\rho_1$, and the fact that $r_v \le 1$ and $C_{24}\rho_y (1 + r_v) >0$. Then, for (\ref{eq46}), we have 
\begin{align*}
L_{12} &< \alpha \bar{C}_1^2 - \rho_2 \left( (1 + r_w)C_{23} + 2 C_{24}\right) \overset{(i)}{\le} 0,
\end{align*}
where $(i)$ uses the definition of $\alpha$.

Finally, we prove $L_{13}<0$. Based on the definition of $L_{13}$, we have 
\begin{align} \label{eq47}  \nonumber
L_{13} &= \alpha L^2 + 2L^2 \left(\frac{\alpha}{2} - L_{11}\right) + \rho_2(1 + r_w)C_{21} - \rho_2 \\   \tag{47}
& \overset{(i)}{<} 2 \alpha L^2 + (1-\eta \mu)(1 + r_w)^3 \rho_2 - \rho_2, \\   \nonumber
\end{align}
where $(i)$ holds since $L_{11}>0$, and 
\begin{equation*}
C_{21} = (1 - \eta \mu)^J(1 + r_w)^2 < (1 - \eta \mu)(1 + r_w)^2.
\end{equation*}
Furthermore, since 
\begin{align*}
(1 - \eta \mu)(1 + r_w)^3 -1 &= (1 - \eta \mu)(1 + 3r_w + 3r_w^2 + r_w^3) - 1\\
& = (3r_w + 3r_w^2 + r_w^3) - \eta \mu  (1 + 3r_w + 3r_w^2 + r_w^3) \\
& \overset{(i)}{<} -\frac{1}{2}\eta \mu (1 + 3r_w + 3r_w^2 + r_w^3) < - \frac{\eta \mu}{2},
\end{align*}
and $(i)$ follows from the definition of $r_w$ and the fact that $(1 - {\eta \mu}/{2})(3r_w + 3r_w^2 + r_w^3)< {7(2 - \eta \mu)r_w}/{2} ={\eta \mu}/{2}$. Then, for (\ref{eq47}), using the definition of $\alpha$, we have 
\begin{align*}
L_{13}&< 2 \alpha L^2 - \frac{\rho_2}{2}\eta \mu \le 0.
\end{align*}

Since $L_{11}>0$, $L_{12}< 0$, and $L_{13}<0$, for (\ref{eq45}), we have 
\begin{align}  \label{eq48}  \nonumber
\mathbb{E}[H^{k+1}] - \mathbb{E}[H^k] &\le- \frac{\alpha}{2}\mathbb{E}\| \nabla \Phi(\boldsymbol{x}^k)\|^2 \\  \tag{48}
&\qquad  + \left(\frac{\alpha}{2} + \rho_2(1 + r_w)\frac{\eta}{\mu}\right)\frac{\sigma_f^2}{D_f} +  \frac{M^2}{\mu^2}\alpha \frac{\sigma_{g, 1}^2}{D_g} \\   \nonumber
& \qquad + (\rho_1 + C_{24}\rho_2)\frac{\mu + L}{2\mu L}\beta \frac{1}{S}\sigma_g^2 + \rho_2 (1 + r_w)\eta \frac{M^2}{\mu^3}\frac{\sigma_{g, 2}^2}{D}.  \nonumber
\end{align}

Telescoping (\ref{eq48}) over $k$ from $0$ to $K-1$ and multiplying both sides by ${2}/{(K\alpha)}$, we have 
\begin{align} \label{eq49}  \nonumber
\frac{1}{K} \sum_{k=0}^{K-1} \mathbb{E}\| \nabla \Phi(\boldsymbol{x}^k)\|^2 &\le \frac{2}{K\alpha}\left(\mathbb{E}[H^0] - \mathbb{E}[H^K]\right) + \frac{2}{\alpha}\left(\frac{\alpha}{2} + \rho_2(1 + r_w)\frac{\eta}{\mu}\right)\frac{\sigma_f^2}{D_f} + \frac{2 M^2}{\mu^2} \frac{\sigma_{g, 1}^2}{Dg}\\   \nonumber
&\qquad  + \frac{1}{\alpha}(\rho_1 + C_{24}\rho_2)\frac{\mu + L}{\mu L}\beta \frac{1}{S}\sigma_g^2  + \frac{2}{\alpha}\rho_2(1 + r_w) \eta \frac{M^2}{\mu^3}\frac{\sigma_{g, 2}^2}{D}\\  \nonumber
& \overset{(i)}{\le} \frac{2}{K\alpha}\left(\mathbb{E}[H^0] - \min\limits_{\boldsymbol{x}}\Phi(\boldsymbol{x})\right) + \frac{2}{\alpha}\left(\frac{\alpha}{2} + \rho_2(1 + r_w)\frac{\eta}{\mu}\right)\frac{\sigma_f^2}{D_f} + \frac{2 M^2}{\mu^2} \frac{\sigma_{g, 1}^2}{Dg}\\  \tag{49}
& \qquad + \frac{1}{\alpha}(\rho_1 + C_{24}\rho_2)\frac{\mu + L}{\mu L}\beta \frac{1}{S}\sigma_g^2  + \frac{2}{\alpha}\rho_2(1 + r_w) \eta \frac{M^2}{\mu^3}\frac{\sigma_{g, 2}^2}{D}, \nonumber
\end{align}
where $(i)$ holds since for $H^k$ in (\ref{eq38}), we have $H^k \ge \Phi(\boldsymbol{x}^k)$ for any $k\in \{0, \ldots, K-1\}$.

Under the setting of Theorem 1, $r_w = \mathcal{O}(\kappa^{-1})$, $C_{22} = \mathcal{O}(\kappa^7)$, $C_{23}= \mathcal{O}(\kappa^5)$, $C_{24}= \mathcal{O}(\kappa^5)$, $\rho_y(1 + r_v)-r_v = {\mu L}/{(L^2 + \mu^2)} = \mathcal{O}(\kappa^{-1})$, $\rho_1 = \mathcal{O}(\kappa^2)$, $L_{31} = \mathcal{O}(\kappa^5)$, $\alpha = \mathcal{O}(\kappa^{-5})$. Then, for (\ref{eq49}), we have 
\begin{equation}  \label{eq50}
\frac{1}{K} \sum_{k=0}^{K-1} \mathbb{E}\| \nabla \Phi(\boldsymbol{x}^k)\|^2 \le \mathcal{O}\left(\frac{\kappa^5}{K} +\kappa^2 \frac{\sigma_f^2}{D_f} + \kappa^2 \frac{\sigma_{g, 1}^2}{D_g} + \kappa^8 \frac{1}{S}\sigma_g^2 + \kappa^4 \frac{\sigma_{g, 2}^2}{D}\right). \tag{50}
\end{equation}

To reach an $\epsilon$-accurate stationary point, i.e., $\mathbb{E}\left\| \nabla \Phi(\bar{\boldsymbol{x}})\right\|^2 \le \epsilon$ with $\bar{\boldsymbol{x}}$ chosen from $\boldsymbol{x}^0, \ldots, \boldsymbol{x}^{K-1}$ uniformly at random, it suffices to choose $K= \mathcal{O}(\kappa^5 \epsilon^{-1})$, $D_f = \mathcal{O}(\kappa^2 \epsilon^{-1})$, $D_g = \mathcal{O}(\kappa^2 \epsilon^{-1})$, $S= \mathcal{O}(\kappa^8 \epsilon^{-1})$, $D= \mathcal{O}(\kappa^4 \epsilon^{-1})$. Then, the gradient complexities are given by Gc$(F, \epsilon) = K(J+1)D_f = \mathcal{O}(\kappa^7 \epsilon^{-2}J)$, Gc$(G, \epsilon)= KTS = \mathcal{O}(\kappa^{13}\epsilon^{-2}T)$; the Jacobian and Hessian vector product complexities are given by Jv$(G, \epsilon)= KD_g = \mathcal{O}(\kappa^7 \epsilon^{-2})$ and Hv$(G, \epsilon) = KJD = \mathcal{O}(\kappa^9 \epsilon^{-2}J)$, respectively. In total, the computational complexity is $\mathcal{O}(\kappa^{13}\epsilon^{-2}\max\{T, J\})$. Then, the proof is complete.

\subsection{Proof of Theorem 2}

For Theorem 2, since $\eta$, $\beta$, $\rho_1$, $\rho_2$, $r_v$, $r_w$, and $\boldsymbol{v}^0$ satisfy the conditions of Lemma $\ref{lemma:8}$, we have 
\begin{align}  \label{eq51}  \nonumber
\mathbb{E}[H^{k+1}] - \mathbb{E}[H^k] &\le - \frac{\alpha}{2}\mathbb{E}\| \nabla \Phi(\boldsymbol{x}^k)\|^2 - L_{11}\mathbb{E}\| \bar{\boldsymbol{h}}_f^k \|^2 + L_{12}\mathbb{E}\| \boldsymbol{y}^k - \boldsymbol{y}^*(\boldsymbol{x}^k)\|^2  \\ \nonumber
&\qquad + L_{13}\mathbb{E}\|\boldsymbol{v}^k - \boldsymbol{v}(\boldsymbol{x}^k, \boldsymbol{y}^*(\boldsymbol{x}^k))\|^2 + \left(\frac{\alpha}{2} - L_{11}\right)\frac{\sigma_f^2}{D_f} + \left(\frac{\alpha}{2}-L_{11}\right) \frac{2M^2}{\mu^2}\frac{\sigma_{g, 1}^2}{D_g} \\   \tag{51}
&\qquad + (\rho_1 + C_{24}\rho_2)\frac{\mu + L}{2\mu L}\beta \frac{1}{S}\sigma_g^2 + \rho_2 (1 + r_w)\eta \frac{M^2}{\mu^3}\frac{\sigma_{g, 2}^2}{D} + \rho_2(1 + r_w)\frac{\eta}{\mu}\frac{\sigma_f^2}{D_f},   \nonumber
\end{align}
where $L_{11}$, $L_{12}$, and $L_{13}$ are defined in Lemma $\ref{lemma:8}$.

In the following, we will prove $L_{11}\ge 0$, $L_{12}\le 0$, and $L_{13}\le 0$.

Firstly, we prove $L_{11}\ge 0$. Based on the definition of $L_{11}$, we have 
\begin{equation}   \label{eq52}   \tag{52}
L_{11} = \frac{\alpha}{2} - \frac{L_{\Phi}}{2}\alpha^2 - 2\left((\rho_1 + C_{24}\rho_2)(1 - \rho_y)^T\left(1 + \frac{1}{r_v}\right)\frac{L^2}{\mu^2}\alpha^2 + \rho_2(1 + r_w)C_{22}\alpha^2\right).
\end{equation}
Note that $\rho_y = {2\beta \mu L}/{(\mu + L)}$ and $\beta = 1/(2L)$. We have 
\begin{equation}   \label{eq53} \tag{53}
(1 - \rho_y)^T = \left(\frac{L}{\mu + L}\right)^T \le \frac{\rho_1}{8(\rho_1 + C_{24}\rho_2)}
\end{equation}
since $T \ge {\log(\rho_1/(8(\rho_1 + C_{24}\rho_2)))}/{\log(L/(\mu + L))}$. Then, for (\ref{eq52}), we have 
\begin{equation*}
L_{11} \overset{(i)}{\ge} \frac{\alpha}{2} - \frac{L_{\Phi}}{2}\alpha^2 - 2\left(\frac{\rho_1}{4} \frac{L^2}{\mu^2}\alpha^2 + \rho_2(1 + r_w)C_{22}\alpha^2\right) \overset{(ii)}{=} \frac{\alpha}{2} - L_{31}\alpha^2 \overset{(iii)}{\ge} 0,
\end{equation*}
where $(i)$ follows from the definition of $r_v$ and the inequality in (\ref{eq53}), $(ii)$ uses the definition of $L_{31}$, and $(iii)$ uses the definition of $\alpha$.

Secondly, we prove $L_{12}\le 0$. Based on the definition of $L_{12}$, we have 
\begin{align*}
L_{12} &= \alpha\bar{C}_1^2 + \left(\rho_1 + C_{24}\rho_2\right)(1 - \rho_y)^T (1 + r_v) + \rho_2 (1 + r_w)C_{23} - \rho_1  \\
& \overset{(i)}{\le} \alpha \bar{C}_1^2 + \frac{\rho_1}{4} - \rho_2(1 + r_w)C_{23} \overset{(ii)}{=}\alpha \bar{C}_1^2 - \frac{\rho_2}{2}(1 + r_w)C_{23} \overset{(iii)}{\le} 0,  
\end{align*}
where $(i)$ follows from the inequality in (\ref{eq53}), the fact that $r_v=1$, and the definition of $\rho_1$, $(ii)$ follows from the definition of $\rho_1$, and $(iii)$ uses the definition of $\alpha$.

Finally, we prove $L_{13}\le 0$. Based on the definition of $L_{13}$, we have 
\begin{align*}
L_{13} &= \alpha L^2 + 2L^2 \left(\frac{\alpha}{2} - L_{11}\right) + \rho_2(1 + r_w)C_{21} - \rho_2 \\
& \overset{(i)}{\le } 2 \alpha L^2 + 2\rho_2 C_{21} - \rho_2  \overset{(ii)}{\le} 2 \alpha L^2 - \frac{\rho_2}{2} \overset{(iii)}{\le} 0,
\end{align*}
where $(i)$ holds since $L_{11}\ge0$ and $r_w=1$, $(ii)$ holds since $C_{21}= (1 - \eta \mu)^J(1 + r_w)^2 \le \frac{1}{4}$ when $J \ge {\log(1/(4(1 + r_w)^2))}/{\log(1-\eta \mu)}$, and $(iii)$ uses the definition of $\alpha$.

Since $L_{11}\ge 0$, $L_{12}\le  0$, $L_{13}\le 0$, similar to the discussion in Theorem 1, we have  
\begin{align} \label{eq54}   \nonumber
\frac{1}{K} \sum_{k=0}^{K-1} \mathbb{E}\| \nabla \Phi(\boldsymbol{x}^k)\|^2 &\le \frac{2}{K\alpha}\left(\mathbb{E}[H^0] - \min\limits_{\boldsymbol{x}}\Phi(\boldsymbol{x})\right) + \frac{2}{\alpha}\left(\frac{\alpha}{2} + \rho_2(1 + r_w)\frac{\eta}{\mu}\right)\frac{\sigma_f^2}{D_f} \\   \tag{54}
&\qquad + \frac{2 M^2}{\mu^2} \frac{\sigma_{g, 1}^2}{Dg} + \frac{1}{\alpha}(\rho_1 + C_{24}\rho_2)\frac{\mu + L}{\mu L}\beta \frac{1}{S}\sigma_g^2  + \frac{2}{\alpha}\rho_2(1 + r_w) \eta \frac{M^2}{\mu^3}\frac{\sigma_{g, 2}^2}{D}.   \nonumber
\end{align}

Under the setting of Theorem 2, $C_{22} = \mathcal{O}(\kappa^6)$, $C_{23}= \mathcal{O}(\kappa^4)$, $C_{24}= \mathcal{O}(\kappa^4)$, $\rho_1 = \mathcal{O}(\kappa)$, $L_{31} = \mathcal{O}(\kappa^3)$, $\alpha = \mathcal{O}(\kappa^{-3})$. Then, for (\ref{eq54}), we have 
\begin{equation*} 
\frac{1}{K} \sum_{k=0}^{K-1} \mathbb{E}\| \nabla \Phi(\boldsymbol{x}^k)\|^2 \le \mathcal{O}\left(\frac{\kappa^3}{K} +\kappa \frac{\sigma_f^2}{D_f} + \kappa^2 \frac{\sigma_{g, 1}^2}{D_g} + \kappa^5 \frac{1}{S}\sigma_g^2 + \kappa^3 \frac{\sigma_{g, 2}^2}{D}\right).
\end{equation*}

To reach an $\epsilon$-accurate stationary point, i.e., $\mathbb{E}\left\| \nabla \Phi(\bar{\boldsymbol{x}})\right\|^2 \le \epsilon$ with $\bar{\boldsymbol{x}}$ chosen from $\boldsymbol{x}^0, \ldots, \boldsymbol{x}^{K-1}$ uniformly at random, it suffices to choose $K= \mathcal{O}(\kappa^3 \epsilon^{-1})$, $D_f = \mathcal{O}(\kappa \epsilon^{-1})$, $D_g = \mathcal{O}(\kappa^2 \epsilon^{-1})$, $S= \mathcal{O}(\kappa^5 \epsilon^{-1})$, $D= \mathcal{O}(\kappa^3 \epsilon^{-1})$, $T= \mathcal{O}(\kappa)$, $J=\mathcal{O}(\kappa)$. Then, the gradient complexities are given by Gc$(F, \epsilon) = K(J+1)D_f = \mathcal{O}(\kappa^5 \epsilon^{-2})$, Gc$(G, \epsilon)= KTS = \mathcal{O}(\kappa^{9}\epsilon^{-2})$; the Jacobian and Hessian vector product complexities are given by Jv$(G, \epsilon)= KD_g = \mathcal{O}(\kappa^5 \epsilon^{-2})$ and Hv$(G, \epsilon) = KJD = \mathcal{O}(\kappa^7 \epsilon^{-2})$, respectively. In total, the computational complexity is $ \mathcal{O}(\kappa^{9}\epsilon^{-2})$. Then, the proof is complete.

\section{Analysis of Hypergradient Estimation Methods}\label{appendix:2}

For the iterations in Algorithm 1, we define 
\begin{align*}
&\mathcal{F}_{k,1}^{'} = \sigma\{ \boldsymbol{x}^0, \ldots, \boldsymbol{y}^{0, T}, \ldots, \boldsymbol{x}^k, \ldots, \boldsymbol{y}^{k,T}\},\\ 
&\mathcal{F}_{k,2}^{'} = \sigma\{ \boldsymbol{x}^0, \ldots, \boldsymbol{y}^{0, T}, \ldots, \boldsymbol{v}^{0, J}, \ldots, \boldsymbol{x}^k, \ldots, \boldsymbol{y}^{k,T}\}, 
\end{align*}
where $\sigma\{\cdot\}$ denotes the $\sigma$-algebra generated by random variables.

\begin{lemma}  \label{lemma:9}
Apply Algorithm 1 to solve problem (2). Suppose Assumptions 1, 2, and 3 hold. let $\mu$, $L$ in Assumptions 1, 2 be the numbers satisfying $\mu < 1$, $L<1$. For $T$, $J$, $\beta$, $\eta$ in Algorithm 1, let $T$, $J$ be any positive integers, $\beta$, $\eta$ be the stepsizes satisfying $\beta = \eta < 1$, and $\prod_{i=J}^{J-1}(\cdot) = \mathbf{I}$, where $\mathbf{I}$ is the identity matarix. Then 
\begin{itemize}
\item For Stochatic BP, for $\boldsymbol{h}_f^k$ in (6) and $\bar{\boldsymbol{h}}_f^k:= \mathbb{E}[\boldsymbol{h}_f^k|\mathcal{F}_{k, 1}^{'}]$, if sample sets $\{\mathcal{D}_F, \mathcal{S}_t(t = 0, \ldots, T-1)\}$ are mutually independent, then, we have 
\begin{equation}  \label{eq55} \tag{55}
\bar{\boldsymbol{h}}_f^k = \nabla_{\boldsymbol{x}} f(\boldsymbol{x}^k, \boldsymbol{y}^{k, T}) - \mathbf{H}_T \nabla_{\boldsymbol{y}}f(\boldsymbol{x}^k, \boldsymbol{y}^{k, T})
\end{equation}
where 
\begin{equation}  \label{eq56} \tag{56}
\mathbf{H}_T= \eta \sum\limits_{t=0}^{T-1} \nabla_{\boldsymbol{x}}\nabla_{\boldsymbol{y}} g(\boldsymbol{x}^k, \boldsymbol{y}^{k, t}) \prod_{i=t+1}^{T-1}(\mathbf{I} - \eta \nabla_{\boldsymbol{y}}^2 g(\boldsymbol{x}^k, \boldsymbol{y}^{k, i}));
\end{equation}
\item For Stochastic NS, for $\boldsymbol{h}_f^k$ in (11) and $\bar{\boldsymbol{h}}_f^k:= \mathbb{E}[\boldsymbol{h}_f^k|\mathcal{F}_{k, 1}^{'}]$, if sample sets $\{\mathcal{D}_F, \mathcal{D}_G, \mathcal{B}_j(j=1, \ldots, J-1)\}$ are mutually independent, then, we have 
\begin{equation}  \label{eq57} \tag{57}
\bar{\boldsymbol{h}}_f^k = \nabla_{\boldsymbol{x}} f(\boldsymbol{x}^k, \boldsymbol{y}^{k, T}) - \mathbf{H}_J \nabla_{\boldsymbol{y}} f(\boldsymbol{x}^k, \boldsymbol{y}^{k, T})
\end{equation}
where 
\begin{equation*} 
\mathbf{H}_J = \eta \sum\limits_{t=0}^{J-1} \nabla_{\boldsymbol{x}}\nabla_{\boldsymbol{y}} g(\boldsymbol{x}^k, \boldsymbol{y}^{k, T}) \prod_{i=t+1}^{J-1}(\mathbf{I} - \eta \nabla_{\boldsymbol{y}}^2 g(\boldsymbol{x}^k, \boldsymbol{y}^{k,T}));
\end{equation*}
\item For SGD-based Estimation, for $\boldsymbol{h}_f^k$ in (14) and $\bar{\boldsymbol{h}}_f^k:= \mathbb{E}[\boldsymbol{h}_f^k|\mathcal{F}_{k, 2}^{'}]$, if sample sets $\{\mathcal{D}_F, \mathcal{D}_G, \mathcal{B}_j(j=0,\ldots, J-1), \mathcal{D}_{F, j}(j=0,\ldots, J-1)\}$ are mutually independent, then, we have 
\begin{equation}   \label{eq58} \tag{58}
\bar{\boldsymbol{h}}_f^k = \nabla_{\boldsymbol{x}} f(\boldsymbol{x}^k, \boldsymbol{y}^{k, T}) - \nabla_{\boldsymbol{x}}\nabla_{\boldsymbol{y}}g(\boldsymbol{x}^k, \boldsymbol{y}^{k, T}) \boldsymbol{v}_J 
\end{equation}
where 
\begin{equation}  \label{eq59} \tag{59}
\boldsymbol{v}_J = (\mathbf{I} - \eta \nabla_{\boldsymbol{y}}^2 g(\boldsymbol{x}^k, \boldsymbol{y}^{k, T}))^J \boldsymbol{v}^{k, 0} + \eta \sum\limits_{t=0}^{J-1}(\mathbf{I} - \eta \nabla_{\boldsymbol{y}}^2 g(\boldsymbol{x}^k, \boldsymbol{y}^{k, T}))^{J-1-t}\nabla_{\boldsymbol{y}} f(\boldsymbol{x}^k, \boldsymbol{y}^{k, T}).
\end{equation}
\end{itemize}
\end{lemma}

\begin{IEEEproof}
The first two conclusions can be easily obtained by the formulas of $\boldsymbol{h}_f^k$ in (6) and (11).

Next, we prove the third conclusion. By the formula in (13), it is easy to know that 
\begin{equation*}   
\boldsymbol{v}^{k, J} = \prod_{i=0}^{J-1}(\mathbf{I} - \eta \nabla_{\boldsymbol{y}}^2 G(\boldsymbol{x}^k, \boldsymbol{y}^{k, T}; \mathcal{B}_i))\boldsymbol{v}^{k, 0} 
+ \eta \sum\limits_{t=0}^{J-1}\prod_{i=t+1}^{J-1} (\mathbf{I} - \eta \nabla_{\boldsymbol{y}}^2 G(\boldsymbol{x}^k, \boldsymbol{y}^{k, T}; \mathcal{B}_i))\nabla_{\boldsymbol{y}}F(\boldsymbol{x}^k, \boldsymbol{y}^{k, T}; \mathcal{D}_{F,t}).
\end{equation*}
Then, it is easy to prove that the third conclusion holds.
\end{IEEEproof}

In the following, we bound $\|\nabla \Phi(\boldsymbol{x}^k) - \bar{\boldsymbol{h}}_f^k\|$.

\begin{lemma}  \label{lemma:10}
Under the conditions of Lemma \ref{lemma:9}, let $\mathbf{H}_0 = \mathbf{0}$, $\boldsymbol{v}_0 = \boldsymbol{v}^{k, 0}$, and define
\begin{equation*}
\mathbf{H}^*  = \nabla_{\boldsymbol{x}}\nabla_{\boldsymbol{y}} g(\boldsymbol{x}^k, \boldsymbol{y}^*(\boldsymbol{x}^k))(\nabla_{\boldsymbol{y}}^2 g(\boldsymbol{x}^k, \boldsymbol{y}^*(\boldsymbol{x}^k)))^{-1}, \qquad \boldsymbol{v}^* = \nabla_{\boldsymbol{y}}^2 g(\boldsymbol{x}^k, \boldsymbol{y}^*(\boldsymbol{x}^k))^{-1}\nabla_{\boldsymbol{y}} f(\boldsymbol{x}^k, \boldsymbol{y}^*(\boldsymbol{x}^k))
\end{equation*}
where $\boldsymbol{v}^{k, 0}$ is the initial value to solve problem (12)(see line 13 in Algorithm 1), and $\boldsymbol{y}^*(\boldsymbol{x}^k)$ is the optimal solution to the LL problem of problem (2) with $\boldsymbol{x} = \boldsymbol{x}^k$. Then  
\begin{itemize}
\item For Stochastic BP, for $\bar{\boldsymbol{h}}_f^k$ in (\ref{eq55}), we have
\begin{align*}
\|\nabla \Phi(\boldsymbol{x}^k) - \bar{\boldsymbol{h}}_f^k\| & \le L\bigg(1 +\frac{L}{\mu}\bigg) \|\boldsymbol{y}^{k, T} - \boldsymbol{y}^*(\boldsymbol{x}^k)\| + M(1 - \eta \mu)^T \|\mathbf{H}^*\| \\
&\qquad + M\eta\bigg(\frac{L}{\mu}\rho + \tau\bigg) \sum\limits_{t=0}^{T-1}(1 - \eta \mu)^t \|\boldsymbol{y}^{k, T-1-t} - \boldsymbol{y}^*(\boldsymbol{x}^k)\|;
\end{align*}
\item For Stochastic NS, for $\bar{\boldsymbol{h}}_f^k$ in (\ref{eq57}), we have 
\begin{align*}
\|\nabla \Phi(\boldsymbol{x}^k) - \bar{\boldsymbol{h}}_f^k\| & \le L\bigg(1 +\frac{L}{\mu}\bigg) \|\boldsymbol{y}^{k, T} - \boldsymbol{y}^*(\boldsymbol{x}^k)\| + M(1 - \eta \mu)^J \|\mathbf{H}^*\| \\
& \qquad + M\eta\bigg(\frac{L}{\mu}\rho + \tau\bigg) \sum\limits_{t=0}^{J-1}(1 - \eta \mu)^t \|\boldsymbol{y}^{k, T} - \boldsymbol{y}^*(\boldsymbol{x}^k)\|;
\end{align*}
\item For SGD-based Estimation, for $\bar{\boldsymbol{h}}_f^k$ in (\ref{eq58}), we have 
\begin{align*}
\|\nabla \Phi(\boldsymbol{x}^k) - \bar{\boldsymbol{h}}_f^k\| & \le \bigg(L +\frac{M}{\mu}\tau \bigg) \|\boldsymbol{y}^{k, T} - \boldsymbol{y}^*(\boldsymbol{x}^k)\| + L(1 - \eta \mu)^J \| \boldsymbol{v}^{k, 0} - \boldsymbol{v}^* \| \\
&\qquad + L \eta\bigg(\frac{M}{\mu}\rho + L \bigg) \sum\limits_{t=0}^{J-1}(1 - \eta \mu)^t \|\boldsymbol{y}^{k, T} - \boldsymbol{y}^*(\boldsymbol{x}^k)\|;
\end{align*}
\end{itemize}
$\nabla \Phi(\boldsymbol{x}^k)$ is in (3) with $\boldsymbol{x} = \boldsymbol{x}^k$, $\mu$, $L$, $M$, $\tau$, $\rho$ are given in Assumptions 1, 2. 
\end{lemma}

\begin{IEEEproof}
Firstly, we prove the first conclusion. By the definition of $\nabla \Phi(\boldsymbol{x})$ in (3), it is easy to know that 
\begin{equation*}
\nabla \Phi(\boldsymbol{x}^k) = \nabla_{\boldsymbol{x}} f(\boldsymbol{x}^k, \boldsymbol{y}^*(\boldsymbol{x}^k)) - \mathbf{H}^* \nabla_{\boldsymbol{y}} f(\boldsymbol{x}^k, \boldsymbol{y}^*(\boldsymbol{x}^k)).
\end{equation*}

Then, for $\nabla \Phi(\boldsymbol{x}^k)$ and $\bar{\boldsymbol{h}}_f^k$ in (\ref{eq55}), we have 
\begin{align}   \label{eq60}    \nonumber
\|\nabla \Phi(\boldsymbol{x}^k) - \bar{\boldsymbol{h}}_f^k \| &\le \|\nabla_{\boldsymbol{x}} f(\boldsymbol{x}^k, \boldsymbol{y}^*(\boldsymbol{x}^k)) - \nabla_{\boldsymbol{x}} f(\boldsymbol{x}^k, \boldsymbol{y}^{k, T})\| + \|\mathbf{H}^* \nabla_{\boldsymbol{y}} f(\boldsymbol{x}^k, \boldsymbol{y}^*(\boldsymbol{x}^k)) - \mathbf{H}_T \nabla_{\boldsymbol{y}}f(\boldsymbol{x}^k, \boldsymbol{y}^{k, T})\| \\   \tag{60}
& \le L\|\boldsymbol{y}^{k, T} - \boldsymbol{y}^*(\boldsymbol{x}^k)\| + \|\mathbf{H}^*\|\|\nabla_{\boldsymbol{y}} f(\boldsymbol{x}^k, \boldsymbol{y}^*(\boldsymbol{x}^k)) - \nabla_{\boldsymbol{y}} f(\boldsymbol{x}^k, \boldsymbol{y}^{k, T})\|   \\  \nonumber
&\qquad + \|\mathbf{H}_T - \mathbf{H}^*\| \|\nabla_{\boldsymbol{y}} f(\boldsymbol{x}^k, \boldsymbol{y}^{k, T}))\|   \\ \nonumber
& \overset{(i)}{\le} \bigg(L + \frac{L^2}{\mu}\bigg) \|\boldsymbol{y}^{k, T} - \boldsymbol{y}^*(\boldsymbol{x}^k)\| + M\|\mathbf{H}_T - \mathbf{H}^*\|
\end{align}
where $(i)$ follows from Assumptions 1 and 2.

In the following, we bound the term $\|\mathbf{H}_T - \mathbf{H}^*\|$. By the definition of $\mathbf{H}_T$ in (\ref{eq56}), it is easy to know that 
\begin{align*}
\mathbf{H}_T &= \mathbf{H}_{T-1}(\mathbf{I} - \eta \nabla_{\boldsymbol{y}}^2 g(\boldsymbol{x}^k, \boldsymbol{y}^{k, T-1})) + \eta \nabla_{\boldsymbol{x}}\nabla_{\boldsymbol{y}}g(\boldsymbol{x}^k, \boldsymbol{y}^{k, T-1})  \\
& = (1- \eta)\mathbf{H}_{T-1} + \eta \bigg(\mathbf{H}_{T-1}\big(\mathbf{I} - \nabla_{\boldsymbol{y}}^2 g(\boldsymbol{x}^k, \boldsymbol{y}^{k, T-1})\big) + \nabla_{\boldsymbol{x}}\nabla_{\boldsymbol{y}} g(\boldsymbol{x}^k, \boldsymbol{y}^{k, T-1})\bigg).
\end{align*}

Then, 
\begin{align*}
\|\mathbf{H}_T - \mathbf{H}^*\| & \overset{(i)}{=} \left\|(1 - \eta) \mathbf{H}_{T-1} + \eta \bigg(\mathbf{H}_{T-1}\big(\mathbf{I} - \nabla_{\boldsymbol{y}}^2 g(\boldsymbol{x}^k, \boldsymbol{y}^{k, T-1})\big) + \nabla_{\boldsymbol{x}}\nabla_{\boldsymbol{y}} g(\boldsymbol{x}^k, \boldsymbol{y}^{k, T-1})\bigg) \right. \\
& \qquad \left. -(1 - \eta)\mathbf{H}^* - \eta\bigg( \mathbf{H}^* \big(\mathbf{I} - \nabla_{\boldsymbol{y}}^2 g(\boldsymbol{x}^k, \boldsymbol{y}^*(\boldsymbol{x}^k))\big) +\nabla_{\boldsymbol{x}}\nabla_{\boldsymbol{y}} g(\boldsymbol{x}^k, \boldsymbol{y}^*(\boldsymbol{x}^k))\bigg)   \right\| \\
& \le (1 - \eta) \|\mathbf{H}_{T-1} - \mathbf{H}^*\| +\eta \bigg\|\mathbf{H}_{T-1}\big(\mathbf{I} - \nabla_{\boldsymbol{y}}^2 g(\boldsymbol{x}^k, \boldsymbol{y}^{k, T-1})\big) - \mathbf{H}^*\big(\mathbf{I} - \nabla_{\boldsymbol{y}}^2 g(\boldsymbol{x}^k, \boldsymbol{y}^*(\boldsymbol{x}^k))\big)\bigg\|  \\
& \qquad + \eta \big\|\nabla_{\boldsymbol{x}}\nabla_{\boldsymbol{y}} g(\boldsymbol{x}^k, \boldsymbol{y}^{k, T-1}) - \nabla_{\boldsymbol{x}}\nabla_{\boldsymbol{y}} g(\boldsymbol{x}^k, \boldsymbol{y}^*(\boldsymbol{x}^k))\big\|  \\
& \le (1 - \eta) \|\mathbf{H}_{T-1} - \mathbf{H}^*\| +\eta \|\mathbf{H}_{T-1} - \mathbf{H}^*\| \|\mathbf{I} - \nabla_{\boldsymbol{y}}^2 g(\boldsymbol{x}^k, \boldsymbol{y}^{k, T-1})\| \\
& \qquad + \eta \|\mathbf{H}^*\| \big\|\nabla_{\boldsymbol{y}}^2 g(\boldsymbol{x}^{k}, \boldsymbol{y}^{k, T-1}) - \nabla_{\boldsymbol{y}}^2 g(\boldsymbol{x}^k, \boldsymbol{y}^*(\boldsymbol{x}^k))\big\| \\
& \qquad + \eta \|\nabla_{\boldsymbol{x}}\nabla_{\boldsymbol{y}} g(\boldsymbol{x}^k, \boldsymbol{y}^{k, T-1}) - \nabla_{\boldsymbol{x}}\nabla_{\boldsymbol{y}} g(\boldsymbol{x}^k, \boldsymbol{y}^*(\boldsymbol{x}^k))\| \\
& \overset{(ii)}{\le}(1 -\eta\mu )\|\mathbf{H}_{T-1} - \mathbf{H}^*\| + \eta\bigg(\frac{L}{\mu}\rho + \tau\bigg)\|\boldsymbol{y}^{k, T-1}- \boldsymbol{y}^*(\boldsymbol{x}^k)\|\\
& \overset{(iii)}{\le} (1 - \eta \mu)^T \|\mathbf{H}_0 - \mathbf{H}^*\| + \eta \bigg(\frac{L}{\mu}\rho + \tau \bigg)\sum\limits_{t=0}^{T-1}(1 - \eta \mu )^t \|\boldsymbol{y}^{k, T-1-t} - \boldsymbol{y}^*(\boldsymbol{x}^k)\|
\end{align*}
where $(i)$ uses the definition of $\mathbf{H}^*$, $(ii)$ follows from Assumptions 1 and 2, and $(iii)$ is obtained by using induction.

Combined with (\ref{eq60}), the first conclusion can be proved. Furthermore, the second conclusion can be proved similarly.

Next, we prove the third conclusion. By the definition of $\nabla \Phi(\boldsymbol{x})$ in (3), it is easy to know that 
\begin{equation*}
\nabla \Phi(\boldsymbol{x}^k) = \nabla_{\boldsymbol{x}}f(\boldsymbol{x}^k, \boldsymbol{y}^*(\boldsymbol{x}^k)) - \nabla_{\boldsymbol{x}}\nabla_{\boldsymbol{y}} g(\boldsymbol{x}^k, \boldsymbol{y}^*(\boldsymbol{x}^k))\boldsymbol{v}^*.
\end{equation*}

Then, for $\nabla \Phi(\boldsymbol{x}^k)$ and $\bar{\boldsymbol{h}}_f^k$ in (\ref{eq58}), we have 
\begin{align}   \label{eq61} \nonumber
\|\nabla \Phi(\boldsymbol{x}^k) - \bar{\boldsymbol{h}}_f^k\| &\le \|\nabla_{\boldsymbol{x}} f(\boldsymbol{x}^k, \boldsymbol{y}^*(\boldsymbol{x}^k)) - \nabla_{\boldsymbol{x}} f(\boldsymbol{x}^k, \boldsymbol{y}^{k, T})\| +\|\nabla_{\boldsymbol{x}}\nabla_{\boldsymbol{y}} g(\boldsymbol{x}^k, \boldsymbol{y}^*(\boldsymbol{x}^k))\boldsymbol{v}^* - \nabla_{\boldsymbol{x}}\nabla_{\boldsymbol{y}} g(\boldsymbol{x}^k, \boldsymbol{y}^{k, T})\boldsymbol{v}_J\| \\   \tag{61}
& \le L\|\boldsymbol{y}^{k, T}- \boldsymbol{y}^*(\boldsymbol{x}^k)\| + \|\boldsymbol{v}^*\| \|\nabla_{\boldsymbol{x}}\nabla_{\boldsymbol{y}} g(\boldsymbol{x}^k, \boldsymbol{y}^*(\boldsymbol{x}^k))  - \nabla_{\boldsymbol{x}}\nabla_{\boldsymbol{y}} g(\boldsymbol{x}^k, \boldsymbol{y}^{k, T})\|  \\    \nonumber
& \qquad + \|\boldsymbol{v}_J  - \boldsymbol{v}^*\| \|\nabla_{\boldsymbol{x}}\nabla_{\boldsymbol{y}} g(\boldsymbol{x}^k, \boldsymbol{y}^{k, T})\|  \\   \nonumber
&\overset{(i)}{\le}\bigg(L+ \frac{M}{\mu}\tau\bigg)\|\boldsymbol{y}^{k, T} - \boldsymbol{y}^*(\boldsymbol{x}^k)\| + L\|\boldsymbol{v}_J - \boldsymbol{v}^*\|  
\end{align}
where $(i)$ follows from Assumptions 1 and 2.

In the following, we bound the term $\|\boldsymbol{v}_J - \boldsymbol{v}^*\|$. By the definition of $\boldsymbol{v}_J$ in (\ref{eq59}), it is easy to know that 
\begin{align*}
\boldsymbol{v}_J & = (\mathbf{I} - \eta \nabla_{\boldsymbol{y}}^2 g(\boldsymbol{x}^k, \boldsymbol{y}^{k, T}))\boldsymbol{v}_{J-1} +\eta \nabla_{\boldsymbol{y}} f(\boldsymbol{x}^k, \boldsymbol{y}^{k, T})\\
& = (1 - \eta)\boldsymbol{v}_{J-1}+ \eta \bigg((\mathbf{I} - \nabla_{\boldsymbol{y}}^2 g(\boldsymbol{x}^k, \boldsymbol{y}^{k, T}))\boldsymbol{v}_{J-1} + \nabla_{\boldsymbol{y}} f(\boldsymbol{x}^k, \boldsymbol{y}^{k, T})\bigg).
\end{align*}

Then, 
\begin{align*}
\|\boldsymbol{v}_J - \boldsymbol{v}^*\| & \overset{(i)}{=} \left\|(1 - \eta) \boldsymbol{v}_{J-1} + \eta \bigg((\mathbf{I} - \nabla_{\boldsymbol{y}}^2 g(\boldsymbol{x}^{k}, \boldsymbol{y}^{k, T}))\boldsymbol{v}_{J-1} + \nabla_{\boldsymbol{y}} f(\boldsymbol{x}^k, \boldsymbol{y}^{k, T})\bigg)  \right. \\
& \qquad \left. - (1-\eta) \boldsymbol{v}^* - \eta \bigg(\big(\mathbf{I} - \nabla_{\boldsymbol{y}}^2 g(\boldsymbol{x}^k, \boldsymbol{y}^*(\boldsymbol{x}^k))\big) \boldsymbol{v}^* + \nabla_{\boldsymbol{y}} f(\boldsymbol{x}^k, \boldsymbol{y}^*(\boldsymbol{x}^k))\bigg)\right\| \\
& \le (1 - \eta)\|\boldsymbol{v}_{J-1} - \boldsymbol{v}^*\| + \eta \bigg\|\big(\mathbf{I} - \nabla_{\boldsymbol{y}}^2 g(\boldsymbol{x}^k, \boldsymbol{y}^{k, T})\big)\boldsymbol{v}_{J-1} - \big(\mathbf{I} - \nabla_{\boldsymbol{y}}^2 g(\boldsymbol{x}^k, \boldsymbol{y}^*(\boldsymbol{x}^k))\big)\boldsymbol{v}^*  \bigg\|  \\
& \qquad + \eta \big\|\nabla_{\boldsymbol{y}} f(\boldsymbol{x}^k, \boldsymbol{y}^{k, T}) - \nabla_{\boldsymbol{y}} f(\boldsymbol{x}^k, \boldsymbol{y}^*(\boldsymbol{x}^k))\big\|  \\
& \le (1 - \eta) \|\boldsymbol{v}_{J-1}- \boldsymbol{v}^*\| + \eta \|\mathbf{I} - \nabla_{\boldsymbol{y}}^2 g(\boldsymbol{x}^k, \boldsymbol{y}^{k, T})\| \|\boldsymbol{v}_{J-1}- \boldsymbol{v}^*\| \\
& \qquad + \eta \|\boldsymbol{v}^*\| \|\nabla_{\boldsymbol{y}}^2 g(\boldsymbol{x}^k, \boldsymbol{y}^*(\boldsymbol{x}^k)) - \nabla_{\boldsymbol{y}}^2 g(\boldsymbol{x}^k, \boldsymbol{y}^{k, T})\| +\eta \|\nabla_{\boldsymbol{y}} f(\boldsymbol{x}^k, \boldsymbol{y}^{k, T}) - \nabla_{\boldsymbol{y}}f(\boldsymbol{x}^k, \boldsymbol{y}^*(\boldsymbol{x}^k))\| \\
& \overset{(i)}{\le} (1 - \eta \mu)\|\boldsymbol{v}_{J-1} - \boldsymbol{v}^*\| + \eta \bigg(\frac{M}{\mu}\rho + L\bigg) \|\boldsymbol{y}^{k, T} - \boldsymbol{y}^*(\boldsymbol{x}^k)\|\\
& \overset{(ii)}{\le} (1 - \eta \mu)^J \|\boldsymbol{v}_0 - \boldsymbol{v}^*\| + \eta\bigg(\frac{M}{\mu}\rho + L\bigg) \sum\limits_{t=0}^{J-1} (1 - \eta \mu)^t \|\boldsymbol{y}^{k, T} - \boldsymbol{y}^*(\boldsymbol{x}^k)\|
\end{align*}
where $(i)$ uses the definition of $\boldsymbol{v}^*$, $(ii)$ follows from Assumptions 1 and 2, and (iii) is obtained by using induction.

Combined with (\ref{eq61}), the third conclusion can be proved.
\end{IEEEproof}

\section{Details of Experiments}\label{appendix:3}
\subsection{Synthetic Bilevel Optimization Problems}
For the bilevel problem in (15), we first perform experiments to compare the performance of Algorithm 1 under different hypergradient estimation methods. In the experiments, the batch size is set to be $5$, $T$ is set to be $5$, $J$ is chosen from $\{1, 10\}$, the learning rate $\alpha$ to update $\boldsymbol{x}$ is set to be $0.001$, the learning rate $\beta$ to update the LL variable is set to be $0.1$, and the $\eta$ used to estimate the hypergradient is set to be $0.1$.

Then, for problem (15), we perform experiments to compare the algorithms SSGD, BSA, TTSA, and stocBiO, and study the influence of iteration step $T$, iteration step $J$, batch size, and stepsize $\alpha$ on the convergence behavior of SSGD. In the experiments, for all the algorithms, $J$ is chosen from $\{3, 5, 20\}$, $\eta$ is chosen from $\{0.003, 0.01, 0.1\}$. Moreover, for SSGD and stocBiO, the batch size is set to be $5$, the iteration step $T$ is chosen from $\{5, 7, 10\}$, the learning rate $\alpha$ to update $\boldsymbol{x}$ is chosen from $\{0.001, 0.0001\}$, and the learning rate $\beta$ to update the LL variable is chosen from $\{0.1, 0.01, 0.001\}$.

For BSA, following \cite{7}, batch size is set to be $1$. For each outer iteration indexed by $k$, the learning rate $\alpha_k$ to update $\boldsymbol{x}^k$ is set to be $0.1/(1 + k)^{\frac{1}{2}}$ and the iteration step $T_k$ is set to be $\lceil (k+1)^{1/2} \rceil$. For every $t \in \{ 0, \ldots, T_k - 1\}$, the learning rate $\beta_t$ to update the LL variable is set to be $0.1/(t+2)$.

For TTSA, following \cite{7}, batch size is set to be $1$. For each outer iteration indexed by $k$, the learning rate $\alpha_k$ to update $\boldsymbol{x}^k$ is set to be $0.1/(1 + k)^{3/5}$ and the learning rate $\beta_k$ to update the LL variable is set to be $0.1/(1+k)^{2/5}$.

\subsection{Data Hyper-cleaning}

For MNIST and FashionMNIST datasets, we use 55000 images for training, 5000 images for validation, and 10000 images for testing. For CIFAR10 dataset, we use 45000 images for training, 5000 images for validation, and 10000 images for testing.

For all the algorithms, $\eta$ is set to be 0.001, $J$ is set to be 4. For SSGD and stocBiO, the batch size is set to be 10, the iteration step $T$ is set to be 5, the learning rate $\alpha$ to update $\boldsymbol{x}$ is set to be 0.001, and the learning rate $\beta$ to update the LL variable is chosen from $\{0.001, 0.004\}$. 

For BSA, following \cite{7}, batch size is set to be $1$. For each outer iteration indexed by $k$, the learning rate $\alpha_k$ to update $\boldsymbol{x}^k$ is set to be $d_{\alpha}/(1 + k)^{\frac{1}{2}}$ and the iteration step $T_k$ is set to be $\lceil (k+1)^{1/2} \rceil$. For every $t \in \{ 0, \ldots, T_k - 1\}$, the learning rate $\beta_t$ to update the LL variable is set to be $d_{\beta}/(t+2)$. In addition, both $d_{\alpha}$ and $d_{\beta}$ are chosen from $\{0.001, 0.01, 0.1\}$.

For TTSA, following \cite{7}, batch size is set to be $1$. For each outer iteration indexed by $k$, the learning rate $\alpha_k$ to update $\boldsymbol{x}^k$ is set to be $d_{\alpha}/(1 + k)^{3/5}$ and the learning rate $\beta_k$ to update the LL variable is set to be $d_{\beta}/(1+k)^{2/5}$. In addition, both $d_{\alpha}$ and $d_{\beta}$ are chosen from $\{0.001, 0.01, 0.1\}$.

\section{Additional Results and Experiments}\label{appendix:4}

\subsection{Additional Results and Experiments of Data Hyper-cleaning}

In the following, we first present some additional results for the experiments on the data hyper-cleaning in the main text.

\begin{figure*}[htb]
\centering
\subfloat[]{\includegraphics[scale = 0.2]{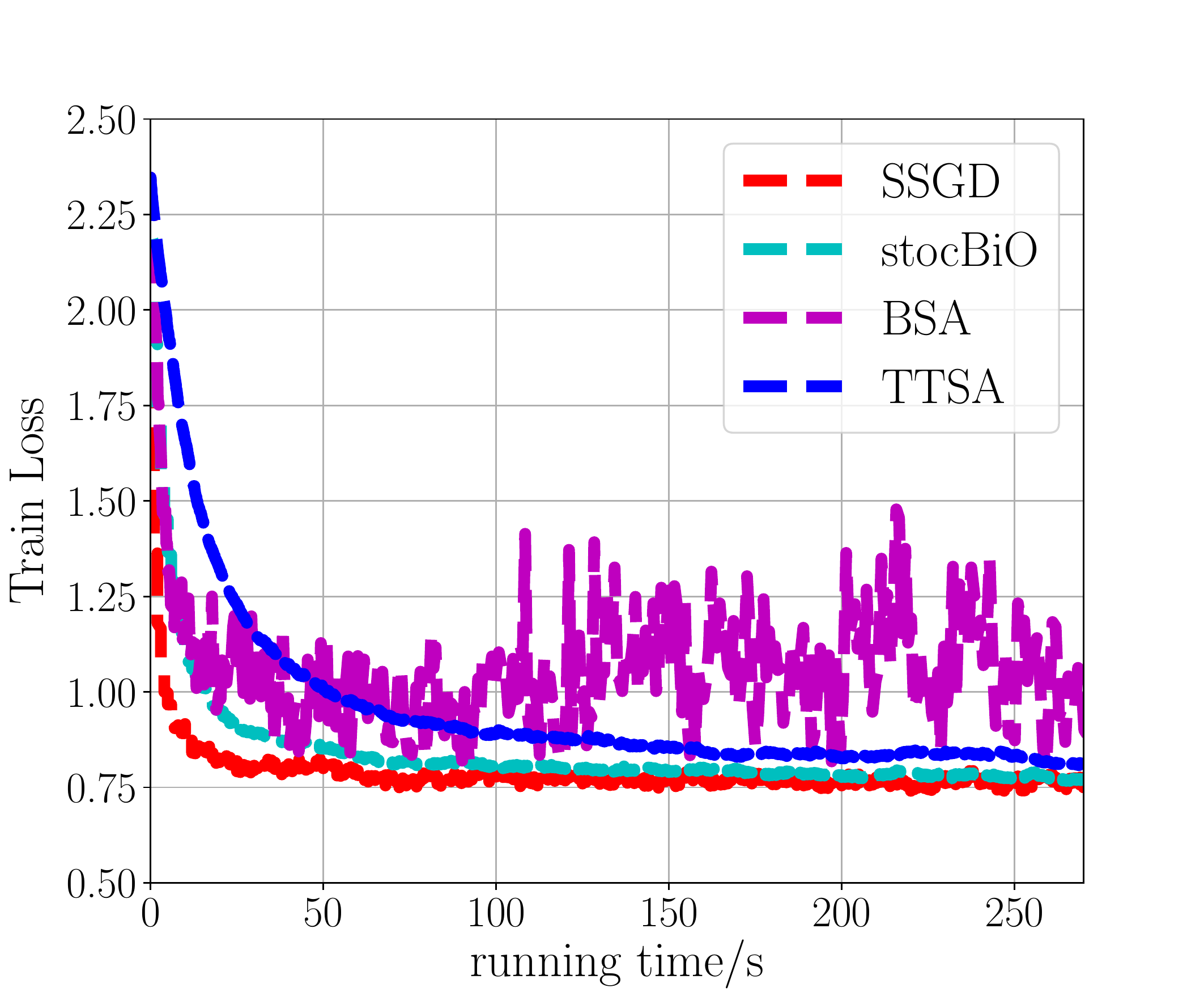}%
\label{fig_first_case}}
\hfil
\subfloat[]{\includegraphics[scale = 0.2]{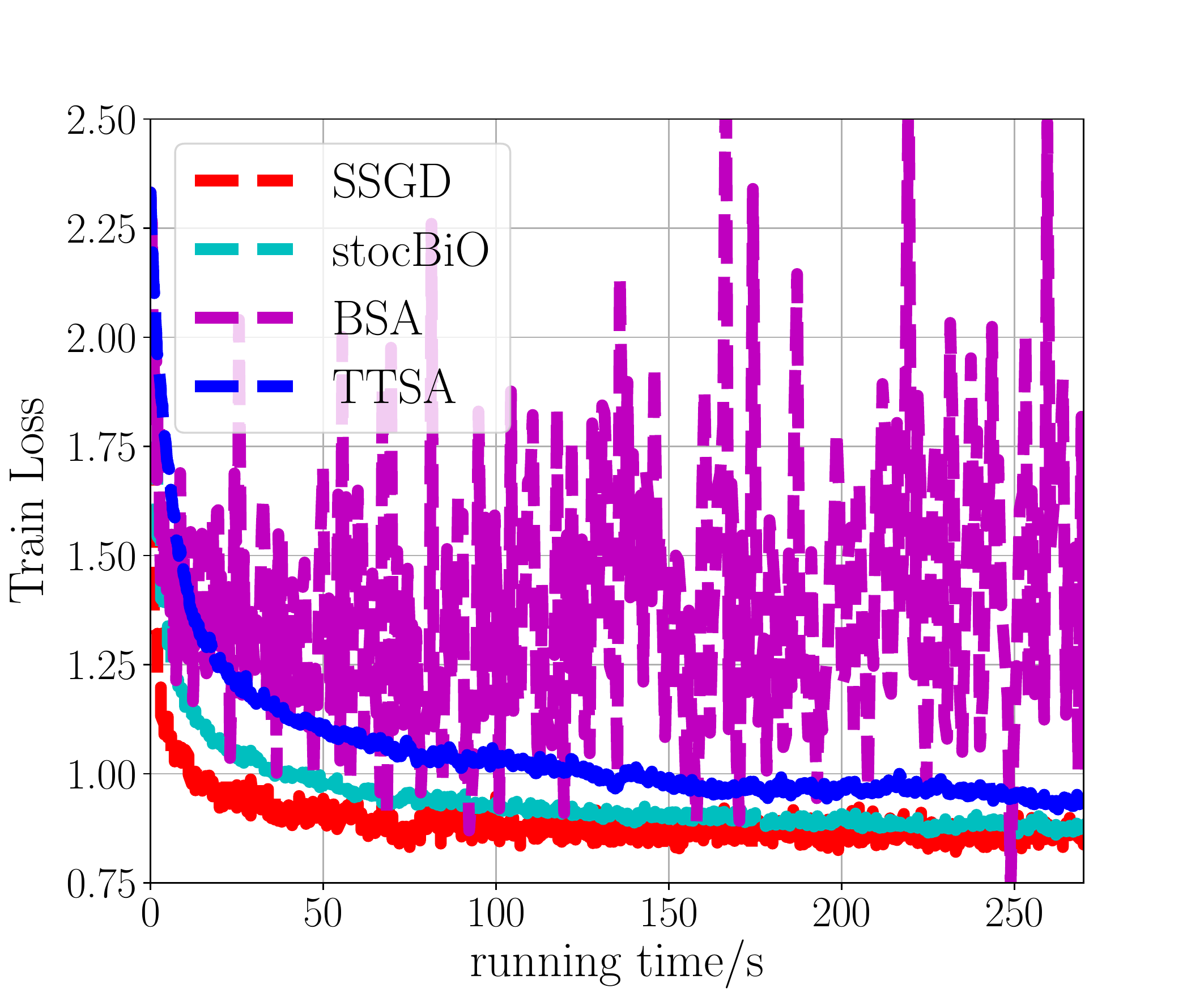}%
\label{fig_second_case}}
\hfil
\subfloat[]{\includegraphics[scale = 0.2]{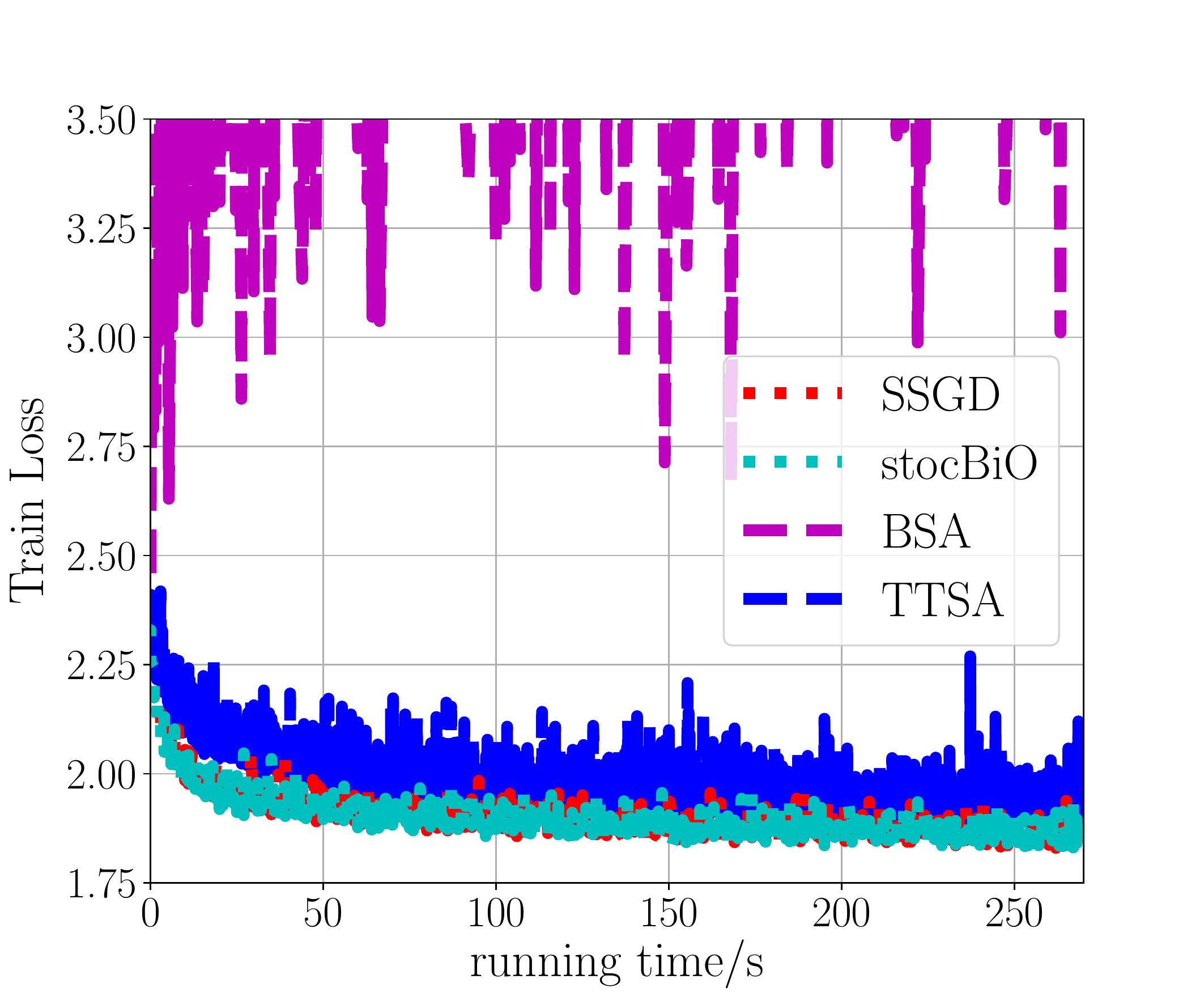}%
\label{fig_third_case}}
\caption{Comparison of the simple SGD type algorithms on MNIST, FashionMNIST, CIFAR10 datasets. (a) shows the resutls on the MNIST dataset. (b) shows the results on the FashionMNIST dataset. (c) shows the results on the CIFAR10 dataset.}
\label{fig:7}   
\end{figure*}

Fig. \ref{fig:7} shows the training loss v.s. running time curves for the experiments in the main text. We observe that SSGD and stocBiO converge faster than BSA and TTSA on all the datasets. Furthermore, SSGD has a slightly faster convergence speed than stocBiO on MNIST and FashionMNIST datasets.

\begin{figure}[htb]
\centering
\subfloat[]{\includegraphics[scale = 0.2]{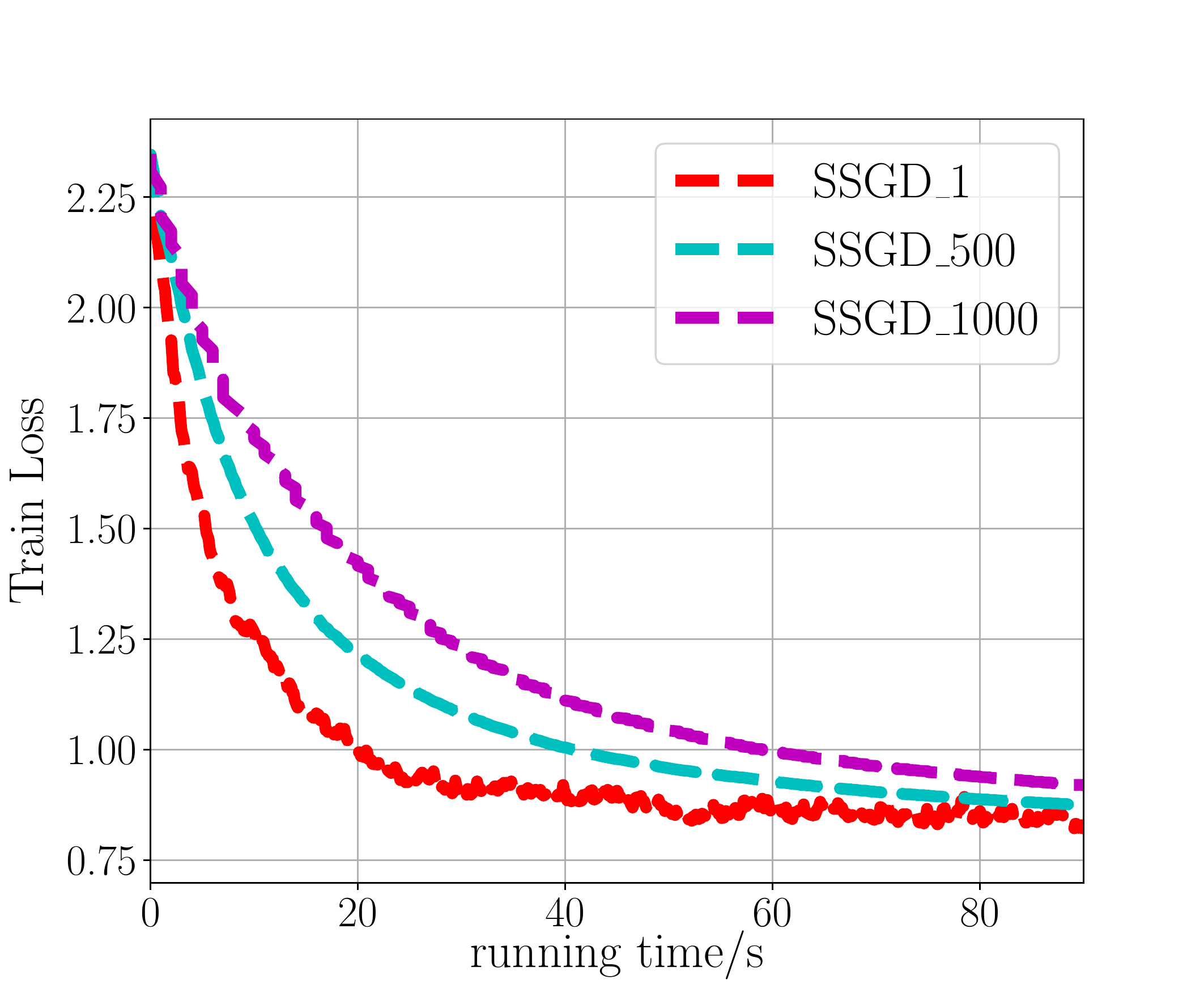}%
\label{fig_first_case}}
\hfil
\subfloat[]{\includegraphics[scale = 0.2]{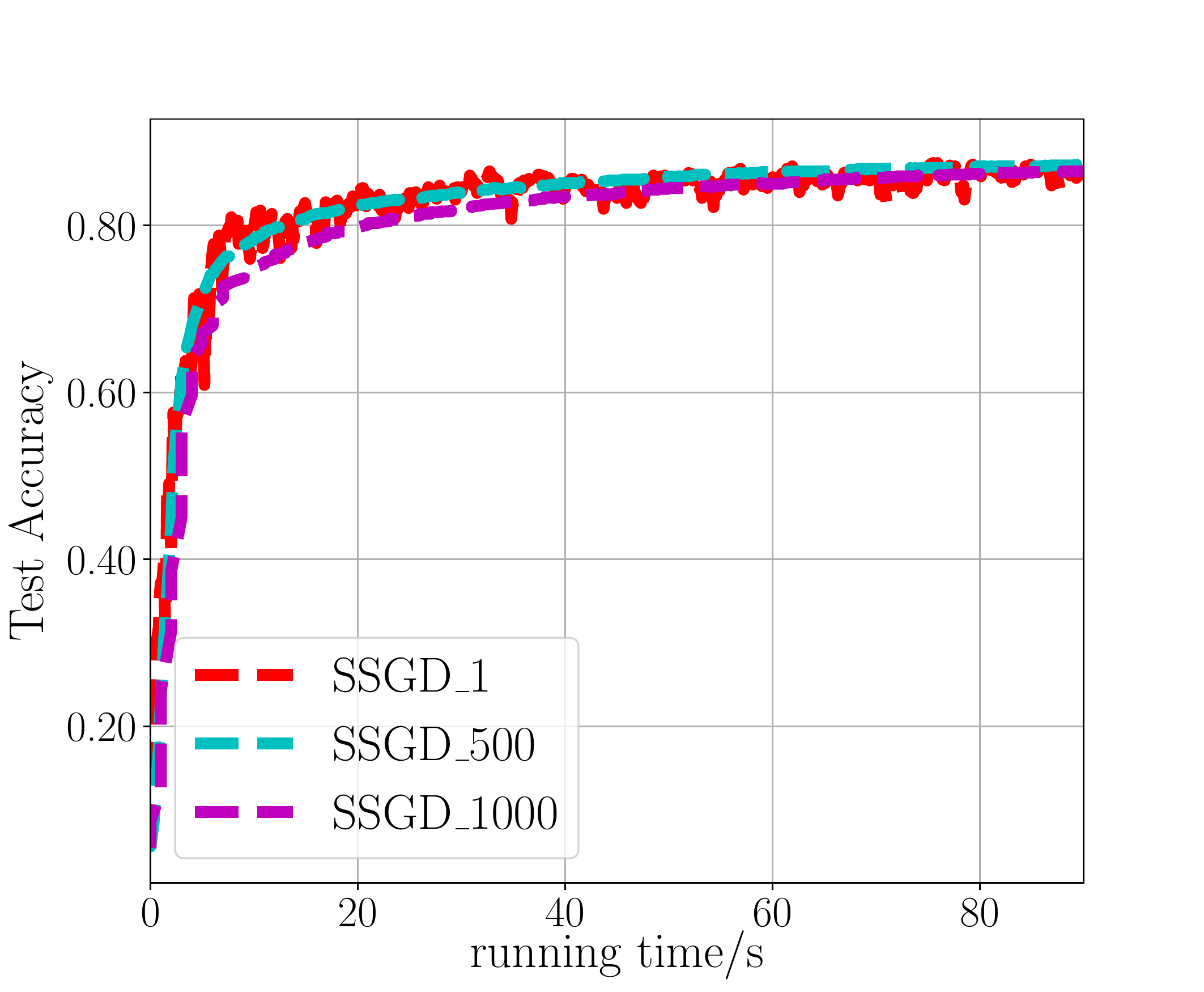}%
\label{fig_second_case}}
\caption{Comparing the performance of SSGD under different batch sizes.}
\label{fig:8}
\end{figure}

Then, we study the influence of batch size on SSGD. Fig. \ref{fig:8}\subref{fig_first_case} shows that SSGD achieves the fastest convergence speed when the batch size is set to be 1. Furthermore, Fig. \ref{fig:8}\subref{fig_second_case} shows that batch size has little influence on the best test accuracy achieved by SSGD.

\begin{figure}[htb]
\centering
\subfloat[]{\includegraphics[scale = 0.2]{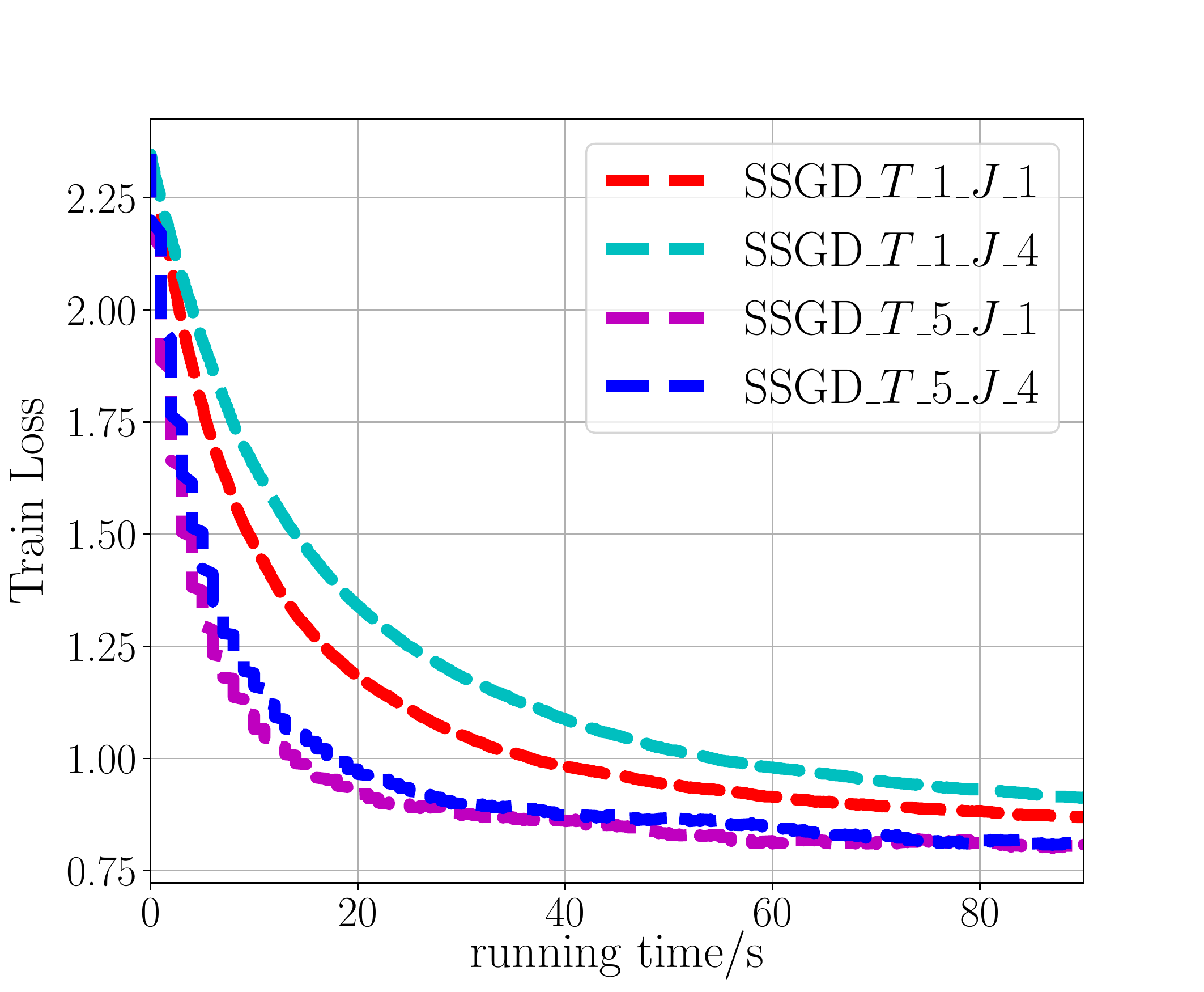}%
\label{fig_first_case}}
\hfil
\subfloat[]{\includegraphics[scale = 0.2]{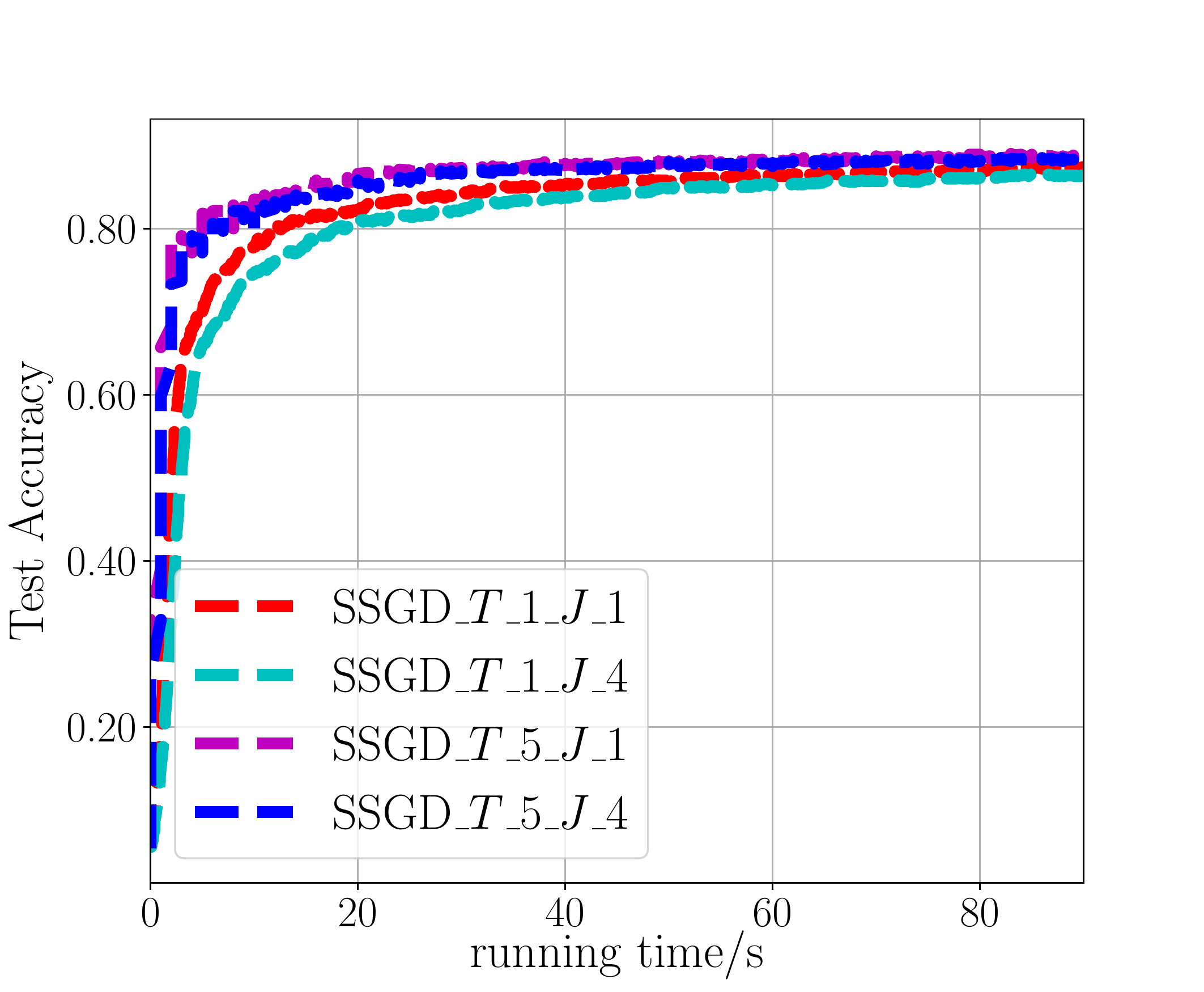}%
\label{fig_second_case}}
\caption{Comparing the performance of SSGD under different choices of $T$ and $J$.}
\label{fig:9}
\end{figure}

Next, we evaluate the performance of SSGD under different choices of $T$ and $J$. Fig. \ref{fig:9} shows that SSGD can obtain faster convergence speed and better test accuracy under the double loop case(i.e. $T>1$) than the single loop case(i.e., $T=1$).

\subsection{Additional Experiments on Meta Learning}

In the following, we evaluate the performance of SSGD on few-shot classification problem, a special case of meta learning, in which each task is a $N$-way classification, and $M$ training samples for each class of each task(ie., $N$-way $M$-shot).

The goal is to learn the parameter $\boldsymbol{x}$ of the feature extraction layer which can handle new tasks well. To be specific, we collect $m$ tasks $\{\mathcal{T}_i, i=1, \ldots, m\}$, and consider the few-shot classification problem whose objective function can be expressed as follows:

\begin{align*}  
&\min\limits_{\boldsymbol{x}} f\left( \boldsymbol{x}, \boldsymbol{y}^*(\boldsymbol{x})\right) := \frac{1}{m}\sum\limits_{i=1}^{m}\frac{1}{|\mathcal{D}_{\text{val}}^i|}\sum\limits_{\xi \in \mathcal{D}_{\text{val}}^i}L(\boldsymbol{x}, \boldsymbol{y}_i^*(\boldsymbol{x}); \xi) \\
&s.t. ~ \boldsymbol{y}^*(\boldsymbol{x}) :=  { \underset {\boldsymbol{y}} { \operatorname {arg\,min} } \, g(\boldsymbol{x}, \boldsymbol{y})} := \frac{1}{m}\sum\limits_{i=1}^m L_{\mathcal{D}_{\text{tr}}^i}(\boldsymbol{x}, \boldsymbol{y}_i; \xi)
\end{align*} 
where 
\begin{equation*}
L_{\mathcal{D}_{\text{tr}}^i}(\boldsymbol{x}, \boldsymbol{y}_i; \xi):= \frac{1}{|\mathcal{D}_{\text{tr}}^i|}\sum\limits_{\xi \in \mathcal{D}_{\text{tr}}^i}L(\boldsymbol{x}, \boldsymbol{y}_i; \xi) + 0.001 \|\boldsymbol{y}_i\|^2,
\end{equation*}
$L(\boldsymbol{x}, \boldsymbol{y}_i; \xi)$ denotes the cross-entropy function over sample $\xi$ for task $\mathcal{T}_i$, $\boldsymbol{y}:= (\boldsymbol{y}_1, \ldots, \boldsymbol{y}_m)$, $\boldsymbol{x}$ is the parameters of the feature extraction layer shared by all the tasks $\{\mathcal{T}_i, i=1, \ldots, m\}$, $\boldsymbol{y}_i$ is the parameters of the linear classifier layer for task $\mathcal{T}_i$, and $\mathcal{D}_{\text{tr}}^i$  and $\mathcal{D}_{\text{val}}^i$ are the training set and validation set for task $\mathcal{T}_i$, respectively. 

During model training, building on the embedded feature obtained from the feature extraction layer with parameter $\boldsymbol{x}$, each task $\mathcal{T}_i$ searches the parameters of the classifier which minimizes its loss over the training set $\mathcal{D}_{\text{tr}}^i$, i.e., $\boldsymbol{y}_i^*(\boldsymbol{x})$. Then, based on the minimizers $\boldsymbol{y}_i^*(\boldsymbol{x})$, $i=1, \ldots, m$, we update the parameters $\boldsymbol{x}$ of the feature extraction layer.

In the experiments, we compare our proposed algorithm SSGD with BSA, TTSA, and stocBiO on two datasets: Omniglot\cite{30} and FC100\cite{31}. Omniglot contains examples of 1623 different handwritten characters from 50 alphabets. FC100 is a dataset derived from CIFAR100, and contains 100 classes. 

\begin{table}[htb]   
\caption{Comparison of test accuracy achieved by the simple SGD type algorithms on few-shot classification tasks. Acc.(\%) denotes the test accuracy}
\label{table3}
\centering
\begin{tabular}{ccccc}\hline
\multicolumn{1}{c}{\multirow{2}{*}{Method}} & \multicolumn{2}{c}{Omniglot} & \multicolumn{2}{c}{FC100}  \\ 
 &5-way 1-shot & 5-way 5-shot &5-way 1-shot & 5-way 5-shot \\  \hline
stocBiO & \textbf{96} & \textbf{98.75} & 60& 72 \\
BSA & 88 & 89& 60& 70\\ 
TTSA & 72 & 74& 50 & 55 \\ 
SSGD & \textbf{96} & \textbf{98.75} & \textbf{61} & \textbf{74} \\ \hline
\end{tabular}
\end{table}

Table \ref{table3} shows the experimental results over 500s. We observe that SSGD achieves the highest test accuracy on all the datasets, which shows that SSGD can obtain competitive results on the few-shot classification problems.

\subsection{Experimental details of meta-learning}

In the experiments, following \cite{13}, we use a 4-layer convolutional neural network as the feature extraction layer. Furthermore, for Omniglot dataset, we use 1200 characters for training, and 423 characters for testing; for FC100 dataset, we use 60 classes for training, 20 classes for validation, and 20 classes for testing. 

For all the algorithms, $\eta$ is set to be 0.01, $J$ is set to be 3. For SSGD and stocBiO, the iteration step $T$ is set to be 3, the learning rate $\alpha$ to update $\boldsymbol{x}$ is set to be 0.03, and the learning rate $\beta$ to update the LL variable is set to be 0.03. 

For BSA, we follow the setting in \cite{7}. For each outer iteration indexed by $k$, the learning rate $\alpha_k$ to update $\boldsymbol{x}^k$ is set to be $d_{\alpha}/(1 + k)^{\frac{1}{2}}$ and the iteration step $T_k$ is set to be $\lceil (k+1)^{1/2} \rceil$. For every $t \in \{ 0, \ldots, T_k - 1\}$, the learning rate $\beta_t$ to update the LL variable is set to be $d_{\beta}/(t+2)$. In addition, both $d_{\alpha}$ and $d_{\beta}$ are chosen from $\{0.001, 0.01, 0.1\}$.

For TTSA, we following the setting in \cite{7}. For each outer iteration indexed by $k$, the learning rate $\alpha_k$ to update $\boldsymbol{x}^k$ is set to be $d_{\alpha}/(1 + k)^{3/5}$ and the learning rate $\beta_k$ to update the LL variable is set to be $d_{\beta}/(1+k)^{2/5}$. In addition, both $d_{\alpha}$ and $d_{\beta}$ are chosen from $\{0.001, 0.01, 0.1\}$.
\bibliographystyle{IEEEtran}
\bibliography{IEEEabrv, mybible}